\documentclass[11pt]{article}

\usepackage{a4wide,latexsym,graphicx,epstopdf,amsmath,varioref}
\DeclareGraphicsExtensions{.eps,.jpg}

\date{19 September 2017}

\newcommand{\numsection}[1]{\section{#1}\setcounter{equation}{0}}
\newcommand{\appnumsection}[1]{\section*{#1}\setcounter{equation}{0}
  \renewcommand{\theequation}{A.\arabic{equation}}
  \renewcommand{\thetheorem}{A.\arabic{theorem}}
  \renewcommand{\thetable}{A.\arabic{table}}
  \renewcommand{\thefigure}{A.\arabic{figure}}
  \renewcommand{\thesection}{A} }

\renewcommand{\theequation}{\arabic{section}.\arabic{equation}}
\renewcommand{\thetable}{\arabic{section}.\arabic{table}}
\renewcommand{\thefigure}{\arabic{section}.\arabic{figure}}

\newcommand{\beqn}[1]{\begin{equation}\label{#1}}
\newcommand{\eeqn}{\end{equation}}

\newcommand{\mat}[2]{\left(\begin{array}{#1}#2\end{array}\right)}

\newcommand{\req}[1]{(\ref{#1})}

\newcommand{\bqa}[1]{\begin{quote} \label{#1} \small \begin{verbatim}}
\newcommand{\bq}{\begin{quote} \small \begin{verbatim}}
\newcommand{\eq}{\end{quote}}

\newcommand{\ms}{\;\;\;\;}

\newcommand{\tim}[1]{\;\; \mbox{#1} \;\;}

\newcommand{\bctable}[1]{\begin{table}[htbp]
                         \begin{center}
                         \begin{tabular}{#1} }
\newcommand{\ectable}[1]{\end{tabular}
                         \caption{#1}
                         \end{center}
                         \end{table}}

\newtheorem{theorem}{Theorem}[section]
\newtheorem{lemma}[theorem]{Lemma}
\newtheorem{proposition}[theorem]{Proposition}
\newtheorem{corollary}[theorem]{Corollary}

\newtheorem{example}{Example}[section]

\newlength{\thmw}
\newcommand{\blem}[2]{\vspace{\baselineskip}\noindent\hbox{%
  \lower\fboxrule\hbox{\vbox{\hrule\hbox{\vrule \kern-\fboxrule \vbox{%
  \vspace{\fboxsep} \noindent\hspace{2\fboxsep}\parbox{\thmw}{
  \begin{lemma}{\rm #1}\end{lemma}\vspace{-\lastskip}}\hspace{\fboxsep}}%
  \kern-\fboxrule \vrule }}}}\newpage \hbox{%
  \lower\fboxrule\hbox{\vbox{\hbox{\vrule \kern-\fboxrule \vbox{%
  \noindent\hspace{2\fboxsep}\parbox{\thmw}{\rm #2}\hspace{\fboxsep}
  \vspace{4\fboxsep}}\kern-\fboxrule \vrule }\hrule }}}\vspace{\baselineskip}
}

\newcommand{\lblem}[3]{\vspace{\baselineskip}\noindent\hbox{%
  \lower\fboxrule\hbox{\vbox{\hrule\hbox{\vrule \kern-\fboxrule \vbox{%
  \vspace{\fboxsep} \noindent\hspace{2\fboxsep}\parbox{\thmw}{
  \begin{lemma}\label{#1}{\rm #2}\end{lemma}\vspace{-\lastskip}}
  \hspace{\fboxsep}}\kern-\fboxrule \vrule }}}}\newpage \hbox{%
  \lower\fboxrule\hbox{\vbox{\hbox{\vrule \kern-\fboxrule \vbox{%
  \noindent\hspace{2\fboxsep}\parbox{\thmw}{\rm #3}\hspace{\fboxsep}
  \vspace{4\fboxsep}}\kern-\fboxrule \vrule }\hrule }}}\vspace{\baselineskip}
}

\newcommand{\bthm}[2]{\vspace{\baselineskip}\noindent\hbox{%
  \lower\fboxrule\hbox{\vbox{\hrule\hbox{\vrule \kern-\fboxrule \vbox{%
  \vspace{\fboxsep} \noindent\hspace{2\fboxsep}\parbox{\thmw}{
  \begin{theorem}{\rm #1}\end{theorem}\vspace{-\lastskip}}\hspace{\fboxsep}}%
  \kern-\fboxrule \vrule }}}}\newpage \hbox{%
  \lower\fboxrule\hbox{\vbox{\hbox{\vrule \kern-\fboxrule \vbox{%
  \noindent\hspace{2\fboxsep}\parbox{\thmw}{\rm #2}\hspace{\fboxsep}
  \vspace{4\fboxsep}}\kern-\fboxrule \vrule }\hrule }}}\vspace{\baselineskip}
}

\newcommand{\lbthm}[3]{\vspace{\baselineskip}\noindent\hbox{%
  \lower\fboxrule\hbox{\vbox{\hrule\hbox{\vrule \kern-\fboxrule \vbox{%
  \vspace{\fboxsep} \noindent\hspace{2\fboxsep}\parbox{\thmw}{
  \begin{theorem}\label{#1}{\rm #2}\end{theorem}\vspace{-\lastskip}}
  \hspace{\fboxsep}}\kern-\fboxrule \vrule }}}}\newpage \hbox{%
  \lower\fboxrule\hbox{\vbox{\hbox{\vrule \kern-\fboxrule \vbox{%
  \noindent\hspace{2\fboxsep}\parbox{\thmw}{\rm #3}\hspace{\fboxsep}
  \vspace{4\fboxsep}}\kern-\fboxrule \vrule }\hrule }}}\vspace{\baselineskip}
}

\newcommand{\bcor}[2]{\vspace{\baselineskip}\noindent\hbox{%
  \lower\fboxrule\hbox{\vbox{\hrule\hbox{\vrule \kern-\fboxrule \vbox{%
  \vspace{\fboxsep} \noindent\hspace{2\fboxsep}\parbox{\thmw}{
  \begin{corollary}{\rm #1}\end{corollary}\vspace{-\lastskip}}\hspace{\fboxsep}}%
  \kern-\fboxrule \vrule }}}}\newpage \hbox{%
  \lower\fboxrule\hbox{\vbox{\hbox{\vrule \kern-\fboxrule \vbox{%
  \noindent\hspace{2\fboxsep}\parbox{\thmw}{\rm #2}\hspace{\fboxsep}
  \vspace{4\fboxsep}}\kern-\fboxrule \vrule }\hrule }}}\vspace{\baselineskip}
}

\newcommand{\lbcor}[3]{\vspace{\baselineskip}\noindent\hbox{%
  \lower\fboxrule\hbox{\vbox{\hrule\hbox{\vrule \kern-\fboxrule \vbox{%
  \vspace{\fboxsep} \noindent\hspace{2\fboxsep}\parbox{\thmw}{
  \begin{corollary}\label{#1}{\rm #2}\end{corollary}\vspace{-\lastskip}}
  \hspace{\fboxsep}}\kern-\fboxrule \vrule }}}}\newpage \hbox{%
  \lower\fboxrule\hbox{\vbox{\hbox{\vrule \kern-\fboxrule \vbox{%
  \noindent\hspace{2\fboxsep}\parbox{\thmw}{\rm #3}\hspace{\fboxsep}
  \vspace{4\fboxsep}}\kern-\fboxrule \vrule }\hrule }}}\vspace{\baselineskip}
}

\newcommand{\bprop}[2]{\vspace{\baselineskip}\noindent\hbox{%
  \lower\fboxrule\hbox{\vbox{\hrule\hbox{\vrule \kern-\fboxrule \vbox{%
  \vspace{\fboxsep} \noindent\hspace{2\fboxsep}\parbox{\thmw}{
 \begin{proposition}{\rm #1}\end{proposition}\vspace{-\lastskip}}\hspace{\fboxsep}}%
  \kern-\fboxrule \vrule }}}}\newpage \hbox{%
  \lower\fboxrule\hbox{\vbox{\hbox{\vrule \kern-\fboxrule \vbox{%
  \noindent\hspace{2\fboxsep}\parbox{\thmw}{\rm #2}\hspace{\fboxsep}
  \vspace{4\fboxsep}}\kern-\fboxrule \vrule }\hrule }}}\vspace{\baselineskip}
}

\newcommand{\lbprop}[3]{\vspace{\baselineskip}\noindent\hbox{%
  \lower\fboxrule\hbox{\vbox{\hrule\hbox{\vrule \kern-\fboxrule \vbox{%
  \vspace{\fboxsep} \noindent\hspace{2\fboxsep}\parbox{\thmw}{
  \begin{proposition}\label{#1}{\rm #2}\end{proposition}\vspace{-\lastskip}}
  \hspace{\fboxsep}}\kern-\fboxrule \vrule }}}}\newpage \hbox{%
  \lower\fboxrule\hbox{\vbox{\hbox{\vrule \kern-\fboxrule \vbox{%
  \noindent\hspace{2\fboxsep}\parbox{\thmw}{\rm #3}\hspace{\fboxsep}
  \vspace{4\fboxsep}}\kern-\fboxrule \vrule }\hrule }}}\vspace{\baselineskip}
}

\newcounter{algo}[section]
\renewcommand{\thealgo}{\thesection.\arabic{algo}}

\newcommand{\algo}[3]{\refstepcounter{algo}
\begin{center}\begin{figure}[htbp]
\framebox[\textwidth]{
\parbox{0.95\textwidth} {\vspace{\topsep}
{\bf Algorithm \thealgo : #2}\label{#1}\\
\vspace*{-\topsep} \mbox{ }\\
{#3} \vspace{\topsep} }}
\end{figure}\end{center}}

\newcommand{\balgo}[4]{\refstepcounter{algo}
\begin{center}\begin{figure}[htbp]
\vspace{\baselineskip}\noindent\hbox{%
  \lower\fboxrule\hbox{\vbox{\hrule\hbox{\vrule \kern-\fboxrule \vbox{%
  \vspace{\topsep} \noindent\hspace{2\fboxsep}\parbox{\thmw}{\vspace{0.5\topsep}
  {\bf Algorithm \thealgo : #2}\label{#1}\\
  \vspace*{-\topsep} \mbox{ }\\
  {\rm #3}\vspace{-\lastskip}}
  \hspace{\fboxsep}}\kern-\fboxrule \vrule }}}}
\end{figure}\end{center}
\newpage \hbox{%
  \lower\fboxrule\hbox{\vbox{\hbox{\vrule \kern-\fboxrule \vbox{%
  \noindent\hspace{2\fboxsep}\parbox{\thmw}{\rm #4}\hspace{\fboxsep}
  \vspace{\fboxsep}}\kern-\fboxrule \vrule }\hrule }}}\vspace{\baselineskip}
}

 \newcommand{\calH}{{\cal H}}

\newcommand{\calM}{{\cal M}}

\renewcommand{\Re}{\hbox{I\hskip -2pt R}}

\newcommand{\Na}{\hbox{I\hskip -1.8pt N}}
  
\newcommand{\sfrac}[2]{{\scriptstyle \frac{#1}{#2}}}
\newcommand{\half}{\sfrac{1}{2}}

\newcommand{\citet}[1]{\citeauthor{#1}, \citeyear{#1}}

\newcommand{\comment}[1]{}

\newcommand{\eqdef}{\stackrel{\rm def}{=}}

\newcommand{\diag}{{\rm diag}}

\newcommand{\ii}[1]{\{1, \ldots, #1 \}}

\title{On the use of the saddle formulation in weakly-constrained 4D-VAR data assimilation}

\author{S. Gratton
  \thanks{Universit\'e de Toulouse, INP, IRIT, Toulouse, France. Email: serge.gratton@enseeiht.fr },
  S. G\"{u}rol
  \thanks{CERFACS,  Toulouse, France. Email: selime.gurol@cerfacs.fr},
  E. Simon
  \thanks{Universit\'e de Toulouse, INP, IRIT, Toulouse, France. Email: ehouarn.simon@enseeiht.fr},
  and Ph. L. Toint
  \thanks{NAXYS, University of Namur, Namur, Belgium. Email: philippe.toint@unamur.be}}

\begin{document}

\maketitle
\begin{abstract}
This paper discusses the practical use of the saddle variational formulation
for the weakly-constrained 4D-VAR method in data assimilation.  It is shown
that the method, in its original form, may produce erratic results or diverge
because of the inherent lack of monotonicity of the produced objective
function values. Convergent, variationaly coherent variants of the algorithm
are then proposed whose practical performance is compared to that of other
formulations.  This comparison is conducted on two data assimilation instances
(Burgers equation and the Quasi-Geostrophic model), using two different
assumptions on parallel computing environment. Because these
variants essentially retain the parallelization advantages of the original
proposal, they often --- but not always --- perform best, even for moderate
numbers of computing processes.
\end{abstract}

\vspace*{1cm}

\section{Introduction}

Data assimilation has long been an integral and important part of weather
forecasting, as new (and often incomplete) meteorological observations are
integrated in the ongoing process of predicting the weather for the next few
days \cite{Bout99}. The question here is that of using the data to
determine a ``best'' current state of the weather system from which elaborate
models may then be evolved in time, providing the desired predictions.  Among
the possible techniques for this task, variational methods have been applied
extensively, typically weighting the use of a priori knowledge (often
materialized by the specification of a background state $x_b$) with the
quality of the fit to the observations. This is the case, in particular, for the
well-known 4D-Var formulation \cite{LeDiTala86,Cour97}. In recent years, it has also
become necessary to take possible model errors into account, thus weighting a
priori knowledge, data fitting and model error reduction, an approach which
leads to the ``weakly-constrained 4D-Var'' formulation of the relevant data
assimilation problem \cite{Zupa97,VidaPiacleDi04,Trem06, Trem07}. In one of the formulations, 
the total time horizon (assimilation window) considered
is split into a number ($N_{sw}$) of time sub-windows, and the problem can be
be written as 
\beqn{wk4DVar}
\min_{x \in\Re^s} J(x) \eqdef
     \half \|x^{(0)} - x_b\|_{B^{-1}}^2
   + \half \sum_{j=0}^{N_{sw}}
     \left\| {\cal H}_j \big( x^{(j)}\big) - y_j \right\|_{{R}_j^{-1}}^2
   + \half \sum_{j=1}^{N_{sw}}
     \|   x^{(j)}- {\cal M}_j(x^{(j-1)}) \|_{{Q}_j^{-1}}^2
\eeqn
where
\begin{itemize}
\item
${x} = ( x^{(0)}, x^{(1)}, \ldots, x^{(N_{sw})})^T \in\Re^s$ ($s=n(N_{sw}+1)$)
is the control variable (with $x^{(j)} = x(t_j)$),
\item  $x_b$  $\in\Re^n$ is the background given at the initial time ($t_0$).
\item  $y_j\in\Re^{m_j}$  is the observation vector over a given time interval
\item ${\cal H}_j$ maps the state vector $x_j$ from model space to observation space
\item  ${\cal M}_j$ represents an integration of the numerical model from time
  $t_{j-1}$ to  $t_j$
\item $B$, $R_j$ and $Q_j$ are the positive-definite covariance matrices for
  background, observation and model error, respectively.
\end{itemize}
The incorporation of possible model errors is achieved by the presence
of the third term in the objective function.

As it is the case for the standard 4D-Var (consisting of the first two terms
in \req{wk4DVar}), the general unconstrained nonlinear least-squares problem
is solved by applying the Gauss-Newton algorithm
\cite{CourThepHoll94,GratLawlNich07}, which iteratively proceeds by
linearizing $\calH$ and $\calM$ at the current iterate and then, often very
approximately, minimizing the resulting quadratic function. A practically
crucial question is then how to approximately perform this minimization.  As
this is a main theme of the present paper, we immediately stress here that
this aim (approximate minimization) is often very different from approximate
gradient/residual reduction (although they coincide if the minimization is
exact).  Key factors for selecting a subproblem solver include the choice of a
quadratic model (re)formulation (known as the variational formulation), the
choice of a preconditioner, the parallelization potential of the resulting
algorithm.

Three formulations are available (state, forcing and saddle) and are detailed
in the next section. The ``saddle formulation'', discussed in
\cite{FishTremAuviTanPoli11,FishGuro17,FishGratGuroVassTrem17}, has recently
attracted interest of practioners because of its appealing potential for
parallel computing while still allowing a wide choice of
preconditioners. However it is fair to say that numerical experience with this
approach remains scarse so far, prompting for a more detailed assessment.

The purpose of the present paper is to propose such an assessment. It will be
shown that, left to its own devices, the original algorithm for the saddle
formulation may produce erratic results or diverge altogether. To circumvent
this problem, a more elaborate variant of the same approach will be proposed
which, at the same time, guarantees convergence of the overall Gauss-Newton
algorithm and essentially retains the excellent parallelization features of
the original method.  The numerical performance will be illustrated and
compared to that of state and forcing formulations on an assimilation example
on the nonlinear Burgers equation and on the two-layers Quasi-Geostrophic (QG)
atmospheric model\footnote{Quasi-geostrophic motion means that, in the
  horizontal direction of the atmospheric flow, the Coriolis force caused by
  the rotation of the Earth, and the pressure gradient force are in
  approximate balance.} provided within OOPS by the European Centre for Medium
Range Weather Forecasts (ECMWF). This model is widely used in theoretical
atmospheric studies, since it is simple enough for numerical experimentations
and yet adequately captures the most relevant large-scale dynamics in the
atmosphere. For more details on the QG model, see
\cite{FishGratGuroVassTrem17,FishGuro17}. This comparison will demonstrate
that its parallel computing features help to explain why the new variant often
outperforms other approaches. Influence of the choice of preconditioner,
detailed operator cost and data organization will also be discussed.

The paper is organized as follows. Section~\ref{formulation-s} provides the
necessary background and notations for the three variational formulations
mentioned above, including some of the associated preconditioning issues and a
discussion of the parallelization bottlenecks. The potentially problematic
behaviour of the original saddle method is then discussed and illustrated in
Section~\ref{saddle-original-s}, and the new algorithmic variants described in
Section~\ref{saddle-new-s}.  A comparison of these new variants on the Burgers
and QG examples is then proposed in Section~\ref{numerics-s} and
\ref{numerics-qg-s}, respectively.  The special case where the inverse of the
correlation matrices is not available is briefly considered in
Section~\ref{approxDinv-s}, while conclusions and perspectives are outlined in
Section~\ref{concl-s}.

{\bf Notations.} The Euclidean product of the vectors $x$ and $y$ is denoted
by $y^Tx$ and the induced Euclidean norm of $x$ by $\|x\|$.

\numsection{Problem formulations and preconditioning}\label{formulation-s}

As indicated above, the formulation of the subproblem at each iteration of
the Gauss-Newton algorithm for solving \req{wk4DVar} is crucial for good
computational performance.  If the operators ${M}_j$ are the linearized $\calM_j$ and ${H}_j$
are the linearized $\calH_j$, this subproblem can be expressed in terms of the increment
$\delta x$ as
\beqn{QuadProb-st}
\min_{ {\delta x} \in \Re^s} \; q_{\rm st}(\delta x) \eqdef
     \half \| L\, {\delta x} - b \|_{D^{-1}}^2
   + \half \| H\, {\delta x} - d \|_{R^{-1}}^2
\eeqn
where
\beqn{Ldef}
L = \mat{ccccc}{
               I  &    &  & &\\
            -M_1  &  I &  & & \\
                  & -M_2  & I & &\\
                  &       & \ddots  & \ddots &\\       
                  &       &         & -M_{N_{sw}} & I\\
},
\eeqn
for suitable ``misfit'' vectors
\[
d = ( d_0^T, d_1^T, \ldots, d_{N_{sw}}^T)^T
\tim{ and }
b = (   b_0^T, c_1^T, \ldots, c_{N_{sw}}^T)^T,
\]
and where
\beqn{H-D-def}
H = \diag(H_0, H_1, \ldots, H_{N_{sw}}),
\ms
D = \diag( B, Q_1, \ldots, Q_{N_{sw}})
\eeqn
and
\[
R = \diag(R_0, R_1, \ldots, R_{N_{sw}}).
\]
(Note the incorporation of the background covariance matrix $B$ in $D$. Also
note that we have eschewed correlation across time windows, as is often done
in practice.) The approximate minimization of the quadratic subproblem is
itself carried out using a Krylov method (often conjugate gradients
\cite{HestStie52}, GMRES \cite{SaadSchu86} or efficient specialized techniques
such as RPCG \cite{GratTshi09} or RSFOM \cite{GratToinTshi11}, see also
\cite{GratGuroToinTshiWeav12}).  These iterative methods typically requires
preconditioning for achieving reasonable computational efficiency.

Three variants of the above problem can then be defined.  In the form
presented above, the formulation is called the ``state formulation'' and
its optimality condition is given by the linear system
\beqn{state-formulation}
(L^T D^{-1} L + H^TR^{-1} H) \; \delta{x}
= L^T D^{-1} b + H^TR^{-1}d.
\eeqn
Another version (called the ``forcing formulation'') may be obtained by
making the change of variables ${\delta p} = {L \delta x }$, then requiring
the solution of the minimization problem
\beqn{QuadProb-fo}
\min_{ {\delta p} \in \Re^s} \; q_{\rm fo}(\delta p) \eqdef
     \half \| {\delta p} -b\|_{D^{-1}}^2
   + \half \|  H L^{-1}\, {\delta p} - d \|_{R^{-1}}^2
\eeqn
whose optimality condition may now be written as
\beqn{forcing-formulation}
(D^{-1} + L^{-T} H^TR^{-1}HL^{-1}) \; \delta{p}
=  D^{-1} {b} + L^{-T} H^TR^{-1}{d} \eqdef b_{\rm fo}.
\eeqn
We immediately note that \req{forcing-formulation} may be obtained as a
two-sided preconditioning of \req{state-formulation} with $L^{-T}$ and
$L^{-1}$. A third version (the ``saddle'' formulation) is obtained by
transforming the terms in \req{QuadProb-st} in a set of equality constraints
and writing the Karush-Kuhn-Tucker conditions for the resulting constrained
problem, leading to the large ``saddle'' linear system
\beqn{saddle-formulation}
   \left( \begin{array}{clc}
        D      & 0     & L \\
        0      & R     & H \\
        L^T    & H^T & 0
    \end{array} \right)
    \left( \begin{array}{c}
         \delta \lambda \\
         \delta \mu \\
         \delta x
    \end{array} \right) =
    \left( \begin{array}{c}
          b \\
          d \\
          0
    \end{array} \right)
\eeqn
where the control vector $[\delta \lambda^T, \delta \mu^T,\delta x^T]^T$ is a ($2s +
m$)-dimensional vector. For the sake of brevity, we do not cover the details
of this latter derivation here (see \cite{FishGuro17, FishGratGuroVassTrem17}):
it is enough to view it as an algebraic ``lifting'' of condition
\req{state-formulation} since this latter condition is recovered
by applying Gaussian block elimination to the first two rows and columns.
Unfortunately, this reformulation is only that of the optimality condition
and it is unclear whether it can be derived from an associated
quadratic minimization problem. Unless exact minimization is considered, this
will turn out to be problematic, as we will see below.

We immediately observe that matrix-vector products $u=L^{-1}v$ are sequential,
because they are defined by the simple recurrence
\beqn{seq} u_0 = v_0, \ms
u_i = v_i + M_iu_{i-1} \ms (i = 1, \ldots, N_{sw})
\eeqn
(a similar recurrence holds for products with $L^{-T}$), which is a serious
drawback in the context of modern computer architectures for high-performance
computing. In this respect, using the forcing formulation can be
computationally cumbersome and, even if the state and saddle point
formulations allow performing matrix-vector products with $L$ in parallel,
their suitable (and often necessary) preconditioners involve the operator
$\tilde{L}^{-1}$, where $\tilde{L}$ is a block bi-diagonal 
approximation of $L$ within which the matrices $M_i$ are replaced by
approximations $\tilde{M}_i$. The choice of such preconditioners is thus restricted to
use operators $\tilde{L}$ whose inversion can be parallelized, or whose
preconditioning efficiency is such that extremely few sequential products are
requested. Given the structure of $L$, this limits the possible options
(see \cite{GratGuroSimoToin17}). We have chosen, for the present paper, to
follow \cite{FishGuro17} and to focus on the two simplest choices: 
\beqn{tildeM}
\tilde{M}_i = 0,
\tim{ and }
\tilde{M}_i = I.
\eeqn
Once $\tilde{L}$ is determined, it remains to decide on the complete form of
the (left\footnote{Right preconditioning is also possible, but depends more on
  the detailed nature and discretizations of the dynamical models, which is
  why it is not considered here.}) preconditioner, depending on which of the
formulations \req{state-formulation}--\req{saddle-formulation} is used. For the
state formulation \req{state-formulation}, we anticipate the background term
to dominate and use the approximate inverse Hessian of the first term, given by
\beqn{S-def}
S^{-1} = \tilde{L}^{-1}D\tilde{L}^{-T}.
\eeqn
If the forcing formulation \req{forcing-formulation} is considered, an obvious
choice is to use $D$ as preconditioner (thus yielding a system which can
be viewed as a low-rank modification of the identity). Its efficiency has been
considered in \cite{ElSa15,ElSaNichLawl17}. Finally, the choice is
more open for the saddle formulation \req{saddle-formulation}.  For the sequel
of this paper, we consider the preconditioners given by the inverse of the
matrices
\beqn{saddle-precs}
P_M =
   \left( \begin{array}{ccc}
        D              & 0      & \tilde{L} \\
        0              & R      & 0 \\
        \tilde{L}^T    & 0      & 0
   \end{array} \right),
   \ms
P_B =   
   \left( \begin{array}{ccc}
        D              & 0      & 0 \\
        0              & R      & 0 \\
        0              & 0      & -S
   \end{array} \right)
   \tim{and}
P_T = 
   \left( \begin{array}{ccc}
        D              & 0      & \tilde{L} \\
        0              & R      &  H \\
        0              & 0      &  S
   \end{array} \right),
\eeqn
$P_M$ being the inexact constraint preconditioner suggested in
\cite{BergGondVentZill07,BergGondVentZill11} and used in \cite{FishGuro17},
and $P_B$ and $P_T$ being the triangular and block-diagonal ones
inspired by \cite{BenzWath08} (see also \cite{Wath15}).

\numsection{The original saddle method}\label{saddle-original-s}

Armed with these concepts and notation, we may now consider the original
saddle technique as discussed in \cite{FishGuro17} (and also used in
\cite{FreiGree17}).  It is outlined as Algorithm~\ref{saddle-original},
where we define
\[
r(\delta \lambda, \delta \mu, \delta x )
  =  \left( \begin{array}{clc}
        D      & {0}   & L \\
        0      & R     & H \\
        L^T    & H^T   & 0
     \end{array} \right)
     \left(\begin{array}{c}
        \delta \lambda \\
        \delta \mu \\
        \delta x
     \end{array}\right)
     - \left(\begin{array}{c}
       b\\
       d\\
       0
     \end{array}\right).
\]
     
\algo{saddle-original}{SADDLE-original (SAQ0)}
{
  An initial $x_0$ is given as well as the correlation matrices
  $D$ and $R$, a maximum number of inner iterations $n_{\rm inner}$ and a
  relative residual accuracy threshold $\epsilon_r \in(0,1)$. Set $k=0$.\\
  While (not converged):
  \begin{enumerate}
  \item Compute $J(x_k)$ and $g_k = \nabla_x J(x_k)$.
  \item Apply the preconditioned GMRES algorithm \cite{SaadSchu86} to reduce
    $\|r(\delta \lambda, \delta \mu, \delta x )\|$
    using one of the left preconditioners given by \req{saddle-precs}.
    Terminate the GMRES iteration at inner iteration $j$ with $(\delta
    \lambda, \delta \mu, \delta x )$ if
    \beqn{gmres-term-orig}
    \|r(\delta \lambda, \delta \mu, \delta x )\| \leq \epsilon_r \Big( ||b|| + ||d||\Big)
    \tim{ or }
    j = n_{\rm inner},
    \eeqn
    yielding a step $\delta x_k = \delta x$.
  \item Set $x_{k+1} = x_k +\delta x_k$ and increment $k$ by one.
  \end{enumerate}
}

In practice the operational constraints often impose moderate values of
$n_{\rm inner}$ (a few tens) as well as a small number of outer iterations
(ten or less), and the improvement obtained by this fast procedure is often
sufficient for producing a reasonable forecast. Thus the concept of
convergence should be taken with a grain of salt in this context. However, the
monotonicity of the values of $J(x_k)$ underlying the convergence idea
remains important as a theoretical guarantee that the method is meaningful
from the statistical and numerical points of view, preventing the algorithm to
produce unreliable results.

In the state formulation, the monotonic decrease of the $J(x_k)$ is
promoted by the fact that the method used to minimize
\req{QuadProb-st} (conjugate gradient or one of its variants) is itself a
monotonic algorithm.  As a consequence, any decrease obtained for
\req{QuadProb-st} translates into a decrease for \req{wk4DVar} provided $q_{\rm
  st}$ reasonably approximates $J$ in the neighbourhood of $x_k$, as is often
the case (or can be enforced by a trust-region scheme \cite{ConnGoulToin00}).
A similar argument applies for the forcing formulation, where a
monotonic algorithm is also used to minimize $q_{\rm fo}(\delta_p)$ (remember
that this formulation can be derived from the state one by suitable two-sided
preconditioning).
However, the GMRES method used in
Algorithm~\ref{saddle-original} merely reduces the residual of the system
\req{saddle-formulation} without any variational interpretation. If $n_{\rm
  inner}$ is large enough for the residual to become ``sufficiently'' small
(producing a ``sufficiently'' accurate solution of \req{saddle-formulation}),
then the equivalence between the optimality conditions \req{state-formulation}
and \req{saddle-formulation} implies that the decrease in $q_{\rm st}$ at the
computed step is comparable to that which would be obtained by minimizing this
quadratic model exactly, thereby ensuring a suitable decrease in $J$.  The
difficulty is to quantify what is meant by ``sufficiently''.

To illustrate this point, let us consider an assimilation problem for the
one-dimensional nonlinear Burgers equation involving 100 discretization states
over 3000 time steps in 50 time sub-windows and 20 observations per subwindows
(the complete description of this problem is given in Appendix~A1). We apply
Algorithm~\ref{saddle-original} to this problem with the number of inner
iterations $n_{\rm inner}$ fixed to 50 and using 10 major Gauss-Newton
iterations, $\tilde{M}_i = 0$ and the preconditioners defined by
\req{saddle-precs}.  We use the abbreviation\footnote{The naming convention
  will become clearer in the sequel of the paper.} SAQ0-P-M to denote the
corresponding algorithmic variants, where P is the choice of preconditioner
type in \req{saddle-precs} and M is the particular choice of the model
approximations $\tilde{M}_i$, which can  be either I or 0.
 Figure~\ref{fig-sad-all-f} shows the resulting
evolutions of the values of $q_{\rm st}$ (dashed curve) and $J$ (continuous
curve) for SAQ0-M-0 and SAQ0-T-0 over all inner iterations\footnote{The values of $J$ at inner
  iterations have been computed for illustration purposes only: they are not
  needed by the algorithm.}, major iterations being indicated by vertical
dotted lines and the true minimum value of $J$ by an horizontal thick black
line.

\begin{figure}[htbp]
\vspace*{2mm}
\centerline{
  \includegraphics[height=6cm,width=14cm]{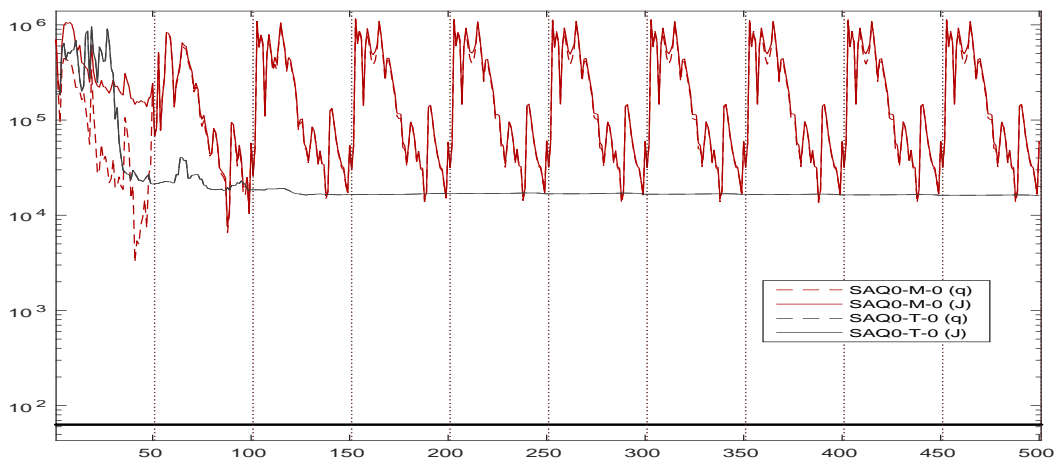} 
}
\caption{\label{fig-sad-all-f} Evolution of $q$ (dashed) and $J$
 (continuous) as a function of the total number of inner iterations for
  SAQ0-M-0 and SAQ0-T-0} 
\end{figure}

\noindent
Several important conclusions follow from the examining this figure.
\begin{enumerate}
\item None of the two methods achieves a significant reduction of the gap
  between $J(x_0)\approx 1.5\times 10^{5}$ and the optimal value ($\approx
  63.11$), the version using $P_M$ even diverging slowly. 

\item The curves for $q_{\rm st}$ and $J$ differ so little for iterations
  beyond the first that they are mostly undistinguishable, indicating a good
  fit between $q_{\rm st}$ and $J$ for moderately small steps. The observed
  stagnation/divergence may therefore not be blamed on the problem's
  nonlinearity.

\item The non-monotonic evolution of both $q_{\rm st}$ and $J$ along inner
  iterations is very obvious.  This is true for both SAQ0-M-0 and SAQ0-T-0 at
  the first iteration and for SAQ0-M-0 at all subsequent ones. We observe in
  particular that the value of $q_{\rm st}$ (and that of $J$) often starts by
  increasing at the first inner iterations. A stopping rule based on a maximum
  number of such iterations or on shortening the step therefore essentially
  relies on luck to produce a decrease in either objectives.
  
\item A qualitatively similar picture is obtained when using $\tilde{M}_i =
  I$, the only significant difference being that SAQ0-T-I now (slowly) diverges,
  despite the fact that the model approximations incorporates more
  information than for the (marginally more efficient) SAQ0-T-0.
\end{enumerate}

The numerical behaviour of for SAQ0-B-0 (the variant using the block diagonal
preconditioner $P_B$ from \req{saddle-precs}) is not shown in
Figure~\ref{fig-sad-all-f} because it would be barely visible.  As it turns
out, it examplifies to an extreme the fundamental difference between reducing
the residual of the system \req{saddle-formulation} and obtaining a decrease
in $q_{\rm st}$ or $J$. In this particular case, the preconditioner at
iteration two is good enough to ensure, with a single step of GMRES, a
relative reduction of the residual norm of the order of $10^7$, thereby
triggering a successful exit from the inner iteration loop.  However, since the first step
of GMRES is colinear with the initial residual $(b^T,d^T,0)^T$ (the right-and
side of \req{saddle-formulation}), this seemingly excellent step does not
alter the values of the state variables at all.  Thus $\delta x_k= 0$ for all
$k \geq 2$, causing the algorithm to stagnate (this would only
show as short horizontal line in the figure). The same undesirable behaviour
is also observed for SAQ0-B-I.

Extensive numerical experience with the Burgers assimilation problem indicate
that the conclusions drawn from this example are typical of many other
problems settings differing by conditioning of the involved correlation
matrices, number of time sub-windows or number of observations. Although they
obviously remain problem-dependent (as will be shown in
Section~\ref{numerics-s}), they show that, in general, \emph{the original
  saddle method} described in Algorithm~\ref{saddle-original} \emph{is
  potentially very inefficient or divergent} and that its efficiency might
decrease even if some algorithm's ingredients (such as model approximations)
are improved\footnote{It is quite remarkable that some practical
  implementations of the original saddle method actually skip Step~1 of
  Algorithm~\ref{saddle-original}. Since nor the objective function values $J$
  nor those of its gradient are ever computed, no
  convergence guarantee can possibly be given for this technique whose link
  with an optimization method becomes somewhat tenuous.}.

\numsection{The globalized SADDLE algorithm}\label{saddle-new-s}

Is it possible to fix this problematic behaviour of the original saddle
algorithm?  The answer is fortunately positive. Since we know that the desired
reduction in $q_{\rm st}$ is obtained if we allow GMRES to fully solve
\req{saddle-formulation}, the main idea is to adapt the GMRES termination rule
so that termination cannot occur before a minimal decrease in $q_{\rm st}$ is
obtained. We now detail this strategy as Algorithm~\ref{saddle-qm} where we use a
generic non-negative sequence $\{\theta_j\}\rightarrow 0$.

\algo{saddle-qm}{Globalized SADDLE (SAQ1)}
{
  An initial $x_0$ is given as well as the correlation matrices
  $D$ and $R$, a target number of inner iterations $n_{\rm inner}$ and a
  relative residual accuracy threshold $\epsilon_r \in(0,1)$.
  A model check frequency $\ell \in \Na$ and a model decrease threshold
  $\epsilon_q \in (0,1)$ are also given. Set $k=0$.\\
  While (not converged):
  \begin{enumerate}
  \item Compute $J(x_k)$ and $g_k = \nabla_x J(x_k)$.
  \item Apply the preconditioned GMRES algorithm  to \req{saddle-formulation},
    using one of the left preconditioners given by \req{saddle-precs}.
    At inner iteration $j$, 
    terminate with $(\delta \lambda,\delta \mu,\delta x)$ if 
    \beqn{gmres-qm-term2}
    q_{\rm st}(0) - q_{\rm st}(\delta x)
    \geq \max\Big[ \epsilon_q \min\left[ 1, \|g_k\|^2 \right], \theta_j \Big]
    \eeqn
    or if \req{saddle-formulation} is solved to full accuracy, yielding a step
    $\delta x_k = \delta x$.
  \item Perform a backtracking linesearch \cite[p. 37]{NoceWrig99} on $J$ along the
    direction $\delta x_k$, yielding
    $x_{k+1} = x_k +\alpha \delta x_k$ for some stepsize $\alpha
    >0$. Increment $k$ by one.
  \end{enumerate}
}

\noindent
It is clear that verifying \req{gmres-qm-term2} requires the (periodic)
evaluation of the quadratic model $q_{\rm st}(\delta x)$, which is an
additional computational cost: one needs to apply the $L, D^{-1}, H$ and $R^{-1}$
operators to obtain $q_{\rm st}(\delta x)$. The GMRES algorithm may also need
more than $n_{\rm inner}$ iterations to terminate, potentially increasing its
cost further. Note also that the linesearch procedure guarantees that
$J(x_{k+1}) \leq J(x_k)$ for all $k$.

Observe now that, because the covariance matrices are positive-definite, the
level set $\{x \in \Re^n \mid J(x) \leq J(x_0)\}$ is compact, and thus, using
the monotonicity of the algorithm, that there exists a constant $\kappa_g \geq 1$
such that $\|g_k\| \leq \kappa_g$, for all $k$. On termination of GMRES, we
therefore obtain that
\[
\epsilon_q \kappa_g^{-2}\|g_k\|^2
\leq q_{\rm st}(0) - q_{\rm st}(\delta x_k)
= -g_k^T\delta x_k - \half (L \delta x_k)^T D^{-1}(L \delta x_k)
                  - \half (H \delta x_k)^T R^{-1}(H \delta x_k).
\]
Using the positive-definite character of $\nabla^2 q_{\rm st}$, the Hessian of
$q_{\rm st}$, we deduce that
\beqn{gr1}
g_k^T\delta x_k \leq -\epsilon_q \kappa_g^{-2} \|g_k\|^2.
\eeqn
In addition, the strict
convexity of $q_{\rm st}$ ensures that \beqn{gr2} \|\delta x_k\| \leq
\frac{2}{\nu_{\min}}\,\|g_k\|, \eeqn where $\nu_{\min} > 0$ is the smallest
eigenvalue of $\nabla^2 q_{\rm st}$.  Thus \req{gr1} and \req{gr2} together
guarantee that $\delta x_k$ is ``gradient related'' in the sense that
\beqn{dxgr}
g_k^T\delta x_k \leq - \kappa_1 \|g_k\|^2
\tim{ and }
\| \delta x_k \| \leq \kappa_2 \|g_k\|
\eeqn for some positive constants $0< \kappa_1\leq \kappa_2$.  In conjunction
with the use of a linesearch, this well-known property of minimization
directions is then sufficient to ensure the monotonic decrease of the sequence
$\{J(x_k)\}$ and, assuming uniformly bounded condition numbers for $\nabla^2
q_{\rm st}(0)$ at all major iterations, the theoretical convergence of the
outer Gauss-Newton iteration 
(see \cite[Section~3.2]{NoceWrig99}, or
\cite[Section~10.1]{ConnGoulToin00}, for instance).  For the test
\req{gmres-qm-term2} to allow for early GMRES termination, it is necessary
that $\epsilon_q$ is chosen not too small.  Yet it should be small enough for
\req{gmres-qm-term2} to be attainable. If it is chosen too large, it may (in
the worst case) force GMRES to solve the system \req{saddle-formulation} to
full accuracy, in which case \req{gmres-qm-term2} is guaranteed with
$\epsilon_q\kappa_g^{-2} = \|[\nabla^2 q_{\rm st}]^{-1}\|$
\cite[Section~10.1]{ConnGoulToin00}. 

While obtaining convergence is often pratically out of reach or much too slow
in practice, the theoretical guarantee provides a strong
reassurance against potentially erratic results. Clearly, this argument also
holds for the choice $\theta_j=0$ for all $j$. The introduction of that
sequence is therefore unnecessary for ensuring mere convergence, but other
choices may be instrumental in speeding up decrease. For our Burgers example,
the choices \beqn{theta-def} \theta_j = \Big(\half q_{\rm st}(0)
\Big)^{\max[1,\sfrac{n_{\rm inner}}{j}]}-1 \tim{ and } \epsilon_g =
\sfrac{1}{100} \eeqn appear to give reasonable results.

We illustrate the behaviour of SAQ1 on the example used in the previous
section to highlight the difficulties of SAQ0. Its performance for the three
saddle preconditioners of \req{saddle-precs} is shown in Figure~\ref{fig-sa1q-all-f}.
  
\begin{figure}[htbp]
\vspace*{2mm}
\centerline{
  \includegraphics[height=6cm,width=14cm]{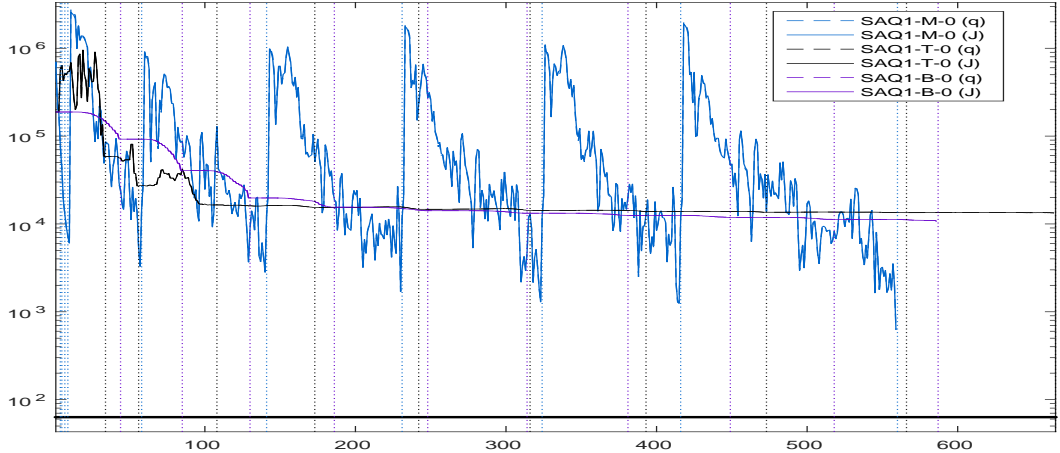} 
}
\caption{\label{fig-sa1q-all-f} Evolution of $q$ (dashed) and $J$
 (continuous) as a function of the total number of inner iterations for SAQ1
  (Burgers example)}
\end{figure}

\noindent
Comparing with the performance of the SAQ0 variants (Fig.~\ref{fig-sad-all-f}),
we may verify that SAQ1 achieves a significant reduction in $J$ in ten
Gauss-Newton iterations, albeit at the price of more inner iterations and the
cost of one valuation of $q_{\rm st}$ per inner iteration. SAQ1-M-0
oscillates most, but gets the best decrease, while the performance of
SAQ1-T-0 remains disappointing despite the theoretical guarantees and the
introduction of the $\{\theta_j\}$.  

If the cost of additional evaluations of $q_{\rm st}$ is high, one may try to
space them out by checking \req{gmres-qm-term2} only every $\ell >1$ inner
iteration.  The inner iteration termination rule \req{gmres-qm-term2} then becomes
\beqn{qm-termination}
\mod(j,\ell) = 0
\tim{ and }
q_{\rm st}(0) - q_{\rm st}(\delta x)
\geq \max\Big[ \epsilon_q \min\left[ 1, \|g_k\|^2 \right], \theta_j \Big],
\eeqn
and the resulting algorithms will be denoted by the abbreviation SAQ$\ell$ in
what follows.

The performance of SAQ25-P-0 (that is SAQ25 with preconditioner $P$ in
\req{saddle-precs} and $\tilde{M}_i=0$) is reported in
Figure~\ref{fig-saq15-all-f} for different preconditioners. As it turns out,
the relative performance of the methods including the preconditioners $P_T$
and $P_B$ is again poor.  This is unfortunately a constant in our experience
and we therefore focus on the use of the more successful $P_M$ only from now
on\footnote{We remark that this is not in contradiction with a comment in
  \cite{Wath15} that, \emph{for the symmetric case}, extended preconditioner
  formulations are unlikely to be efficient because of an alternating property
  in MINRES. Indeed, the square root of $P_M^{-1}$ is not well defined because
  it is indefinite, and the left-preconditioned system matrix is no longer
  symmetric (which is why GMRES is used).}.

We postpone the assessment of the sequential/parallel computational cost of
SAQ$\ell$ as a function of $\ell$ and the number of computing processes $p$ to
Section~\ref{numerics-s}, but we immediately notice that larger values of
$\ell$ may cause SAQ$\ell$ to require more inner iterations (as termination is
checked less often), in turn leading to larger memory and orthogonalisation
costs.  Thus a value of $\ell \leq n_{\rm inner}$ seems most
reasonable. However, more inner iterations may also result in a better
decrease of the quadratic model, and, if the problem is not too nonlinear, of
the overall objective function.

\begin{figure}[htbp]
\vspace*{2mm}
\centerline{
  \includegraphics[height=6cm,width=14cm]{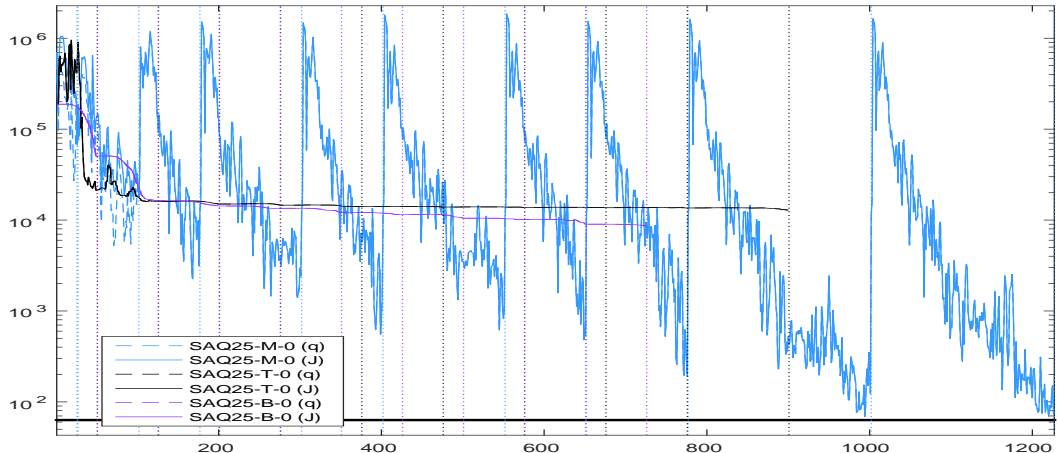} 
}
\caption{\label{fig-saq15-all-f} Evolution of $q$ (dashed) and $J$
  (continuous) as a function of the total number of inner iterations for SAQ25
(Burgers example)}
\end{figure}

If evaluating $q_{\rm st}$ is so costly that it must be avoided altogether, it is
still possible to cure the defects of the original formulation by relaying
directly on checking the gradient-related property of $\delta x$ stated by
\req{dxgr}.  However, the resulting algorithm then suffers more directly from
the need to estimate $\kappa_1$ and $\kappa_2$ a priori, and moreover appears
to be significanlty slower than the SAQ$\ell$ versions, mostly because it cannot
incorporate the forcing sequence $\{\theta_j\}$.  It will therefore not
be discussed here.

\numsection{Numerical comparisons on the Burgers example}\label{numerics-s}

We now turn to a comparative evaluation of the computational costs associated
with the SAQ$\ell$ algorithms, as well as state-of-the-art algorithms for
the state and forcing variational formulation. Our evaluation attempts to
provide conclusions in a context where parallel computing is available.

The comparison will involve several implementations of the state and
forcing formulations.  The first set (which, for now, we call the ST
algorithms) uses a standard Gauss-Newton algorithm where the
quadratic model $q_{\rm st}$ is minimized at each outer iteration using the
left-preconditioned Full Orthogonalization Method (FOM) \cite{Saad96}
expressed in the inner product defined by the inverse
preconditioner\footnote{In order to handle the unsymmetric matrix resulting
  from left-preconditioning, see Algorithm~\ref{foms} in Appendix~A3.}, which is in general preferable
\cite{GratToinTshi11} to the conjugate-gradient algorithm with
reorthogonalization (the practioner's most common choice so far).  The
preconditioner used is given by $S^{-1}=\tilde{L}^{-1}D\tilde{L}^{-T}$, the 
approximate inverse of the Hessian of the first term in 
\req{QuadProb-st}. The implementations of the forcing formulation (which we
call the FO algorithms) use the same left-preconditioned FOM method to
minimize $q_{\rm fo}$ given by \req{QuadProb-fo} as a function of $\delta
p$ using $D$ as a preconditioner and deduce $\delta x$ from
\beqn{fo-backsolve}
\delta x = L^{-1} \delta p.
\eeqn
All these methods use a standard trust-region scheme
\cite{ConnGoulToin00} for ensuring theoretical convergence. In these algorithms, the
value of $q_{\rm st}$ is readily available at the price of a single inner product at
the end of each inner iteration (remember that
$q_{\rm fo}(\delta p) = q_{\rm st}(L^{-1}\delta p)$).
It is therefore most coherent to terminate the inner iterations when
\req{qm-termination}-\req{theta-def} holds or, in the worst case, if the
relevant system has been solved to full  accuracy (i.e. residual norm below
$10^{-12}$ in our tests).

This lead us to relatively large set of algorithms, which differ by four
possible choices: the variational formulation (SA, ST or FO), the frequency
$\ell$ of the quadratic model check for terminating inner iterations in
\req{qm-termination}, the type of preconditioner used (inexact constraint M,
  triangular T or block-diagonal B for SA, the Schur complement S for ST, and
  the block-diagonal D for FO), and finally the choice of model approximation
  $\tilde{M}_i$ used in defining $\tilde{L}$ (0 or I, for SA and ST
  only). Using a naming convention coherent with that already introduced for
  the SAQ$\ell$ algorithms, we will, in the sequel, consider the algorithmic
  variants whose names are of the form AAQ$\ell$-P-M, where AA denotes the
  variational formulation, Q$\ell$ the frequency of the check for quadratic
  decrease in \req{qm-termination}, P the preconditioner type and M the choice
  of $\tilde{M}_i$, as summarized in Table~\ref{naming-t}.

\begin{table}[htbp]
\begin{center}
\begin{tabular}{|c|c|c|c|}
  \hline
  Var. Form. & quad. check freq. & $P$ type    & $\tilde{M}_i$ \\
  (AA)       &     ( Q$\ell$)    &   (P)       &   (M)         \\  
  \hline
  SA         & Q1, Q15, Q25, Q50 & M, T, B, n  &  0, I, M \\
  ST         & Q1, Q15, Q25, Q50 &    S, n     &  0, I, M \\
  FO         & Q1, Q15, Q25, Q50 &    D, n     &          \\
  \hline
\end{tabular}
\caption{\label{naming-t} Naming conventions for the considered algorithmic variants}
\end{center}
\end{table}

\noindent
Thus algorithm STQ15-S-0 uses the state formulation, checks for sufficient
quadratic decrease every 15-th inner iteration, uses the Schur complement
preconditioner in which $\tilde{L}$ is defined using $\tilde{M}_i= 0$.
In order to limit the number of variants, we have chosen $\ell \in \{ 1, 15,
25, 50\}$. We also introduced the 'n' preconditioner type, which stands for
not using preconditioning at all.
Altogether, we therefore obtain a set of 36 different algorithms\footnote{
  SAQ1-n,  SAQ1-M-0,  SAQ1-M-I,  SAQ1-M-M,
  SAQ15-n, SAQ15-M-0, SAQ15-M-I, SAQ15-M-M,
  SAQ25-n, SAQ25-M-0, SAQ25-M-I, SAQ25-M-M,
  SAQ50-n, SAQ50-M-0, SAQ50-M-I, SAQ50-M-M,
  STQ1-n,  STQ1-S-0,  STQ1-S-I,  STQ1-S-M,
  STQ15-n, STQ15-S-0, STQ15-S-I, STQ15-S-M,
  STQ25-n, STQ25-S-0, STQ25-S-I, STQ25-S-M,
  STQ50-n, STQ50-S-0, STQ50-S-I, STQ50-S-M,
  FOQ1-D,  FOQ15-D,   FOQ25-D,   FOQ50-D.}.

Before we embark on the comparison of computational costs, some further
comments on the various methods are in order.  The first is that the use of
preconditioning with the state formulation is not without risks.  Indeed,
\cite{GratGuroSimoToin17c} provides an analysis of the inherent difficulty of
preconditioning weighted least-squares problems, caused by the interplay
between the eigenstructures of $\tilde{L}$ and $D$. To illustrate how
problematic this can be, we borrow the following illustrative example from
this latter paper: let $\alpha >1$ be a parameter and
\[
L = \left( \begin{array}{cc}
  1      & 0\\
  \alpha & 1
\end{array}\right),
\ms
\tilde{L} = \left( \begin{array}{cc}
  1        & 0\\
  2+\alpha & 1
\end{array}\right)
\tim{ and }
D = \left( \begin{array}{cc}
  \alpha & 0 \\
    0    & 1 
\end{array}\right).
\]
Then it can be verified that $\tilde{L}^{-1}\tilde{L}^{-T}$ is a good
preconditioner of $L^TL$ in the sense that the condition number of
$\tilde{L}^{-1}\tilde{L}^{-T}L^TL$ is finite for all $\alpha$ (its value is
equal to 33.97 for all $\alpha>1$) while that
$\tilde{L}^{-1}D\tilde{L}^{-T}L^TD^{-1}L$ tends to infinity when $\alpha$
grows. This discussion suggests that a comparison between the various choices
of $\tilde{M}_i$ within the state formulation can be useful for each specific
problem.  Figure~\ref{fig-stfs-all} and other tests not reported here
indicate that using the choice $\tilde{M}_i = I$ (STQ15-S-I) is, in the Burgers
example, very inefficient compared with the three other choices. By contrast,
the choice of $\tilde{M}_i = M_i$ (STQ15-S-M) expectedly yields maximal
accuracy in a very small number of inner iterations.

\begin{figure}[htbp]
\vspace*{2mm}
\centerline{
    \includegraphics[height=6cm,width=14cm]{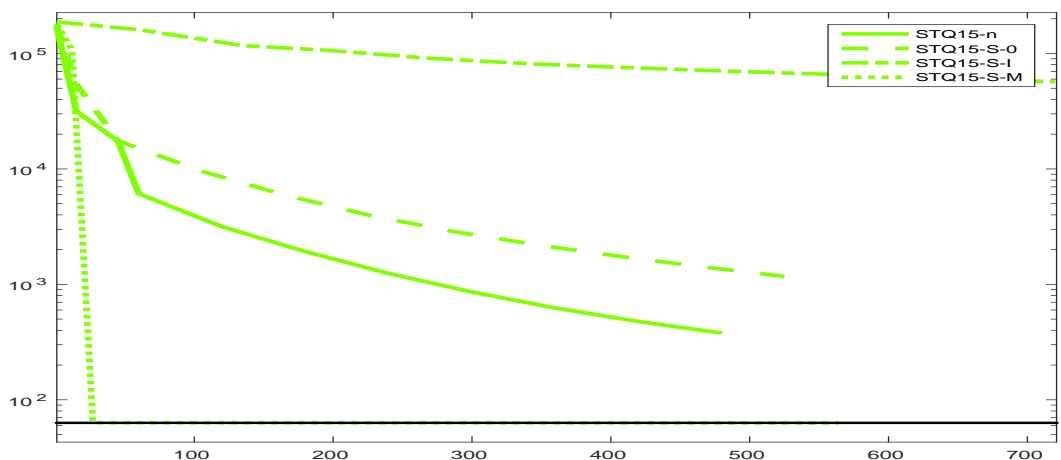}
}
\caption{\label{fig-stfs-all} Evolution of $J$ as a function of the number of
  inner iterations for STQ15-n, STQ15-S-0, STQ15-S-I, and STQ15-S-M (Burgers
  example)} 
\end{figure}

Other significant differences exist between STQ$\ell$ and FOQ$\ell$ algorithms
on one hand, and SAQ$\ell$ algorithms on the other hand.
\begin{enumerate}
\item At each inner iteration, the STQ$\ell$ algorithm compute a matrix-vector
  product with the Hessian of $q_{\rm st}$, which involve using the $D^{-1}$
  operator. The situation is better for the FOQ$\ell$-D variants, since
  (formally) the $D^{-1}$ term in the Hessian of \req{QuadProb-fo} is in this
  case premultiplied by $D$, which is readily simplified not to involve
  $D^{-1}$ at all.
\item The FOQ$\ell$-D methods also require products with $L^{-1}$ and $L^{-T}$
  at each inner iteration.  Unfortunately, these products are inherently
  sequential (see \cite{GratGuroSimoToin17}) and therefore potentially very
  costly. However, the quantity $L^{-1} \delta p$ can be recurred within the
  FOM algorithm itself at marginal cost, making the backsolve
  \req{fo-backsolve} unnecessary. A description of the resulting FOM algorithm
  (stripped from its trust-region enforcing features) is given as
  Algorithm~\ref{foms2} in Appendix~A3.
\item The use of the preconditioners $S^{-1}=\tilde{L}^{-1}D\tilde{L}^{-T}$ for the
  STQ$\ell$ methods is fully parallelizable, given our simple choices for
  $\tilde{M}_i$.
\end{enumerate}
Note that all methods use the $D^{-1}$ operator for jointly computing the
values of $J(x_k)$ and $g_k$ (once per major iteration) as well as for the periodic
evaluations of $q_{\rm st}$ every $\ell$-th inner iteration.
These remarks suggest that the two main parameters influencing the parallel
computing costs of the considered methods are the (parallelizable) cost of
computing $D^{-1}$ and that the purely sequential ones of computing $L^{-1}$
and $L^{-T}$. Due to the form of $L$ in \req{Ldef}, the two latter costs
are bounded by a small multiple of that of integrating the full nonlinear
model over the complete assimilation window. Our objective is therefore to
assess the efficiency of the various algorithms in a parametric study varying
the cost of applying $D^{-1}$.

The model used for the execution of parallel tasks is fairly simple, but it
hoped that it can nevertheless be sufficient for the broad type of analysis
presented. Let us denote by $c_{\rm op}$ the cost of evaluating (possibly in
parallel) the operator op. Then the cost of evaluating the tasks of costs
$c_1,\ldots, c_k$ in parallel on $p$ parallel computing processes is approximated by
\[
\pi_p(c_1,\ldots,c_k) = \max\left[\;
              \left\lceil \frac{k}{p} \right\rceil \frac{1}{k}\sum_{j=1}^k c_j,
              \;\max_{j=1,\ldots,k} c_j
              \;            \right].
\]

Taken alone, this approximation is not enough to provide the
description of a parallel computing environment, as it is crucial to consider
the impact of communications, which is beyond the scope of the present
paper. In what follows, we consider two complementary cases, and discuss their
associated parallel computing costs successively.

\subsection{A fully MPI approach}

A first, if somewhat restrictive, setting is to assume that the computation is
performed in $p$ MPI processes allowing the parallelizable operators to be
executed simultaneously for different time windows.  In
particular, this implies that parallel products with $L_i$, $L_i^T$,
$L_i^{-1}$, $L_i^{-T}$, $D_i$, $D_i^{-1}$ and $H_i$ are excluded
because each of them already uses the full available paralellism. By
the same argument, parallel products with $R_i$, $R_i^{-1}$  and $H_i^T$ are also banned. 

In this framework, the cost of evaluating $q_{\rm st}$ is given by
\beqn{cq-loc}
c_q =  c_L+c_{D^{-1}}+c_H+c_{R^{-1}},
\eeqn
while that of evaluating $J$ and its gradient is 
\beqn{cJ-loc}
c_J= c_\calM + c_\calH + c_{L^T}+c_{D^{-1}}+c_{H^T}+c_{R^{-1}}.
\eeqn
Note that the quantities $D^{-1}b$ and $R^{-1}d$ are both available once $J$
and its gradient have been evaluated.
We first investigate the computational costs of the components of GMRES, for which it
can be verified that the cost of a Krylov iteration for \req{saddle-formulation} is
\beqn{cKsa-loc}
c_{K,sa}= c_L+c_D+c_{L^T}+c_H+c_{H^T}+c_R
\eeqn
while that of applying the saddle preconditioner $P_M$ of \req{saddle-precs} involves
\beqn{cS-loc}
c_{S^{-1}} = c_{\tilde{L}^{-T}} + c_D + c_{\tilde{L}^{-1}},
\eeqn
(the cost of applying $S^{-1} = \tilde{L}^{-1}D\tilde{L}^{-T}$) and is
\beqn{cPM-loc}
c_{P_M} = c_{S^{-1}} + c_{R^{-1}}.
\eeqn
The cost of applying $n_o$ outer iterations of SAQ$\ell$-M for a total of $n_i$
inner iterations may then be approximated by
\beqn{cSAQiM-loc}
c_{SAQ\ell-M} \approx n_o(c_J+ c_{P_M}) + n_i(c_{K,sa}+c_{P_M}) + \frac{n_i}{\ell}c_q,
\eeqn
the second term in the first bracket of the right-hand side accounting for the
preconditioning, at each outer iteration, of the initial inner-iteration
residual, and the last term accounting for the periodic 
evaluations of $q_{\rm st}$ within the termination criterion. 
Similarly, it can be verified that, for the state formulation,
\beqn{cKst-loc}
c_{K,st} = c_L + c_{D^{-1}} + c_{L^T}+ c_H + c_{H^T} + c_{R^{-1}}
\eeqn
and thus, using \req{cS-loc} and 
\beqn{crhsst-loc}
c_{rhs,st} = c_{L^T}+c_{H^T}
\eeqn
the cost of computing the right-hand side of  \req{state-formulation}, that
\beqn{cSTF-loc}
c_{STQ\ell-S} \approx n_o(c_J+ c_{rhs,st} + c_{S^{-1}}) + n_i(c_{K,st}+c_{S^{-1}}).
\eeqn
Finally, for the forcing formulation,
\beqn{cKfo-loc}
c_{K,fo} =  c_D+c_{L^{-1}}+c_H+c_{R^{-1}}+c_{H^T}+c_{L^{-T}},
\ms
c_{rhs,fo} = c_{L^{-T}} + c_{H^T}
\eeqn
and
\beqn{cFOF-loc}
c_{FOQ\ell} \approx n_o(c_J+c_{rhs_fo}+c_D) + n_ic_{K,fo}.
\eeqn

We next assign approximate costs for all building blocks other than $D^{-1}$,
where one unit of cost is given by the integration of the model $\calM$ on the
complete time window\footnote{The cost values are based on discussions with
practitioners.}. Assuming $p$ computing processes are available and defining 
$e_{Nsw}$ to be the vector of all ones and length $N_{sw}$, let
\beqn{bcosts-1}
c_\calM = 1,
\ms
c_\calH = \frac{1}{20\,N_{sw}}\pi_p\Big( e_{Nsw} \Big),
\ms
c_D = \frac{1}{2\,N_{sw}}\pi_p\Big( e_{Nsw} \Big),
\ms
c_R = \frac{1}{100\,N_{sw}}\pi_p\Big( e_{Nsw} \Big),
\eeqn
\beqn{bcosts-2}
c_{R^{-1}} = \frac{1}{100\,N_{sw}}\pi_p\Big( e_{Nsw} \Big),
\ms
c_H = \frac{1}{10\,N_{sw}}\pi_p\Big( e_{Nsw} \Big),
\ms
c_{H^T} = \frac{1}{10\,N_{sw}}\pi_p\Big( e_{Nsw} \Big),
\eeqn
\beqn{bcosts-3}
c_L    =  \frac{2}{N_{sw}}\pi_p\Big( e_{Nsw} \Big),
\ms
c_{L^T} = \frac{4}{N_{sw}}\pi_p\Big( e_{Nsw} \Big),
\ms
c_{L^{-1}} = 2,
\tim{and}
c_{L^{-T}} = 4,
\eeqn
where we have used the block-diagonal structure of $D$, $R$, $R^{-1}$, $H$,
$H^T$, $L$ and $L^T$ to allow their costs to decrease with $p$.  Using these admittedly
fairly rough approximations, we may re-analyze the behaviour of the SAQ$\ell$
algorithms this time as a function of computational effort.

For the specific choice $c_{D^{-1}}= 0.5$, Figure~\ref{fig-seq-sas-all} corresponds to
Figure~\ref{fig-sa1q-all-f} where the horizontal axis now indicates sequential
computational costs, respectively (instead of
inner-iteration counts) and where the evolution of the nonlinear cost $J$ is
only shown (in thicker lines) between major iterations. The linesearch is
active at the first major iteration for all reported variants except
SAQ1-M-0, as is shown by the nearly vertical lines in the top left of the graph.

\begin{figure}[htbp]
\vspace*{2mm}
\centerline{
  \includegraphics[height=6cm,width=14cm]{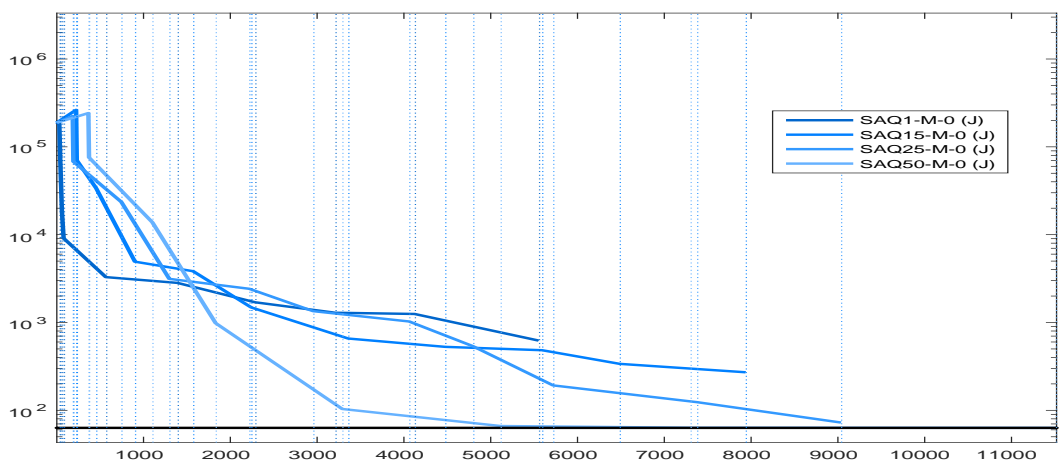} 
}
\caption{\label{fig-seq-sas-all} Evolution of $J$
  as a function of the sequential computational cost, for SAQ$\ell$-M-0,
  ($\ell=1,15,25,50$) using $c_{D^{-1}}= \half$ (Burgers example, fully MPI model)}
\end{figure}

\noindent
If we are now interested in plotting the parallel computational cost for the
same algorithms, the picture looks entirely similar (as all SAQ$\ell$-M-0
scale in the same way), but the maximal total cost appearing on the horizontal
axis shrinks from 11475 to 542 for 50 computing processes, a 21-fold speedup.
This illustrates the excellent parallelization potential of the
SAQ$\ell$-M-I methods, despite the chaotic (but controlled)
evolution of the quadratic model's values.  Similar plots and number could be
presented for the SAQ$\ell$-M-I and STQ$\ell$-S algorithms, all showing
reasonable parallelization potential.  As expected due to the use of the
inherently sequential operators $L^{-1}$ and $L^{-T}$, no such gains can be
obtained with the FOQ$\ell$-D variants, whose costs only vary marginally with
$p$. 

We are now left with the question of choosing a solution algorithm among our 36
SAQ$\ell$, STQ$\ell$ and FOQ$\ell$ variants ($\ell = 1,15,25,50$), depending
on the relative costs of applying $D^{-1}$. To answer this question, we first
applied each of the 36 methods to the Burgers example for
$c_{D^{-1}} \in [\half, 10]$.  We then
discarded all methods for which the total decrease in $J$ differed by a factor
more than $\rho \in (0,1)$ of the optimal decrease, that is
\[
J(x_0) - J(x_f) > \rho ( J(x_0) - J(x_*) )
\]
where $J(x_*)$ was obtained by running STQ1-S-M to full accuracy and $x_f$
denotes the final value of $x$ resulting from the application of the
algorithm. We finally selected the method for which the computational cost was
least for $p=1,10,25,50$. The maps indicating the winning method for each pair
$(c_{D^{-1}},\rho)$ with $ \rho\in [10^{-1},10^{-3}]$ and each $p$ are given
in Figures~\ref{fig-map-loc}. Each such map is accompanied with a picture of
the surface of the minimum computational costs over all $(c_{D^{-1}},\rho)$
pairs. The legend provides a correspondance between colors on the maps (the
left graph in each box) and the algorithmic variants\footnote{The colors used
  by MATLAB\copyright\ for the minimum cost surfaces at the right of each box
  are meaningless here.}.

\begin{figure}[htbp]
    \centerline{\fbox{\parbox[c]{7.2cm}{
        \centerline{$p=1$}
        \vspace*{2mm}
        \includegraphics[height=3.4cm,width=3.4cm,clip=true]{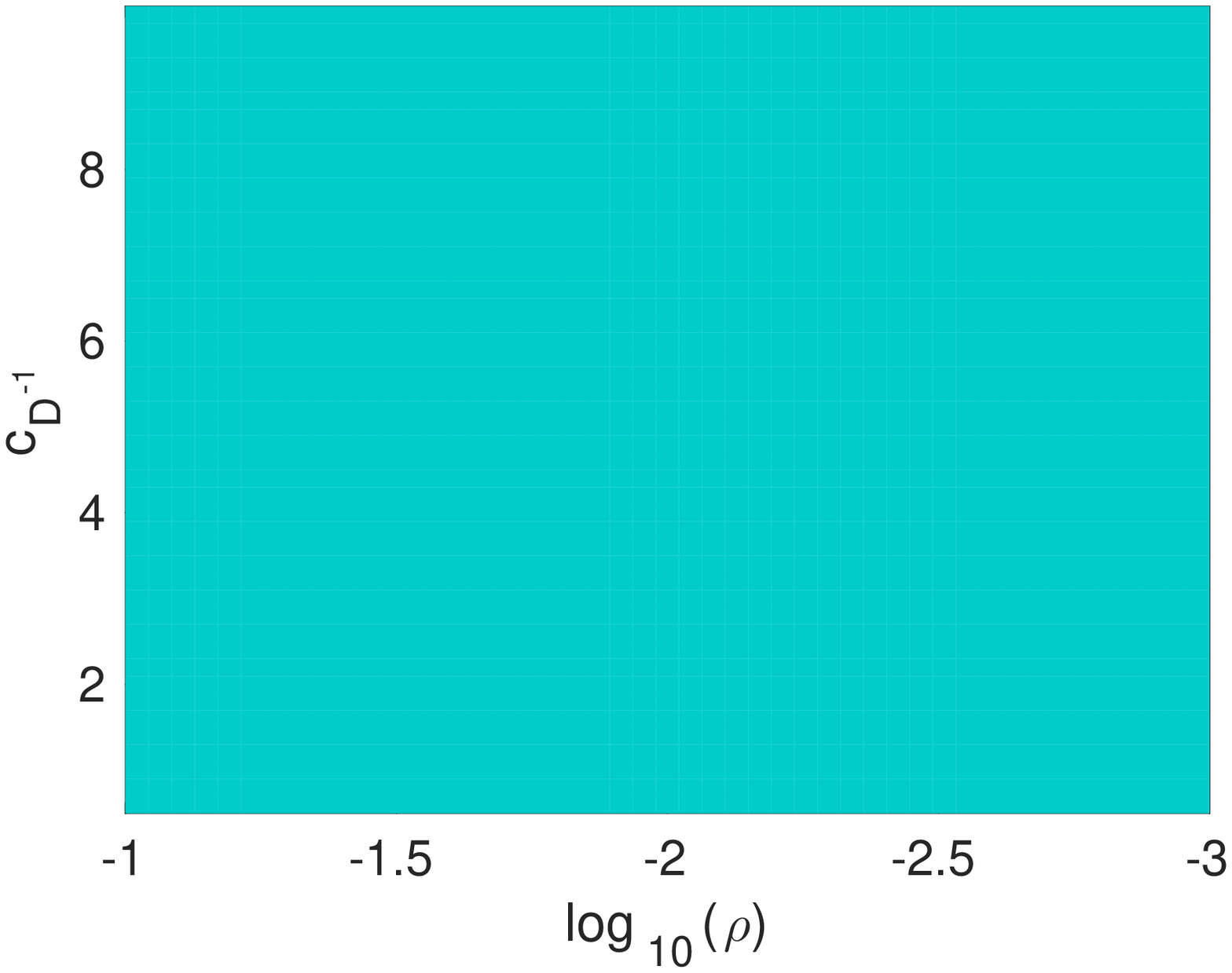}
        \hspace*{2mm}
        \includegraphics[height=3.4cm,width=3.4cm,clip=true]{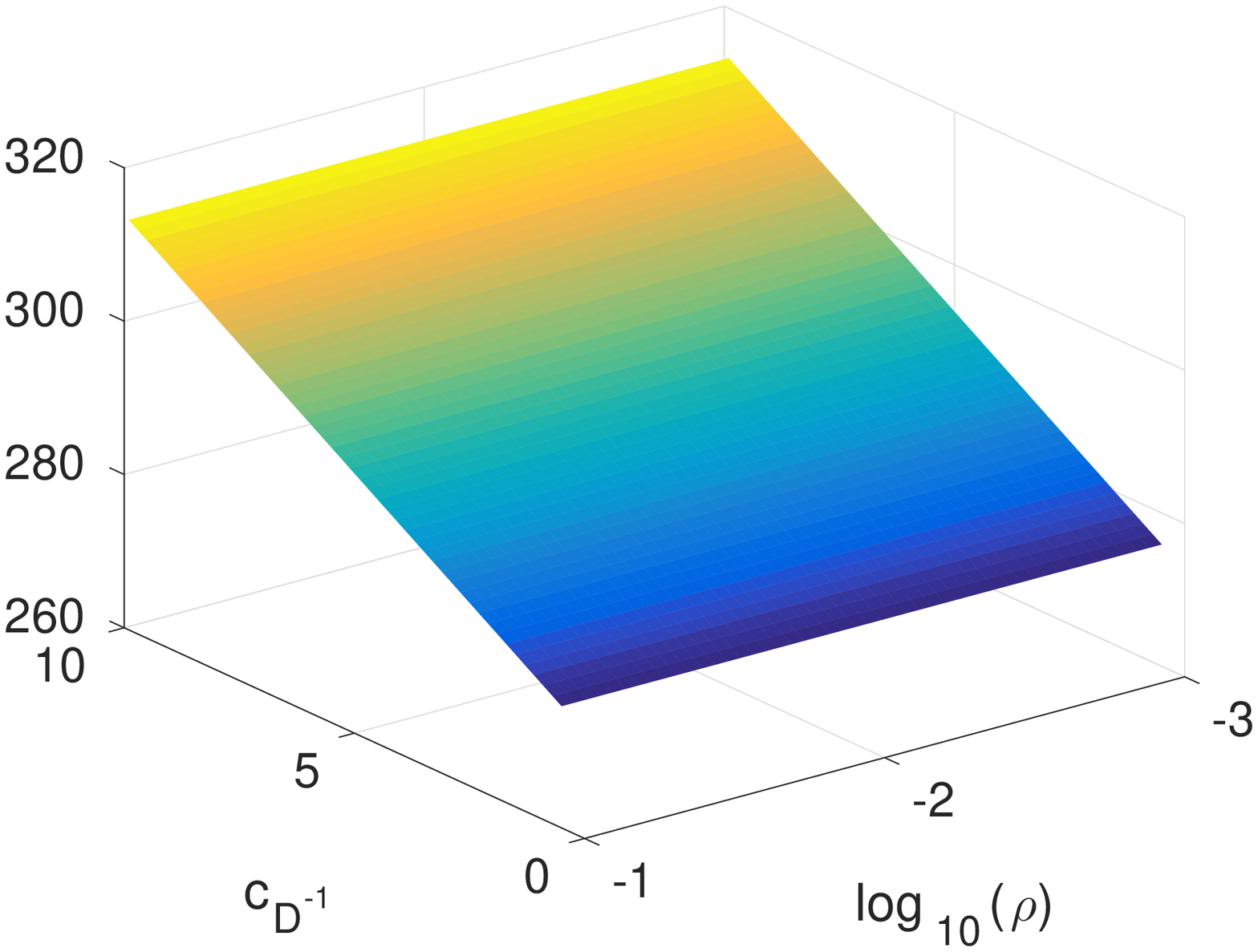}
      }}
    \fbox{\parbox[c]{7.2cm}{
        \centerline{$p=15$}
        \vspace*{2mm}
        \includegraphics[height=3.4cm,width=3.4cm,clip=true]{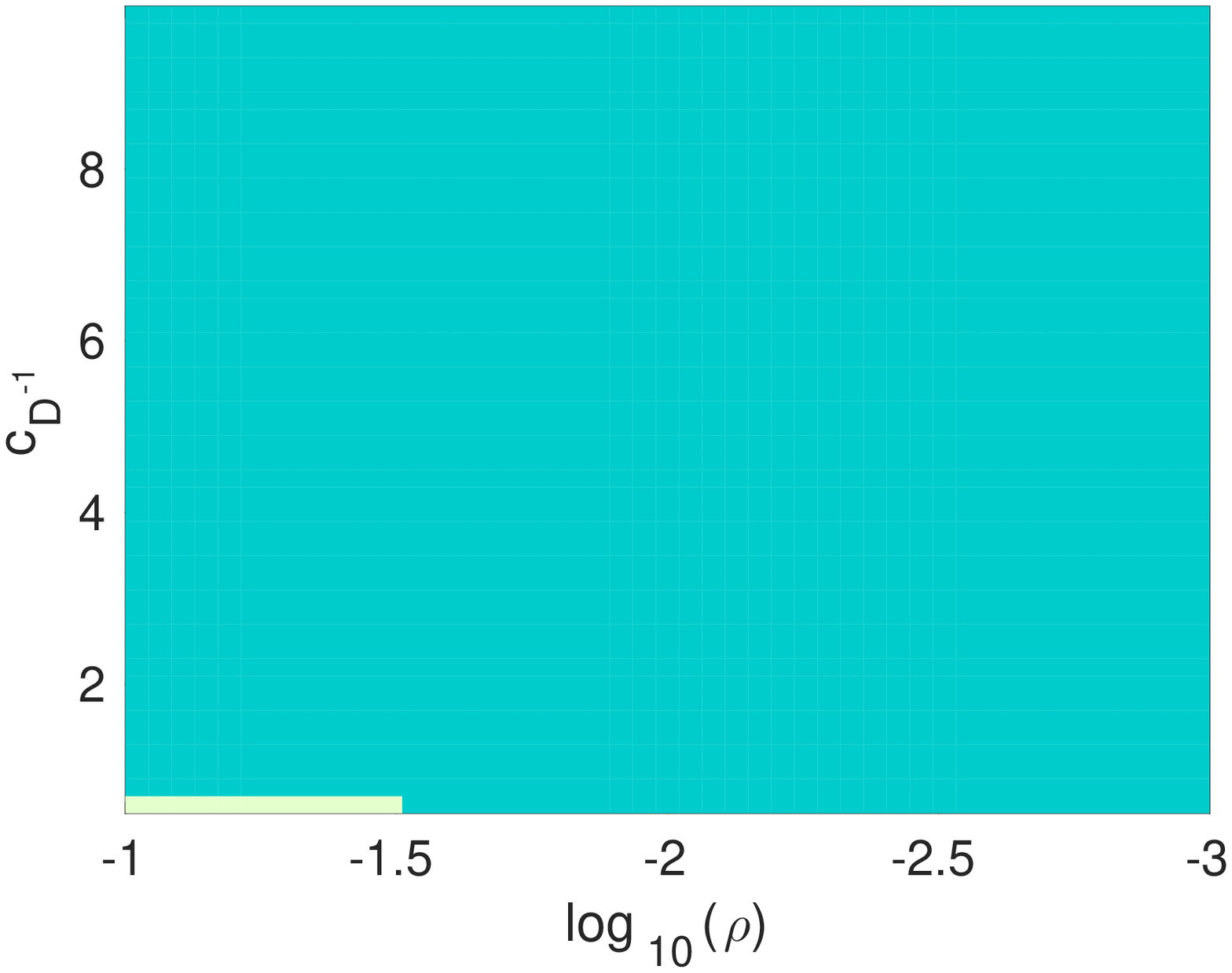} 
        \hspace*{2mm}
        \includegraphics[height=3.4cm,width=3.4cm,clip=true]{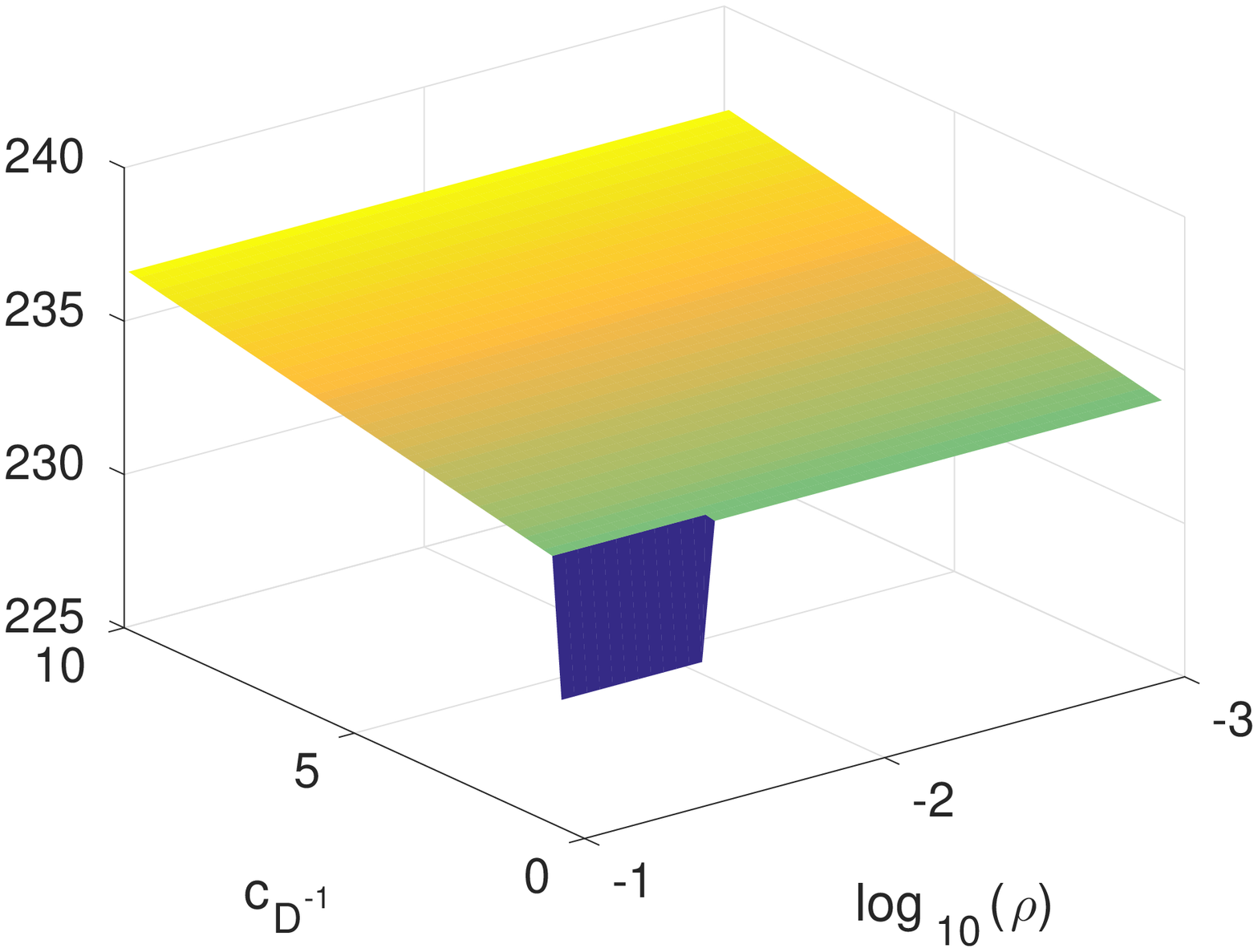}
    }}}
    \vspace*{1mm}
    \centerline{\fbox{\parbox[c]{7.2cm}{
        \centerline{$p=25$}
        \vspace*{2mm}
        \includegraphics[height=3.4cm,width=3.4cm,clip=true]{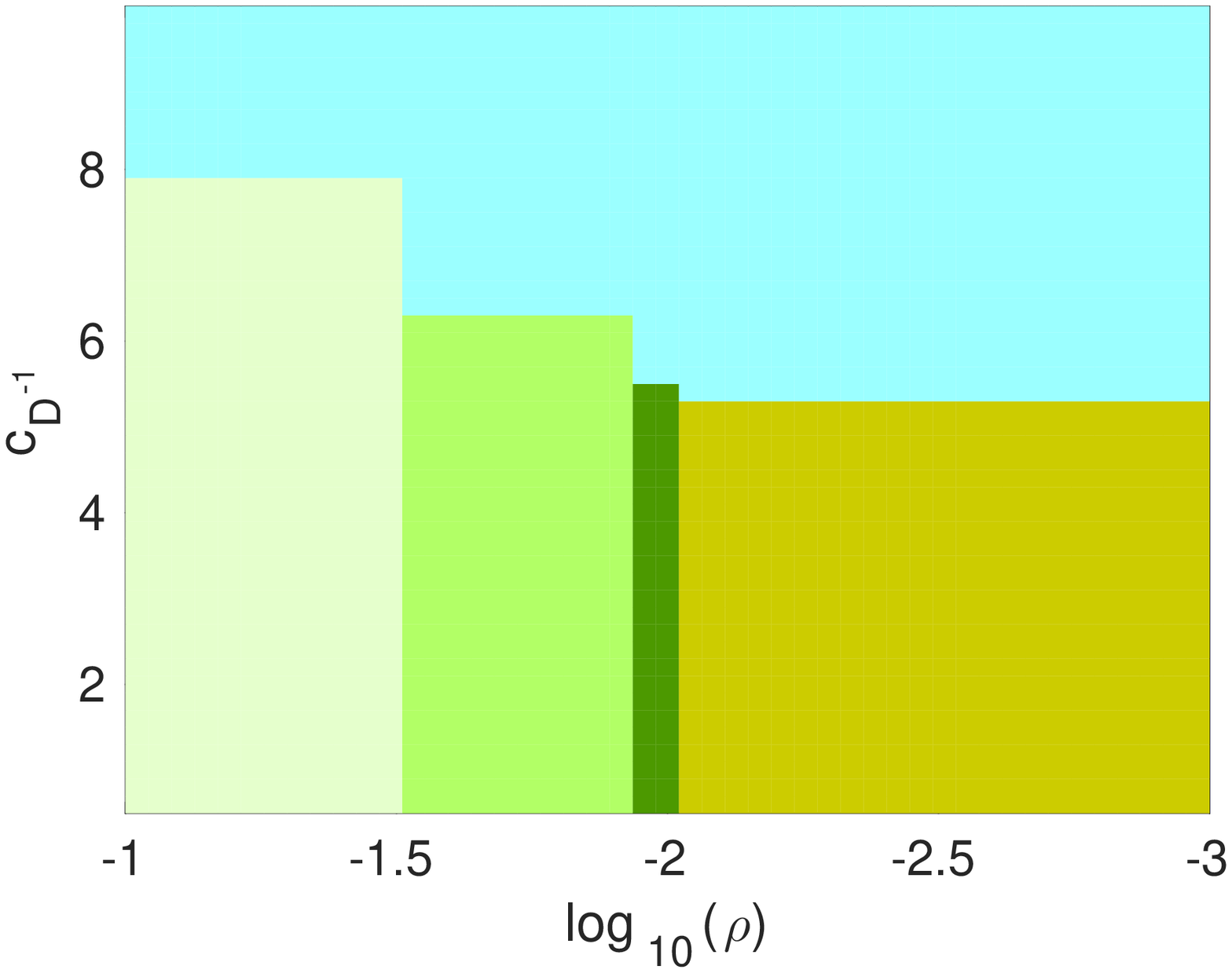} 
        \hspace*{2mm}
        \includegraphics[height=3.4cm,width=3.4cm,clip=true]{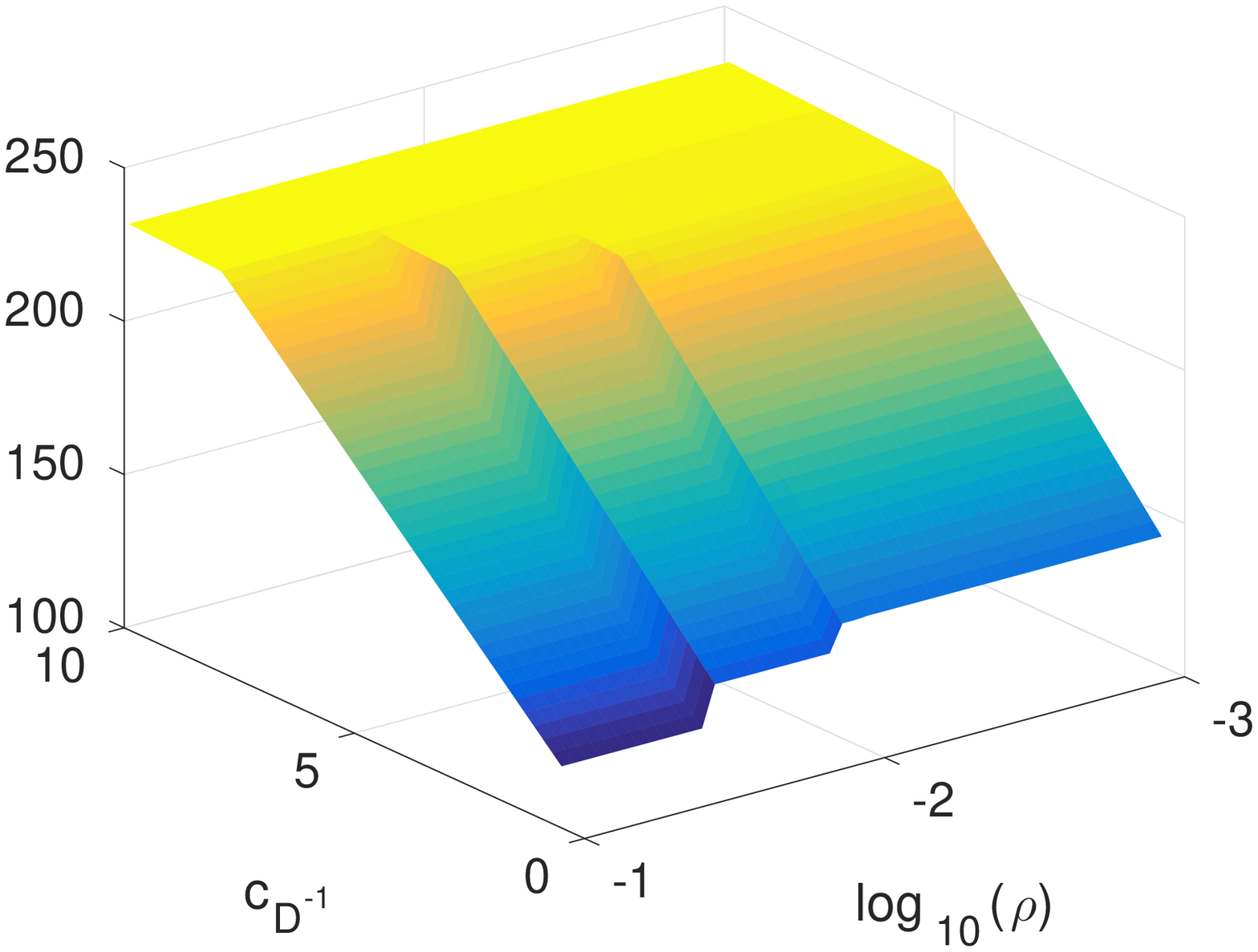}
      }}
    \fbox{\parbox[c]{7.2cm}{
        \centerline{$p=50$}
        \vspace*{2mm}
        \includegraphics[height=3.4cm,width=3.4cm,clip=true]{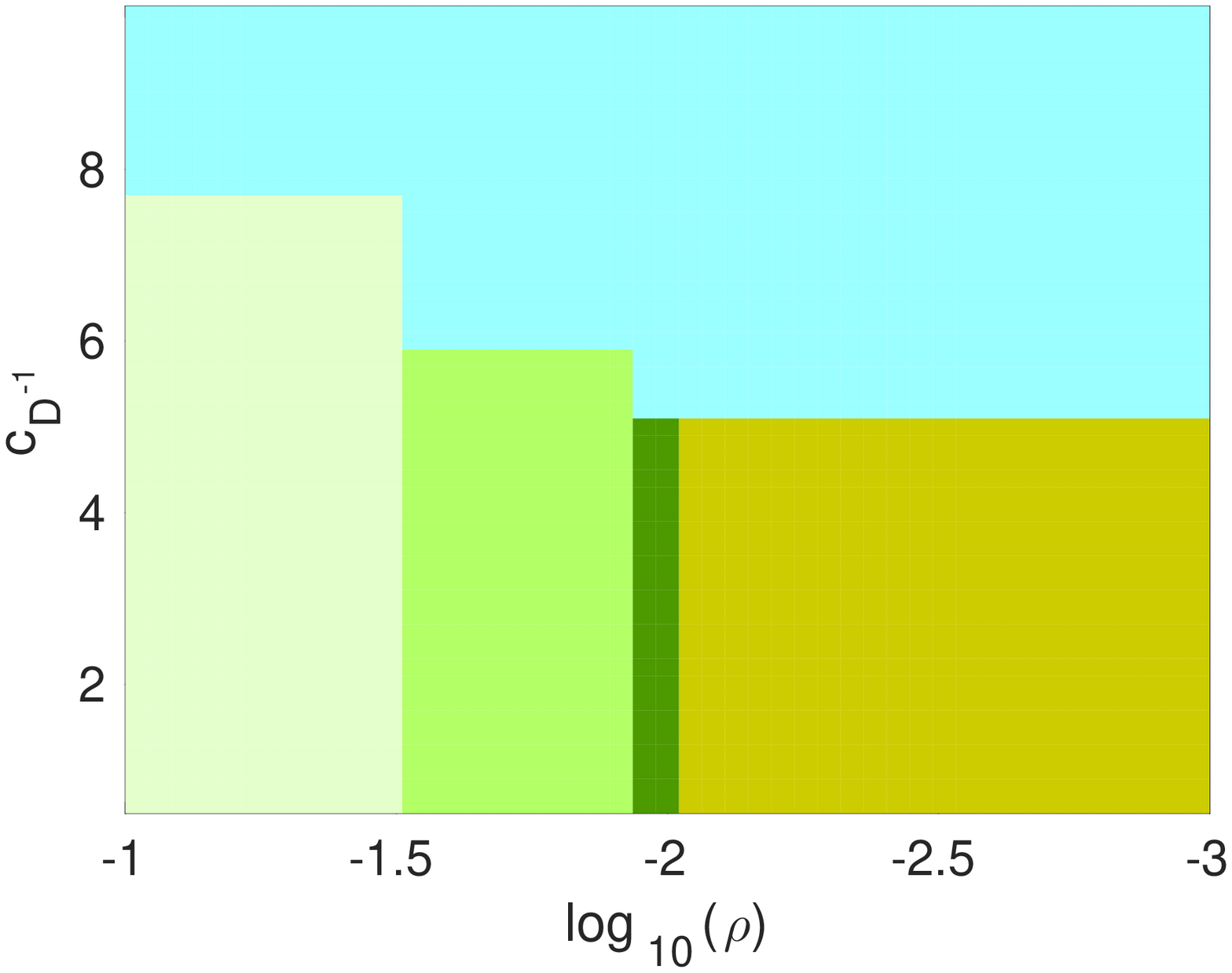} 
        \hspace*{2mm}
        \includegraphics[height=3.4cm,width=3.4cm,clip=true]{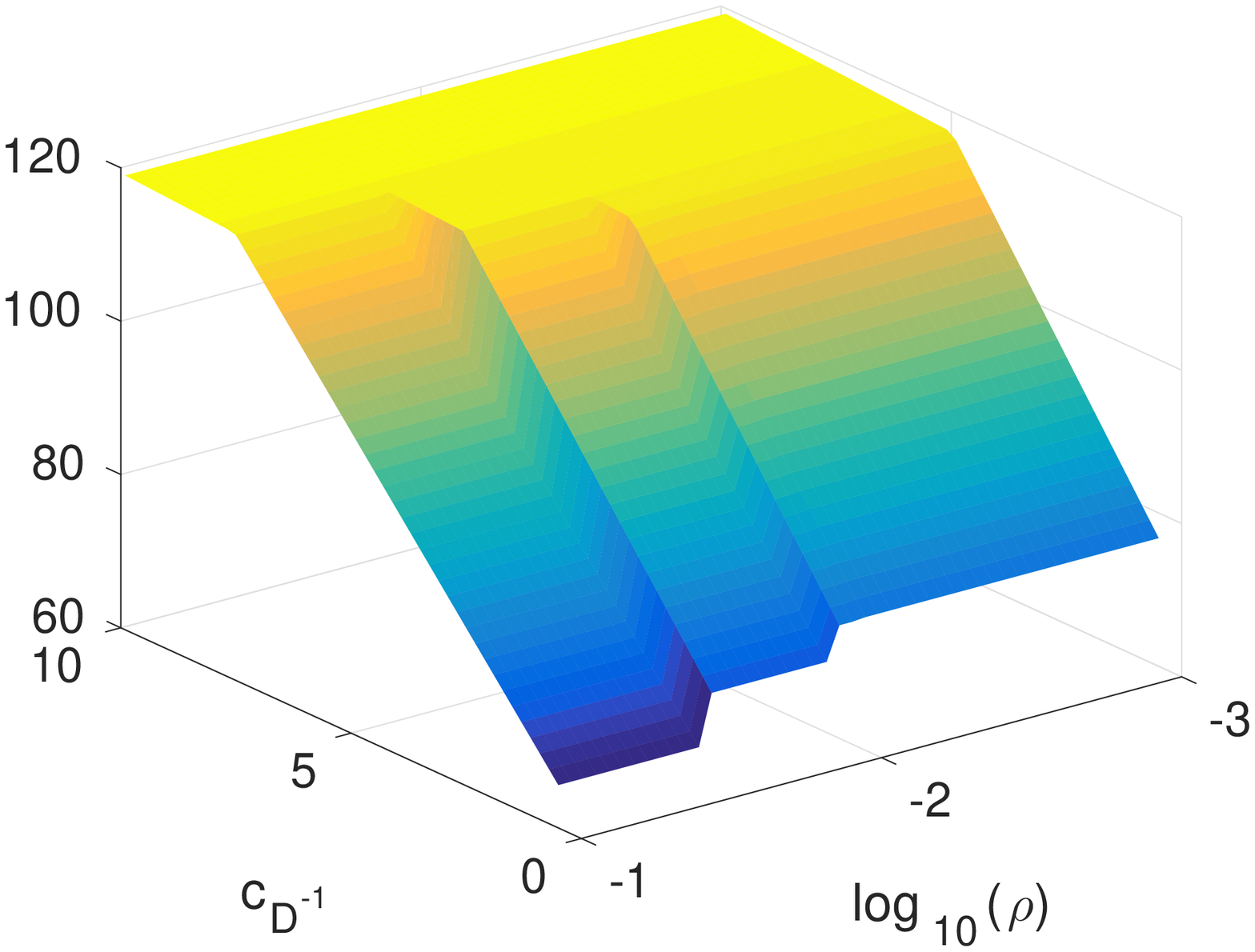}
      }}}
  \vspace*{1mm}
  \begin{center}
   Map colors: 
  \begin{tabular}{llllll}
  \includegraphics[height=3mm,width=3mm]{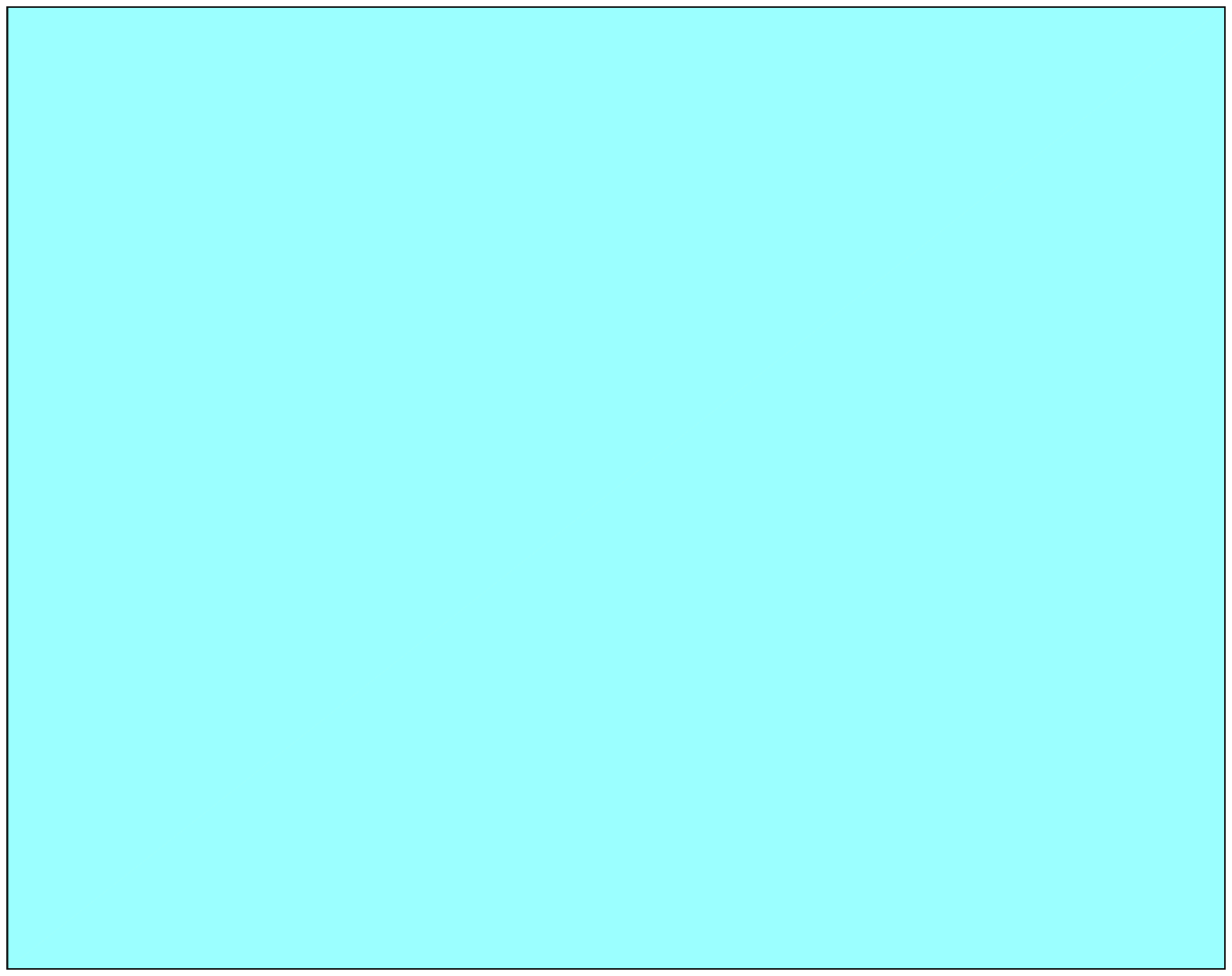} \hspace*{1mm} & SAQ50-M-I &
  \includegraphics[height=3mm,width=3mm]{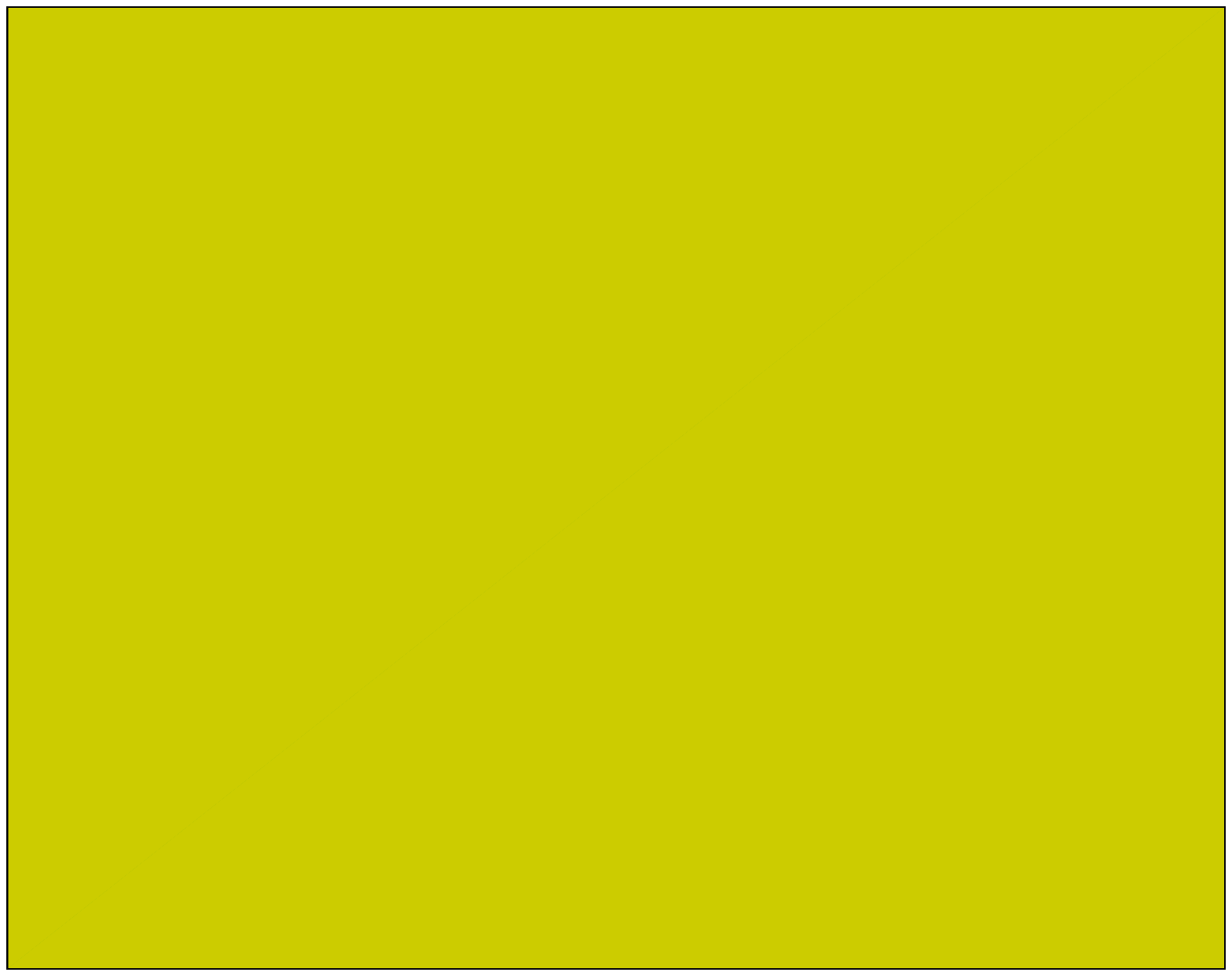}  \hspace*{1mm} & STQ1-n    &
  \includegraphics[height=3mm,width=3mm]{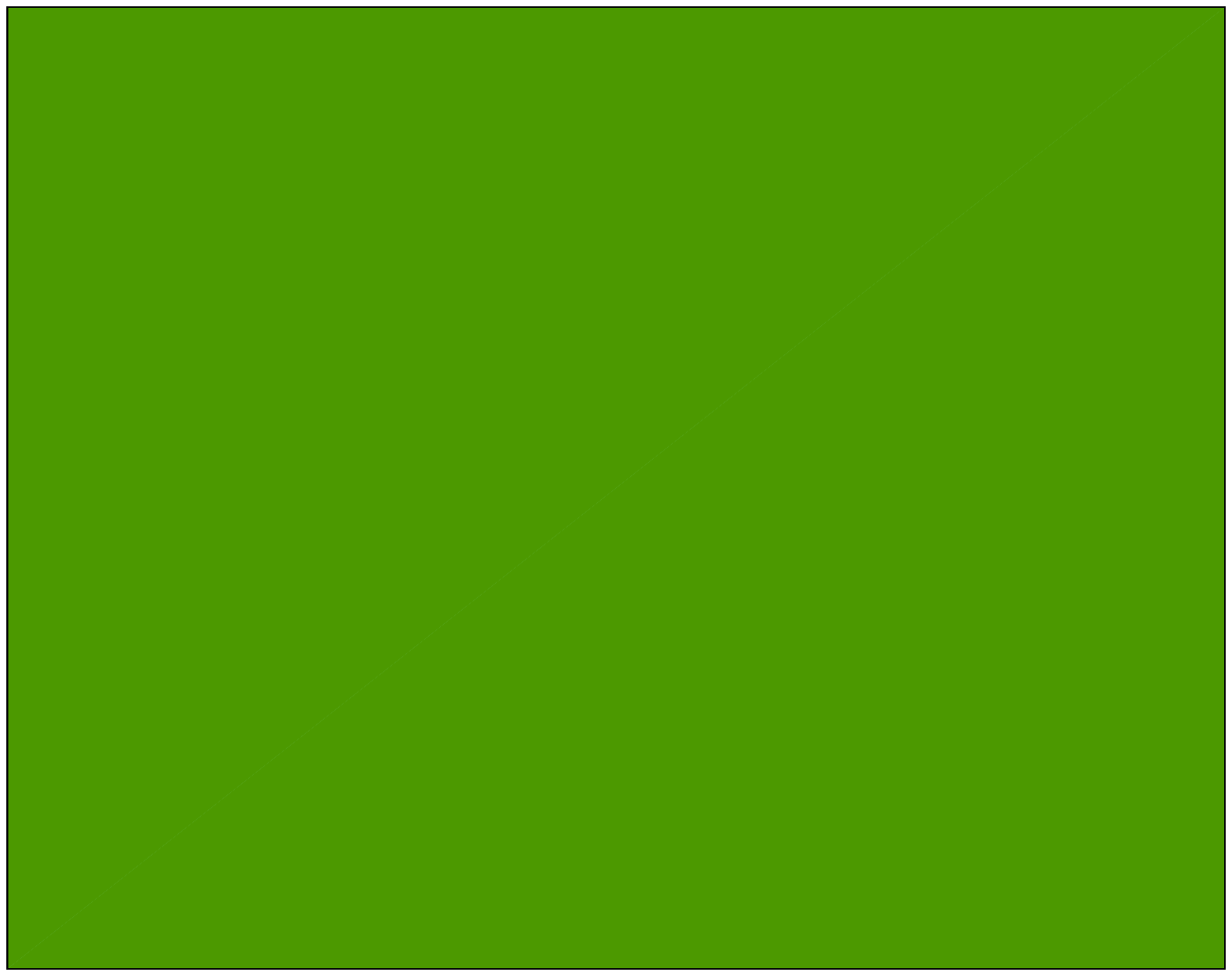}  \hspace*{1mm} & STQ1-S-0  \\
  \includegraphics[height=3mm,width=3mm]{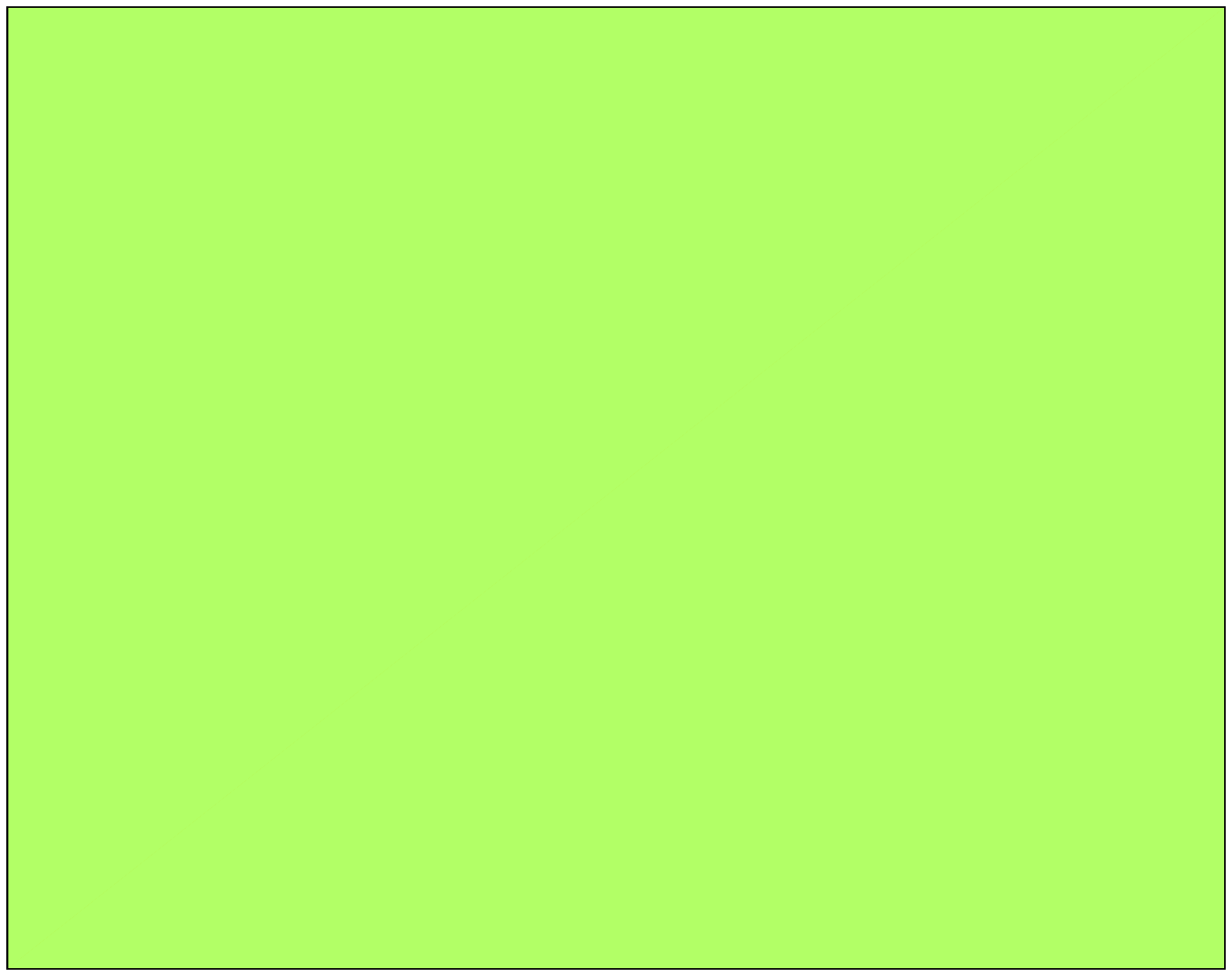} \hspace*{1mm} & STQ25-S-0 &
  \includegraphics[height=3mm,width=3mm]{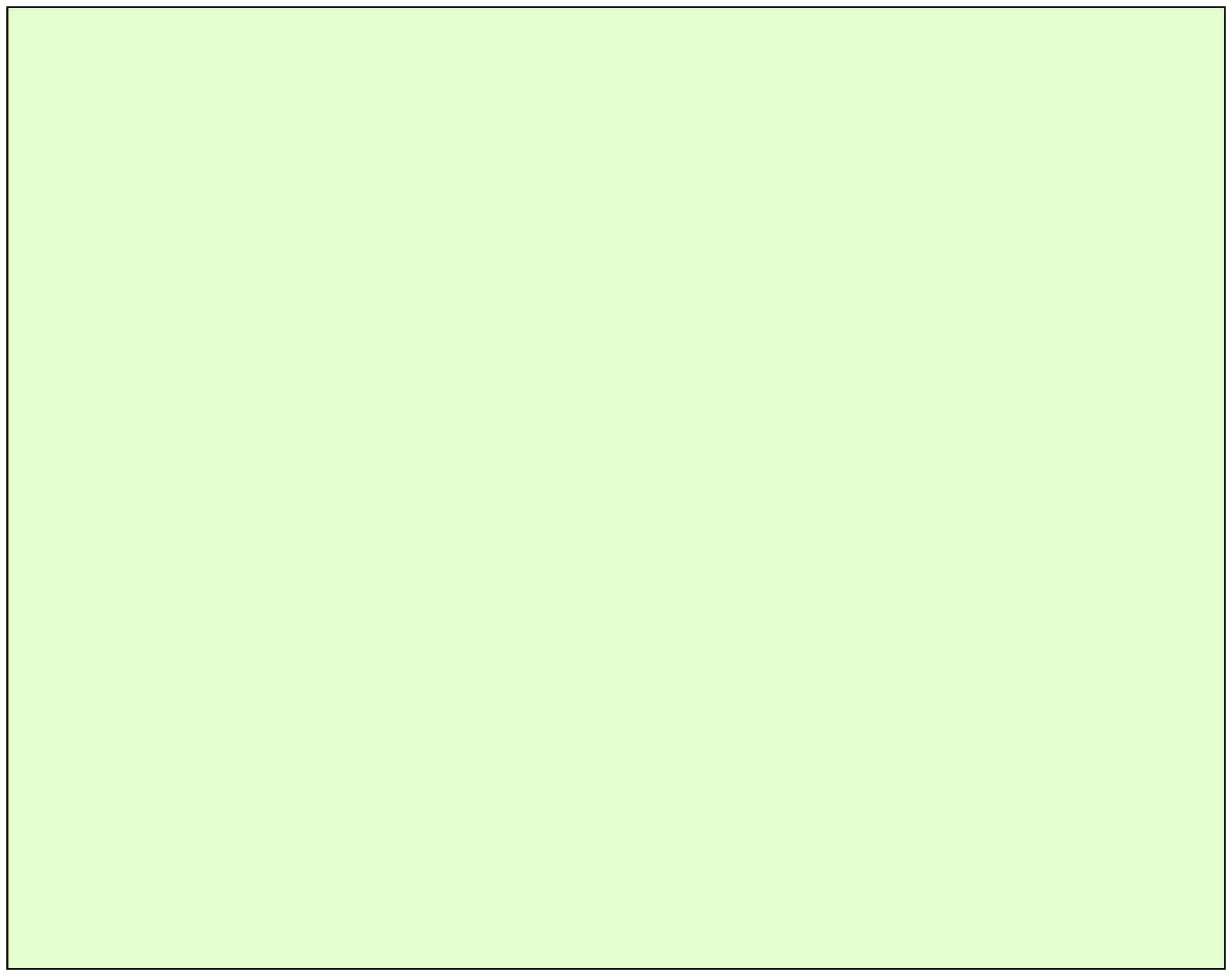} \hspace*{1mm} & STQ50-S-0 &
  \includegraphics[height=3mm,width=3mm]{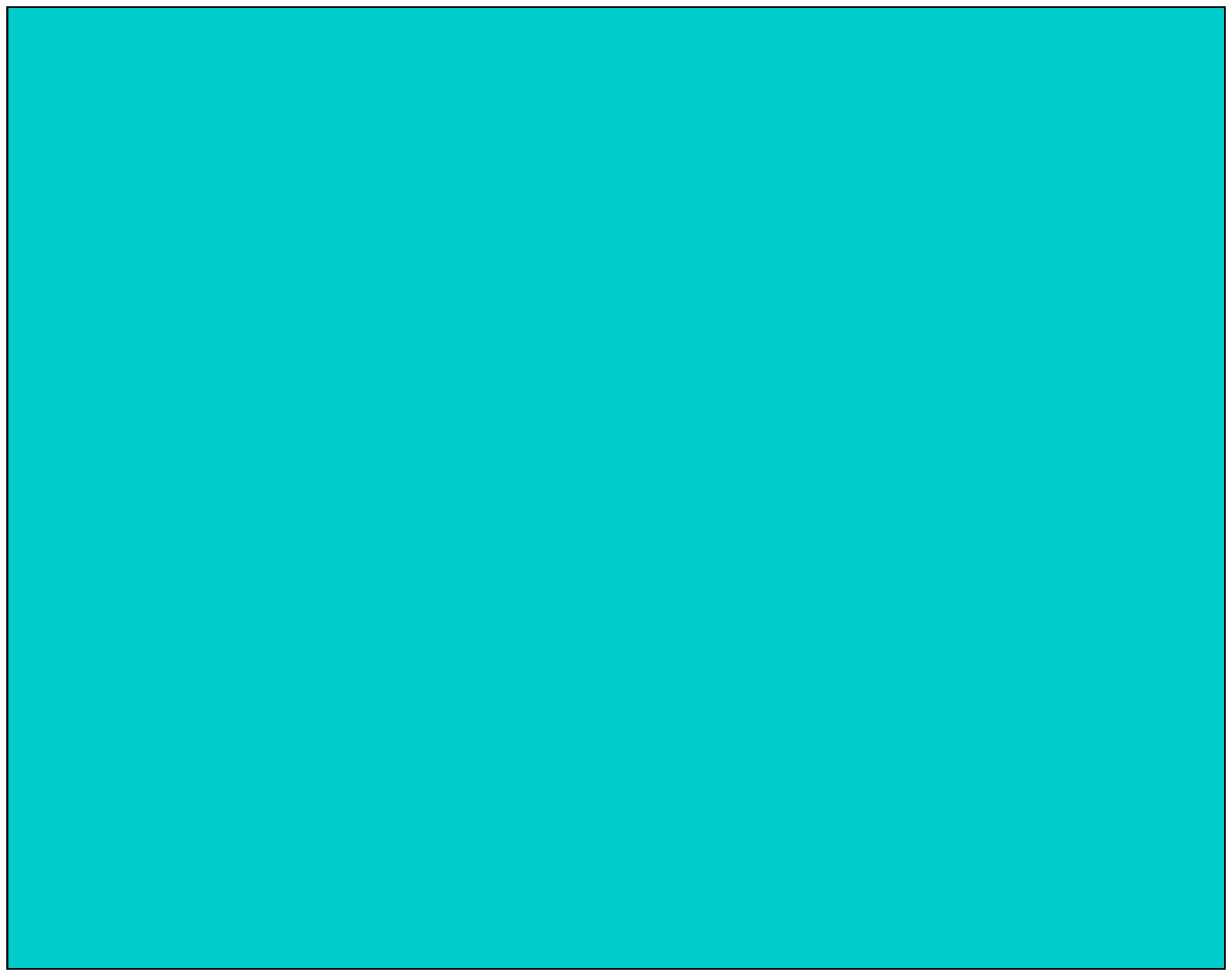}   \hspace*{1mm} & FOQ15-D   \\
  \end{tabular}
  \end{center}
\caption{\label{fig-map-loc}Best algorithmic variants as a function of
  $c_{D^{-1}}$  the number of computing processes $p$ and the reliability
  factor $\rho$ (Burgers example, fully MPI model)}
\end{figure}
        
Several conclusions \emph{for the Burgers example} may be drawn from these maps and
cost surfaces.
\begin{enumerate}
\item Using the forcing formulation dominates all other variants when computations are
  sequential, for a computing cost with maximal values close of 313.5. We note
  the proportionnally faster increase of the total cost with $c_{D^{-1}}$
  compared to requiring higher accuracy . It is interesting that the variant FOQ15-D is
  best, indicating that terminating the inner iterations as soon as possible
  (as is the case with FOQ1-D) may be sub-optimal when the preconditiner is excellent.

  The usefulness of the forcing formulation clearly decreases when the number
  of computing processes increases, as expected.
 
\item Using a saddle-based algorithm clearly supposes a high value of
  $c_{D^{-1}}$ and the availablility of several computing processes. The frequency
  $\ell=50= n_{\rm inner}$ appears to provide the best compromise between good
  decrease on the quadratic model and excessive number of inner iterations.
  
\item When $p=25$ or $p=50$ and $c_{D^{-1}}$ is moderate, the algorithms using
  the state formulation dominate with a frequency $\ell$ diminishing for
  increasing accuracy, the unpreconditioned version being suitable for maximum
  accuracy. This is coherent with our comments in the beginning of this
  section and Figure~\ref{fig-stfs-all}, as is the fact that none of the
  STQ$\ell$-S-I variants ever appears on the podium of best methods.
  
\item Compared to the sequential case, the most parallel of our scenarii
    ($p=50$) provides a reduction of computational costs from 315 to
  119, which corresponds to a speed-up of approximately 2.6.
  
\item The additional computing cost necessary for obtaining improved accuracy
  is negligible, irrespective of the number of computing processes.
\end{enumerate}

\subsection{An hybrid MPI/OpenMP approach}

Let us now assume that a more elaborate parallel computing environment is
considered, such as an hybrid  MPI/OpenMP system, where we assume that $p$
processes are available, each of which with two computing cores.  This means
that we may now apply two time-parallel operators simultaneously.  Then the
computing costs \req{cq-loc}-\req{cSTF-loc} may be rewritten as 
\beqn{cq-nol}
c_q =  \pi_2\Big(c_L+c_{D^{-1}},c_H+c_{R^{-1}}\Big),
\eeqn
\beqn{cJ-nol}
c_J
= c_\calM + c_\calH + \pi_2\Big(c_{L^T}+c_{D^{-1}},c_{H^T}+c_{R^{-1}}\Big),
\eeqn
\beqn{cKsa-nol}
c_{K,sa}= \pi_2\Big(c_L + c_D + c_H, c_{L^T} + c_R + c_{H^T} \Big),
\eeqn
\beqn{cS-nol}
c_{S^{-1}} = c_{\tilde{L}^{-T}} + c_D + c_{\tilde{L}^{-1}},
\ms\ms\ms
c_{P_M} = \pi_2\Big( c_{S^{-1}}, c_{R^{-1}}\Big).
\eeqn

\beqn{cSAQiM-nol}
c_{SAQ\ell-M} \approx n_o(c_J+c_{P_M}) + n_i(c_{K,sa}+c_{P_M}) + \frac{n_i}{\ell}c_q,
\eeqn
\beqn{cKst-nol}
c_{K,st} = \pi_2\Big( c_L + c_{D^{-1}} + c_{L^T}, c_H + c_{R^{-1}} + c_{H^T}  \Big),
\ms
c_{rhs,st} = \pi_2( c_{L^T}, c_{H^T} ),
\eeqn
\beqn{cSTF-nol}
c_{STQ\ell-S} \approx n_o(c_J+c_{rhs,st}+c_{S^{-1}}) + n_i(c_{K,st}+c_{S^{-1}}),
\eeqn
$c_{K,fo}$, $c_{rhs,fo}$ and $c_{FOQ\ell-D}$ being unmodified.

We may then repeat our experiments in this new setting, which allows
further gains compared to the fully MPI case, as shown in
Figure~\ref{fig-dataloc-all} (note the doubly logarithmic axis).

\begin{figure}[htbp]
\vspace*{2mm}
\centerline{
  \includegraphics[height=6cm,width=14cm]{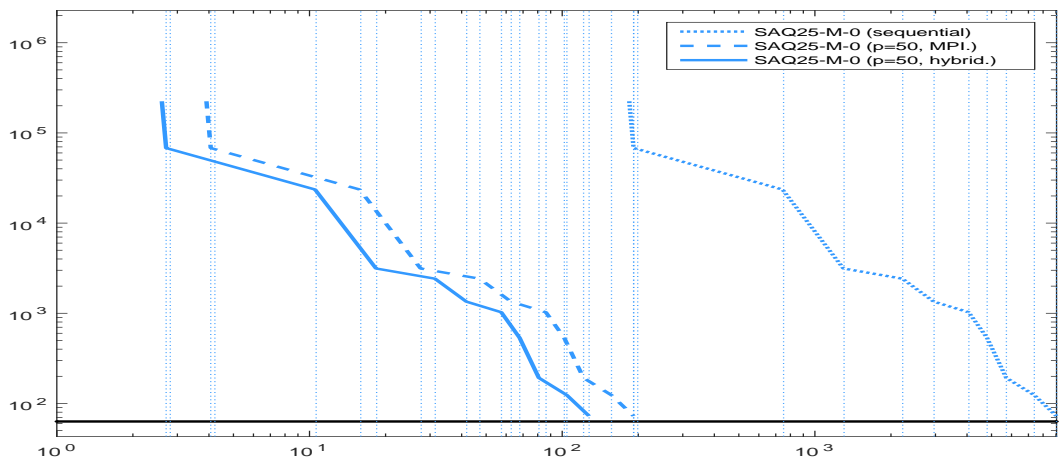}
}
\caption{\label{fig-dataloc-all} Evolution of $J$ as a function of
  computational cost for SAQ25-M-0 for $c_{D^{-1}} = \half$ in the
  sequential, fully MPI and hybrid MPI/OpenMP settings
  (Burgers example)} 
\end{figure}

\noindent
Although the improvement by switching from the fully MPI model to the hybrid
MPI/OpenMP model is far from negligible, we
neverthless note that most of the advantage obtained by parallel processing is
due to the parallelization in time, with a very strong correlation with
$N_{sw}$, the number of subwindows (see \req{bcosts-1}\req{bcosts-3}).  This
already apparent at the very beginning of the computation, as evaluating
$J(x_0)$ already makes a significant difference (both parallel runs of the
algorithm are completed before $J(x_0)$ is evaluated in the sequential mode).
Figure~\ref{fig-map-nloc} then illustrates how these gains in computational
costs are translated in the new 'best method' maps/minimum cost surfaces.

\begin{figure}[htbp]
    \centerline{\fbox{\parbox[c]{7.2cm}{
        \centerline{$p=1$}
        \vspace*{2mm}
        \includegraphics[height=3.4cm,width=3.4cm,clip=true]{map-loc-1-rho.eps}
        \hspace*{2mm}
        \includegraphics[height=3.4cm,width=3.4cm,clip=true]{minc-loc-1-rho.eps}
      }}
    \fbox{\parbox[c]{7.2cm}{
        \centerline{$p=15$}
        \vspace*{2mm}
        \includegraphics[height=3.4cm,width=3.4cm,clip=true]{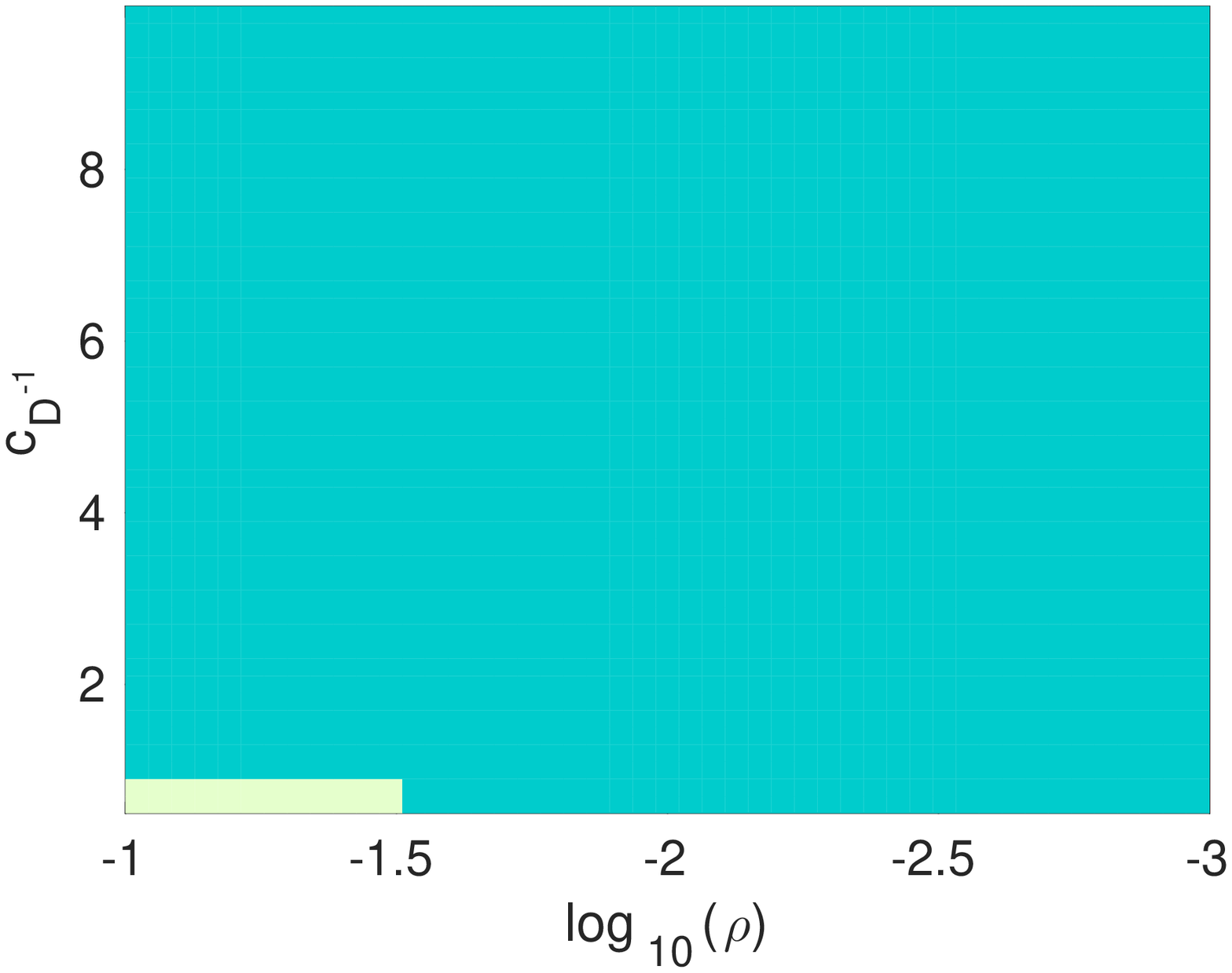} 
        \hspace*{2mm}
        \includegraphics[height=3.4cm,width=3.4cm,clip=true]{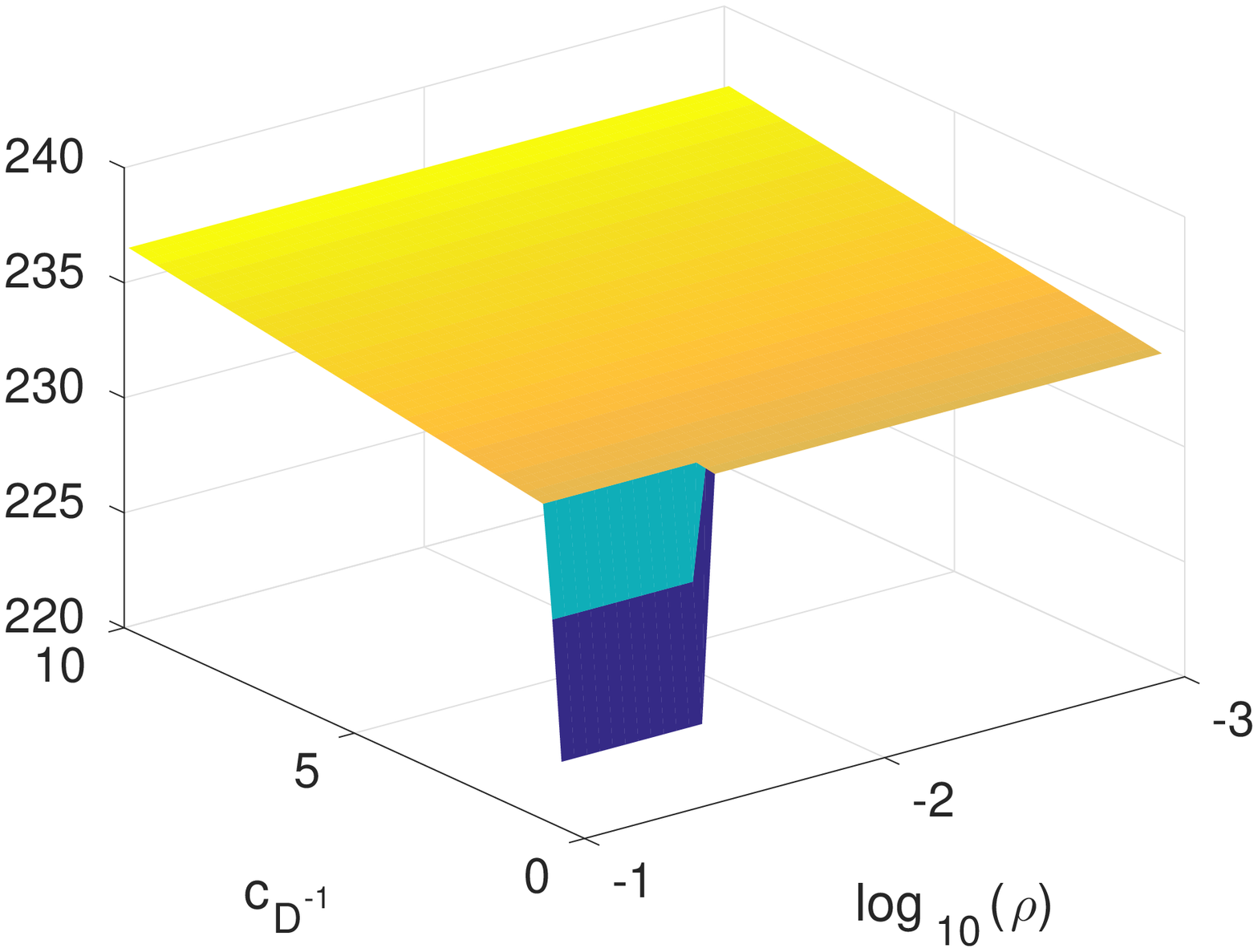}
    }}}
    \vspace*{1mm}
    \centerline{\fbox{\parbox[c]{7.2cm}{
        \centerline{$p=25$}
        \vspace*{2mm}
        \includegraphics[height=3.4cm,width=3.4cm,clip=true]{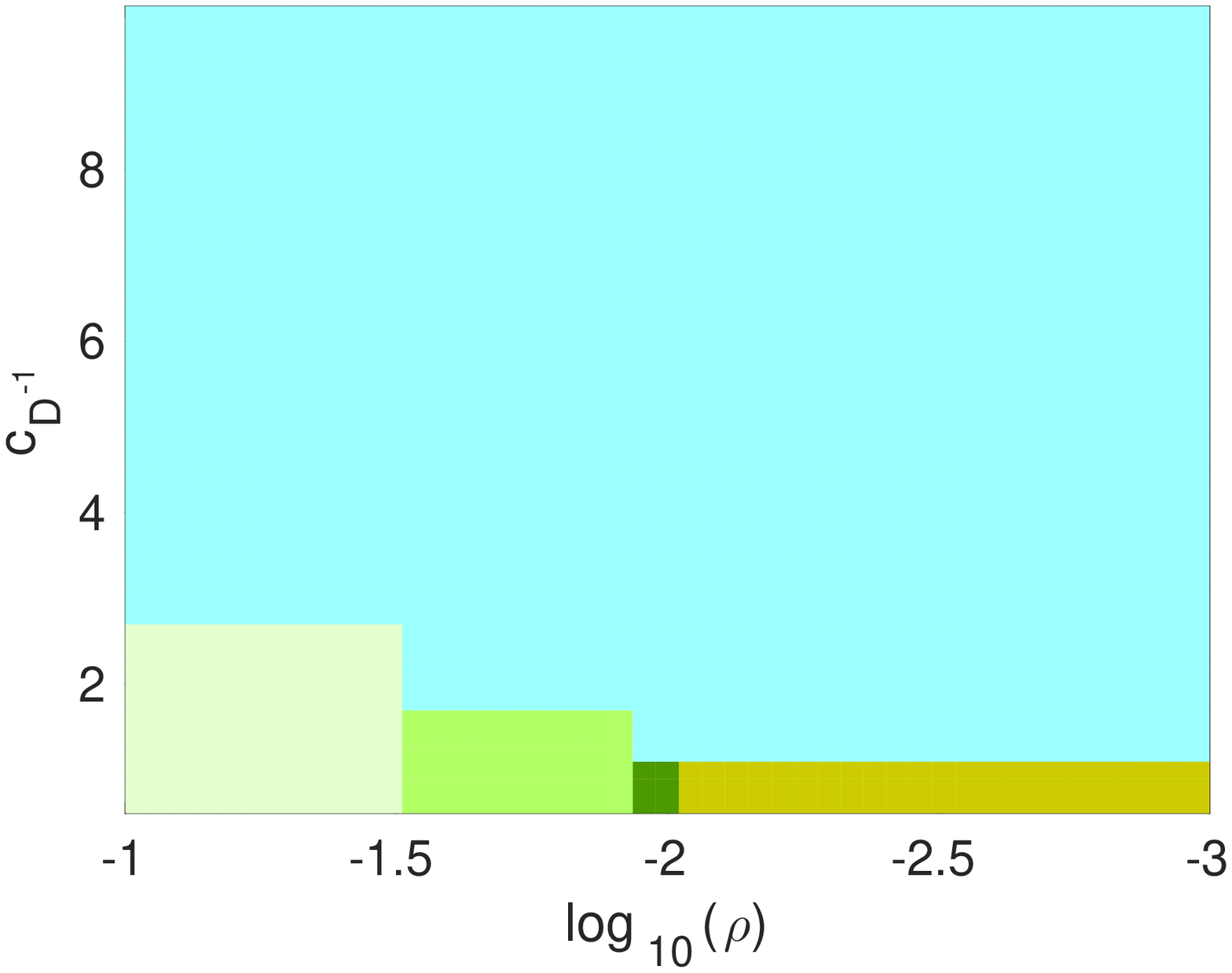} 
        \hspace*{2mm}
        \includegraphics[height=3.4cm,width=3.4cm,clip=true]{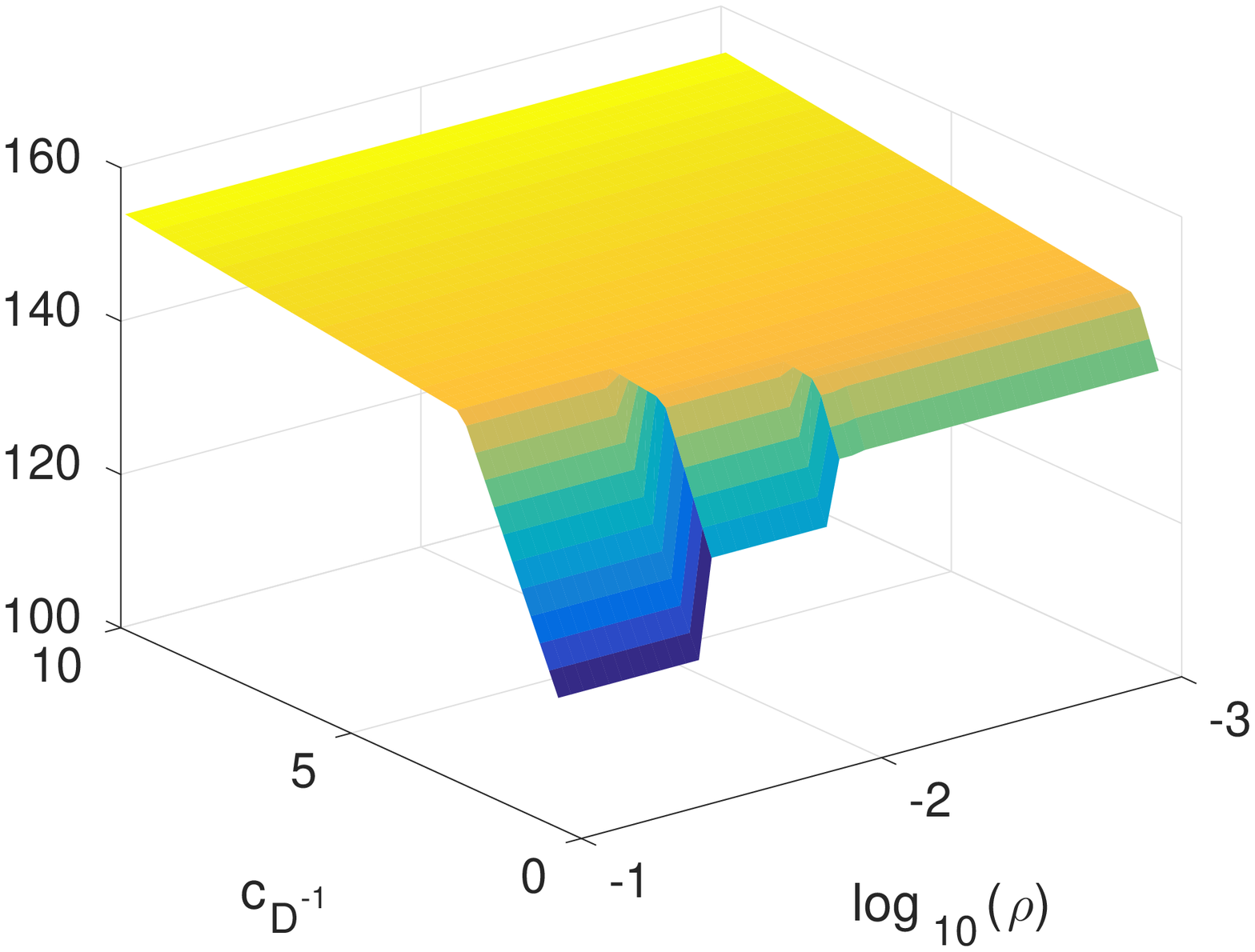}
      }}
    \fbox{\parbox[c]{7.2cm}{
        \centerline{$p=50$}
        \vspace*{2mm}
        \includegraphics[height=3.4cm,width=3.4cm,clip=true]{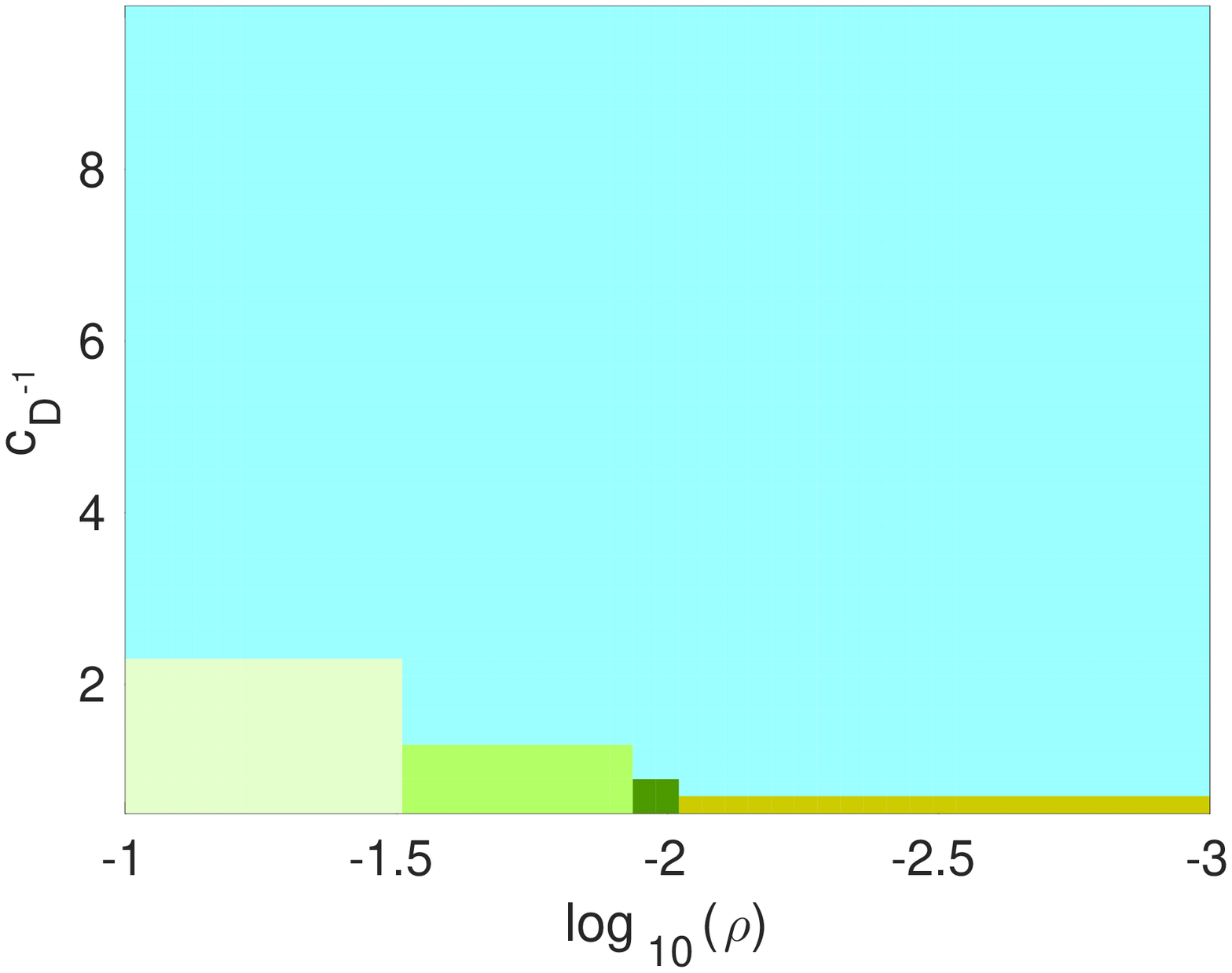} 
        \hspace*{2mm}
        \includegraphics[height=3.4cm,width=3.4cm,clip=true]{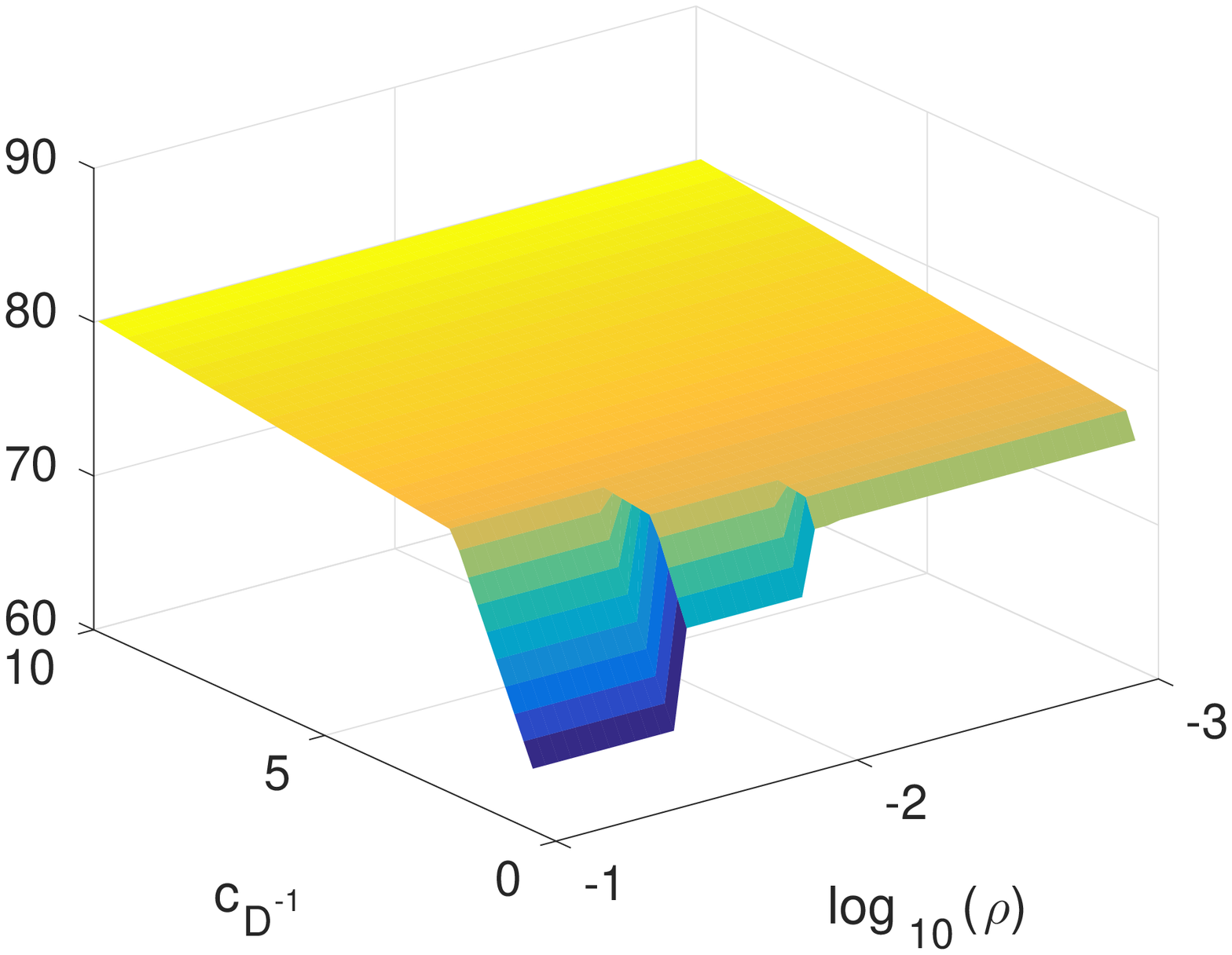}
    }}}
  \vspace*{1mm}
  \begin{center}
   Map colors: 
  \begin{tabular}{llllll}
  \includegraphics[height=3mm,width=3mm]{sa1q50_Mi_color.eps} \hspace*{1mm} & SAQ50-M-I &
  \includegraphics[height=3mm,width=3mm]{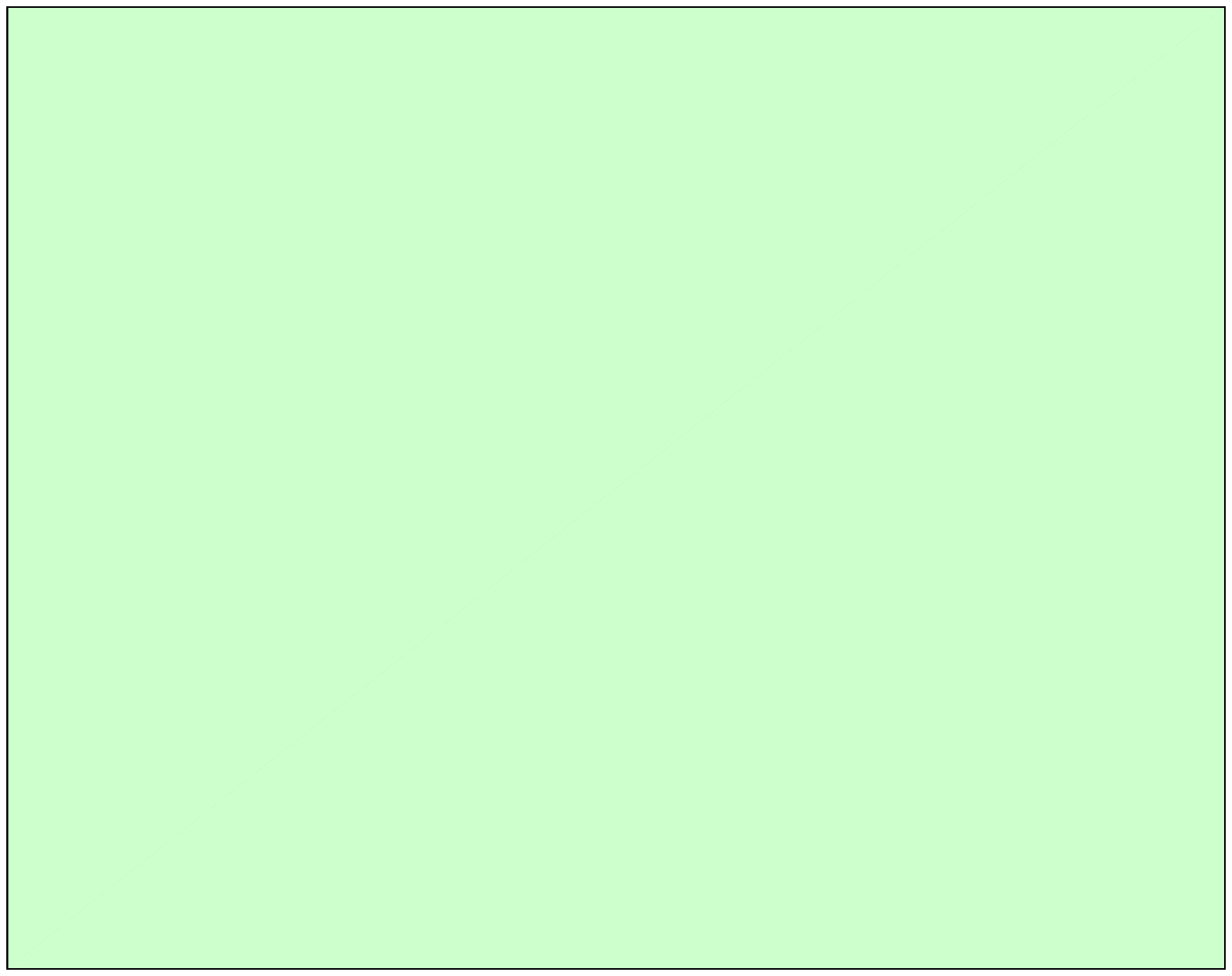}    \hspace*{1mm} & STQ1-n    &
  \includegraphics[height=3mm,width=3mm]{stfq1_Li_color.eps}  \hspace*{1mm} & STQ1-S-0  \\
  \includegraphics[height=3mm,width=3mm]{stfq25_Li_color.eps} \hspace*{1mm} & STQ25-S-0 &
  \includegraphics[height=3mm,width=3mm]{stfq50_Li_color.eps} \hspace*{1mm} & STQ50-S-0 &
  \includegraphics[height=3mm,width=3mm]{sfo2_Li_color.eps}   \hspace*{1mm} & FOQ15-D   \\
  \end{tabular}
  \end{center}
\caption{\label{fig-map-nloc}Best algorithmic variants as a function of
  $c_{D^{-1}}$  the reliability factor $\rho$ and the number of computing
  processes $p$ (Burgers example, hybrid MPI/OpenMP model)}
\end{figure}
        
\noindent
This figure broadly confirms and amplifies the trends already present in
Figure~\ref{fig-map-loc}, with the SAQ50-M-I method becoming more important
especially when $p$ grows.  The reduction of total computational cost is also
noticeable, the smallest cost now being approximately 80 (which corresponds to a
speed-up close to 4).

\numsection{Numerical comparisons on the QG example}\label{numerics-qg-s}

We now turn to the results obtained for the two-layers ECMWF QG example, also
using $n_{\rm inner} = 50$ and at most 5 Gauss-Newton iterations (see
Appendix A2 for a more complete description of the problem).
We first observe in Figure~\ref{fig-QG-forecast} that the fields of interest
do evolve (relatively) slowly over the complete assimilation time, and thus
even more so within each of the 48 subwindows considered.  Hence we may expect
$\tilde{M}_i = I$ to be a reasonable approximation of $M_i$, at variance with
the Burgers case (see Figure~A.14).  In particular, the caveat on
using STQ$\ell$-S-I may no longer apply.

\begin{figure}[htbp]
\centerline{
  \includegraphics[height=9cm,width=10cm]{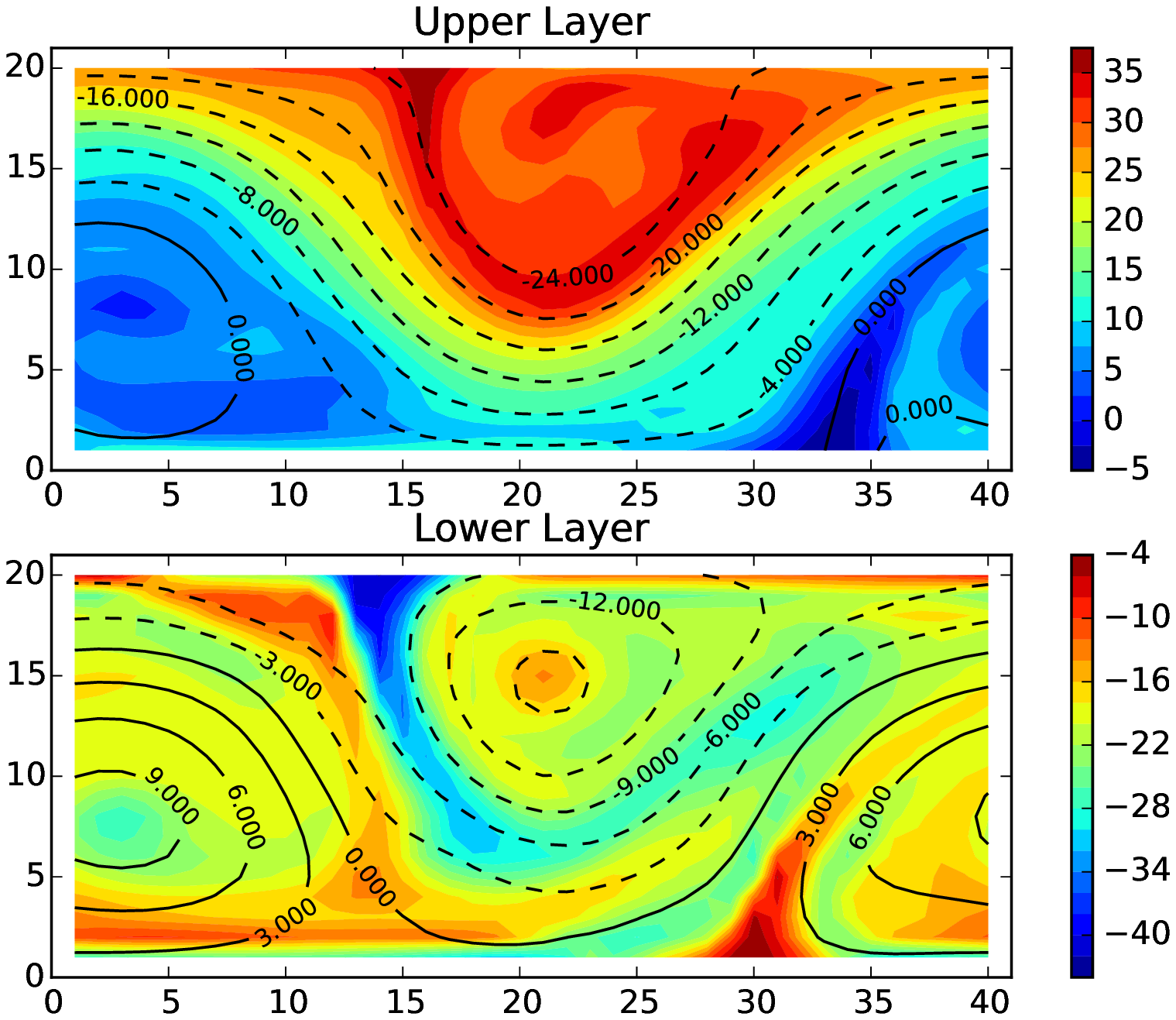}
  \hspace*{-25mm}
  \includegraphics[height=9cm,width=10cm]{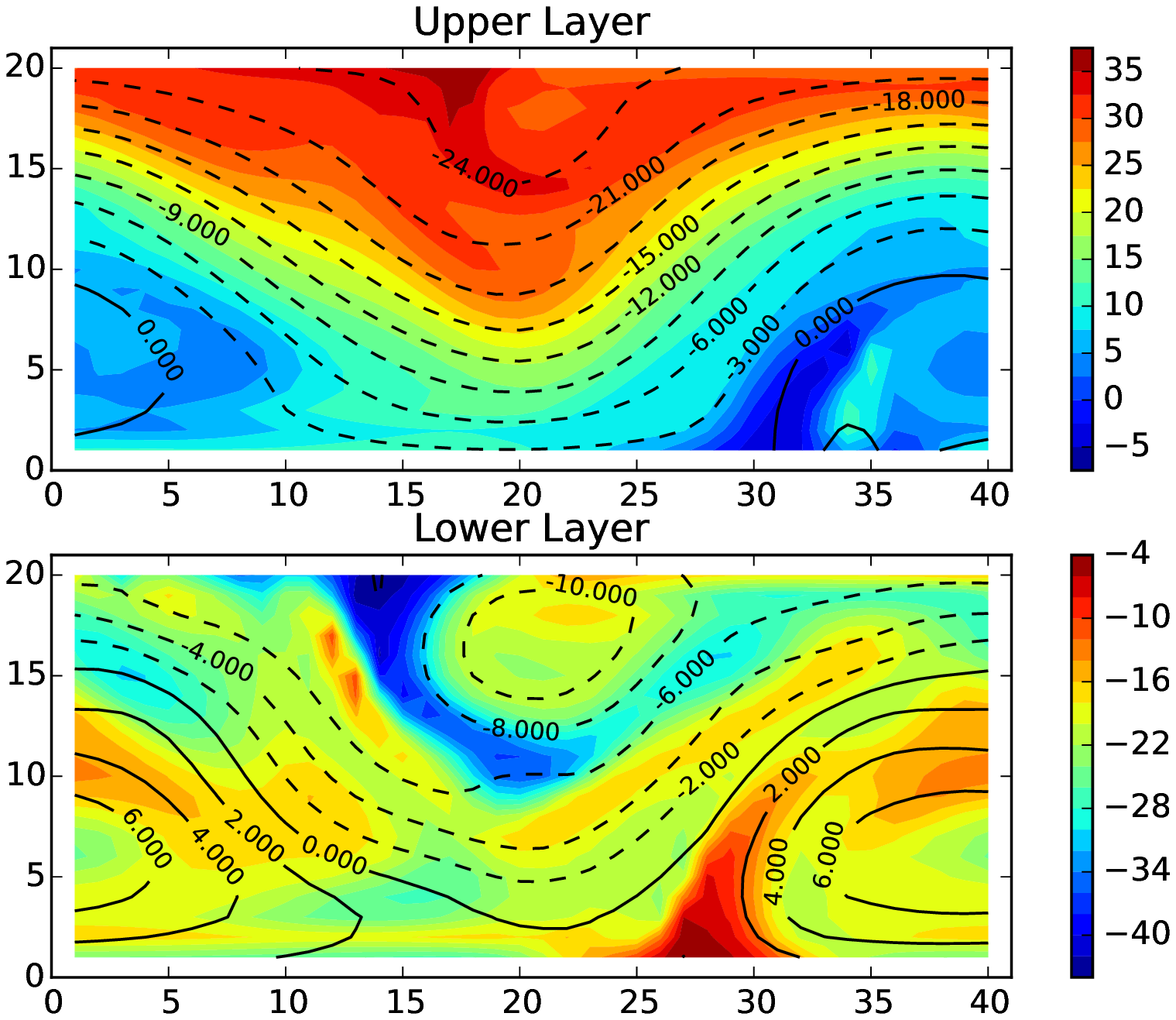}
}
\caption{\label{fig-QG-forecast}Initial (left) and final (right) forecasts using the QG model} 
\end{figure}

Keeping this in mind, we next examine the performance of the original saddle
methods SAQ0-M and SAQ0-M and compare them, first in terms of number of inner
iterations, for the parellizable preconditioners \req{tildeM}, to competing
variants such as SAQ15 or STQ15. The outcome is presented in
Figure~\ref{fig-best-Li} for the choice $\tilde{M}_i = 0$ and
Figure~\ref{fig-best-Mi} for the choice $\tilde{M}_i = I$ .

\begin{figure}[htbp]
\centerline{
  \includegraphics[height=6cm,width=14cm]{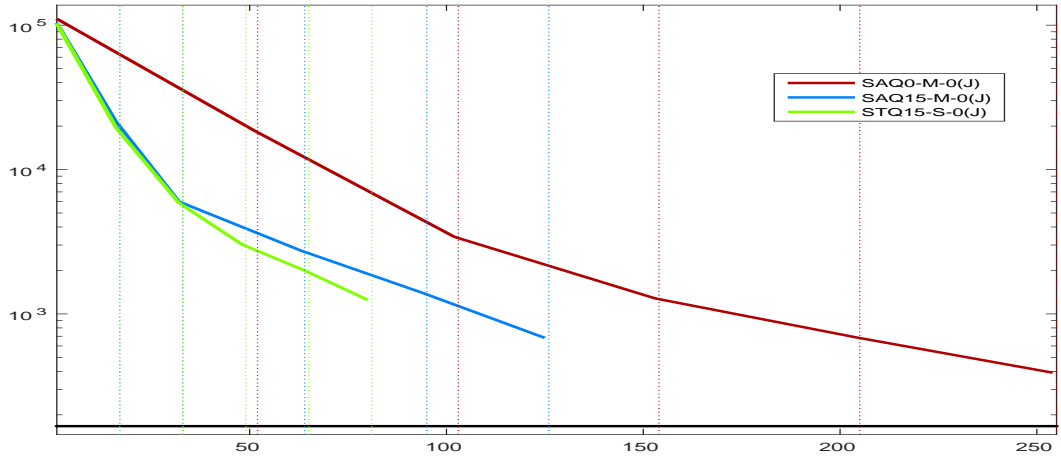} 
}
\caption{\label{fig-best-Li} Evolution of $q$ (dashed) and $J$
 (continuous) as a function of the number of inner iterations for
  SAQ0-M-0, SAQ15-M-0 and STQ15-S-0 (QG example)}
\end{figure}

\begin{figure}[htbp]
\vspace*{2mm}
\centerline{
  \includegraphics[height=6cm,width=14cm]{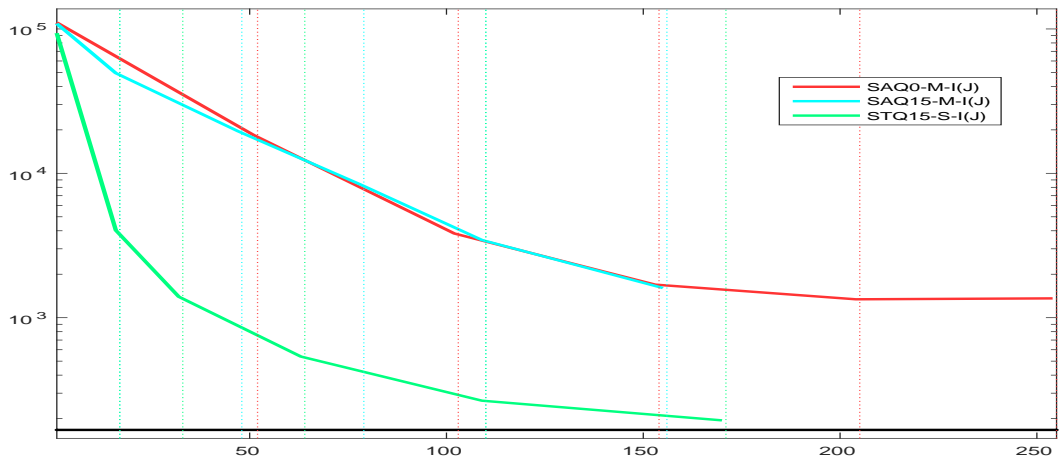} 
}
\caption{\label{fig-best-Mi} Evolution of $q$ (dashed) and $J$
  (continuous) as a function of the number of inner iterations, for 
  SAQ0-M-I, SAQ15-M-I and STQ15-S-I (QG example)}
\end{figure}

\noindent
We see in the first of these figures that the simply preconditioned SAQ0-M-0
performs relatively well, albeit a bit slowly compared to the globalized
saddle SAQ15-M-0 and the state algorithm STQ15-S-0.  The dominance of this
latter formulation with STQ15-S-I is even more obvious in the second figure,
SAQ15-M-I and SAQ0-M-I following the same initial curve, the latter then
levelling off and diverging slowly.

We may now again compute the same 'best method' maps and minimum cost
surfaces, mimicking our analysis for the Burgers case, yielding
Figures~\ref{fig-map-loc-oops} for the fully MPI
model, and Figure~\ref{fig-map-nloc-oops} for the hybrid MPI/OPenMP one.

\begin{figure}[htbp]
    \vspace*{-1mm}
    \centerline{\fbox{\parbox[c]{7.2cm}{
        \centerline{$p=1$}
        \vspace*{2mm}
        \includegraphics[height=3.4cm,width=3.4cm,clip=true]{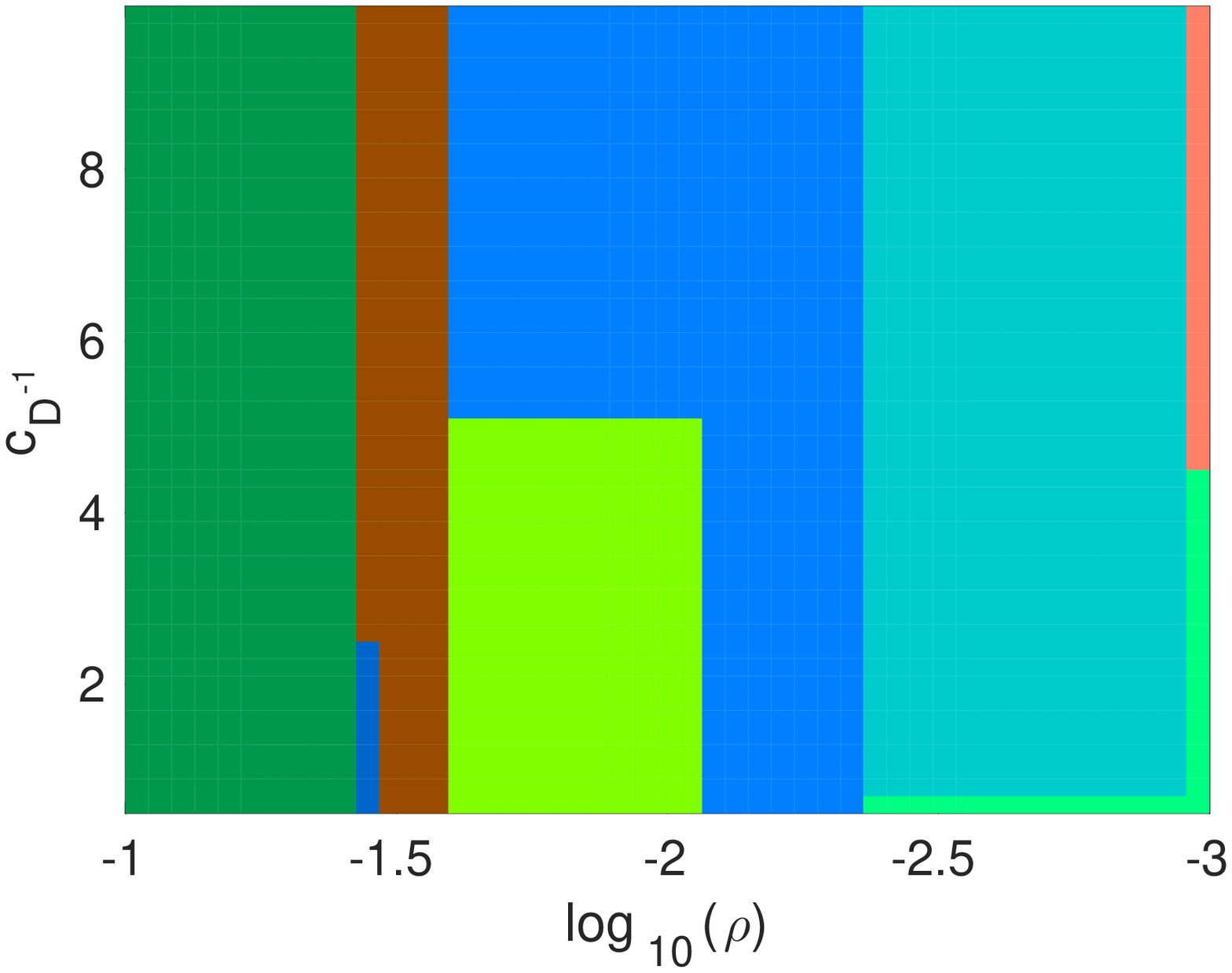}
        \hspace*{2mm}
        \includegraphics[height=3.4cm,width=3.4cm,clip=true]{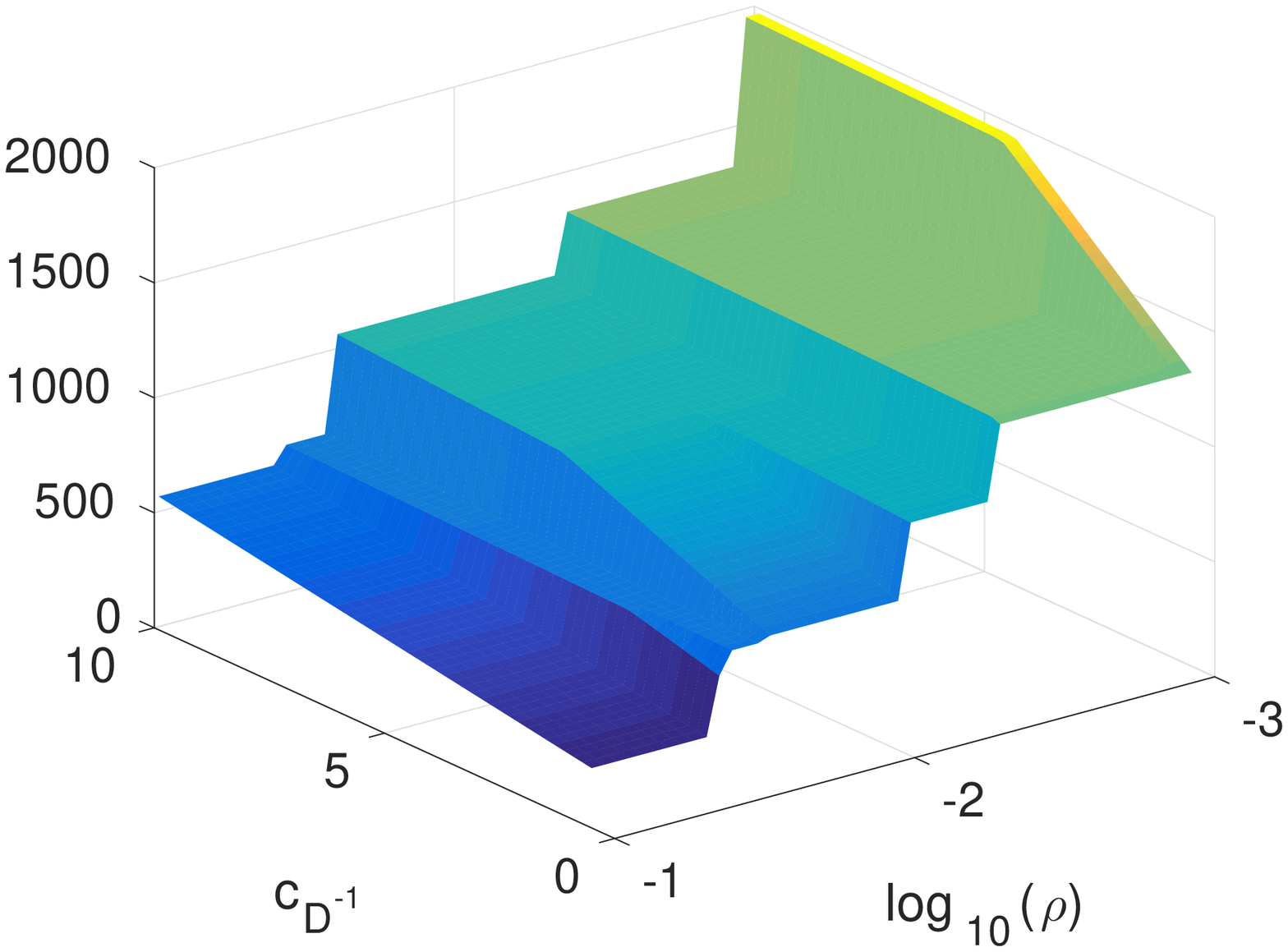}
      }}
    \fbox{\parbox[c]{7.2cm}{
        \centerline{$p=15$}
        \vspace*{2mm}
        \includegraphics[height=3.4cm,width=3.4cm,clip=true]{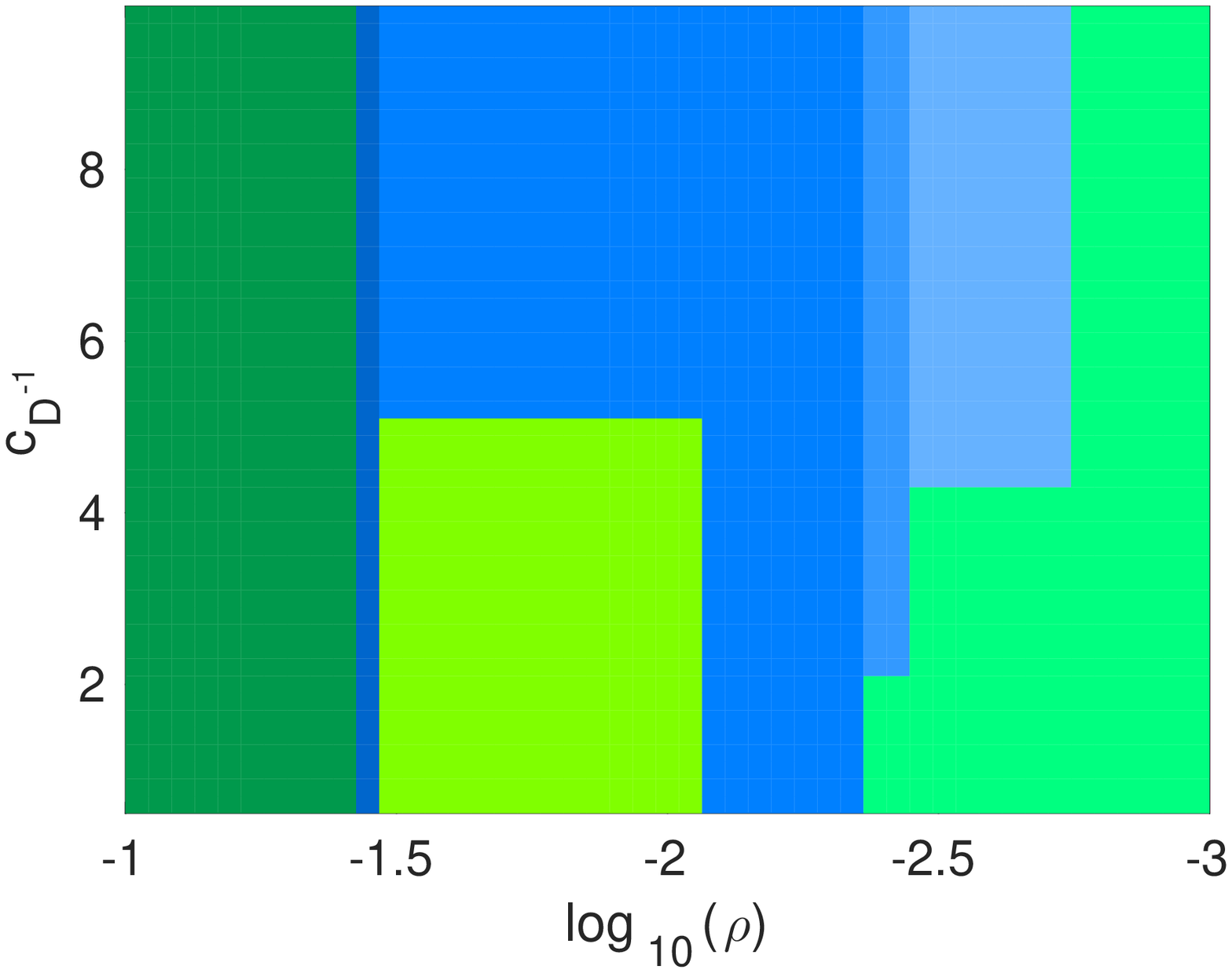} 
        \hspace*{2mm}
        \includegraphics[height=3.4cm,width=3.4cm,clip=true]{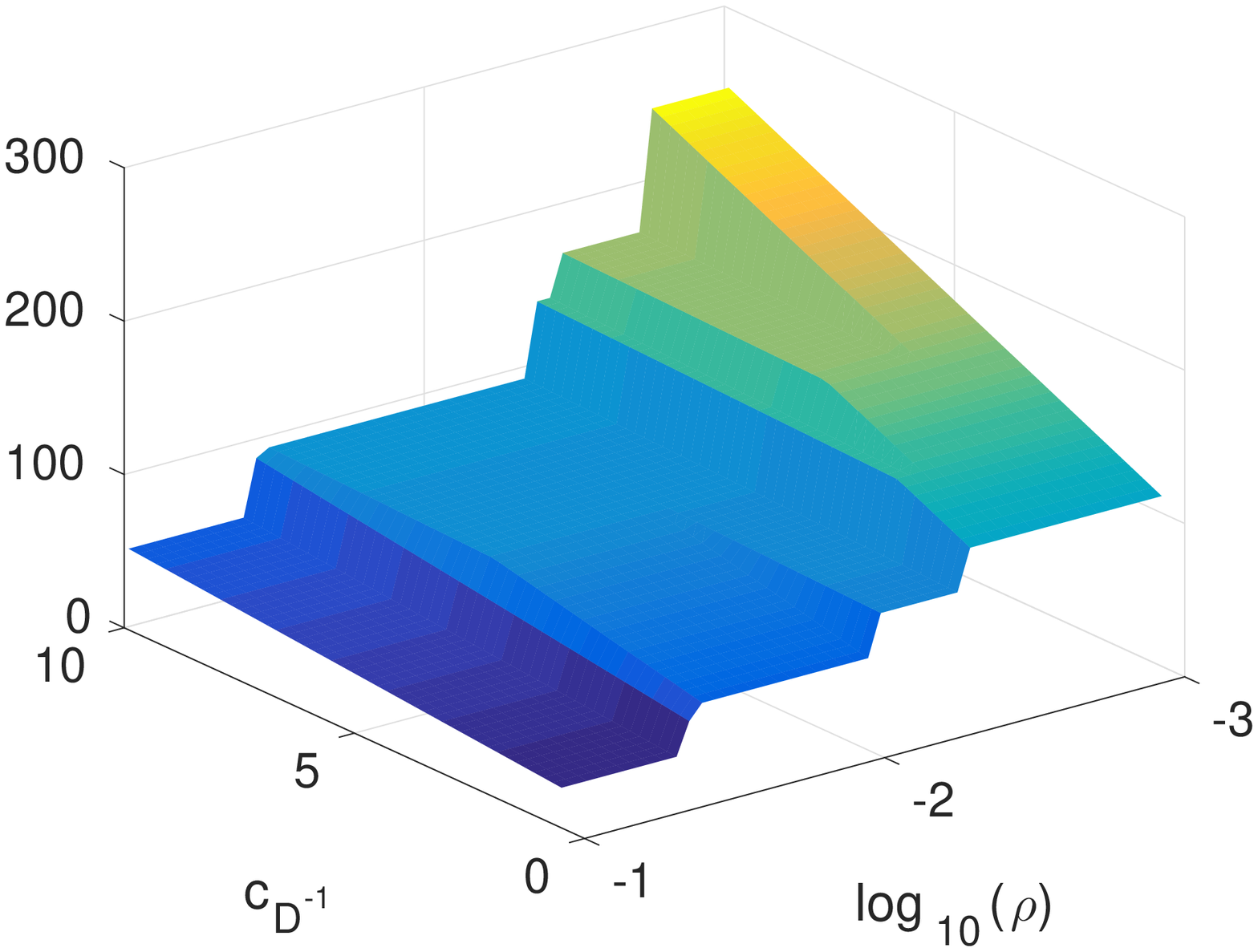}
      }}}
    \vspace*{1mm}
    \centerline{\fbox{\parbox[c]{7.2cm}{
        \centerline{$p=25$}
        \vspace*{2mm}
        \includegraphics[height=3.4cm,width=3.4cm,clip=true]{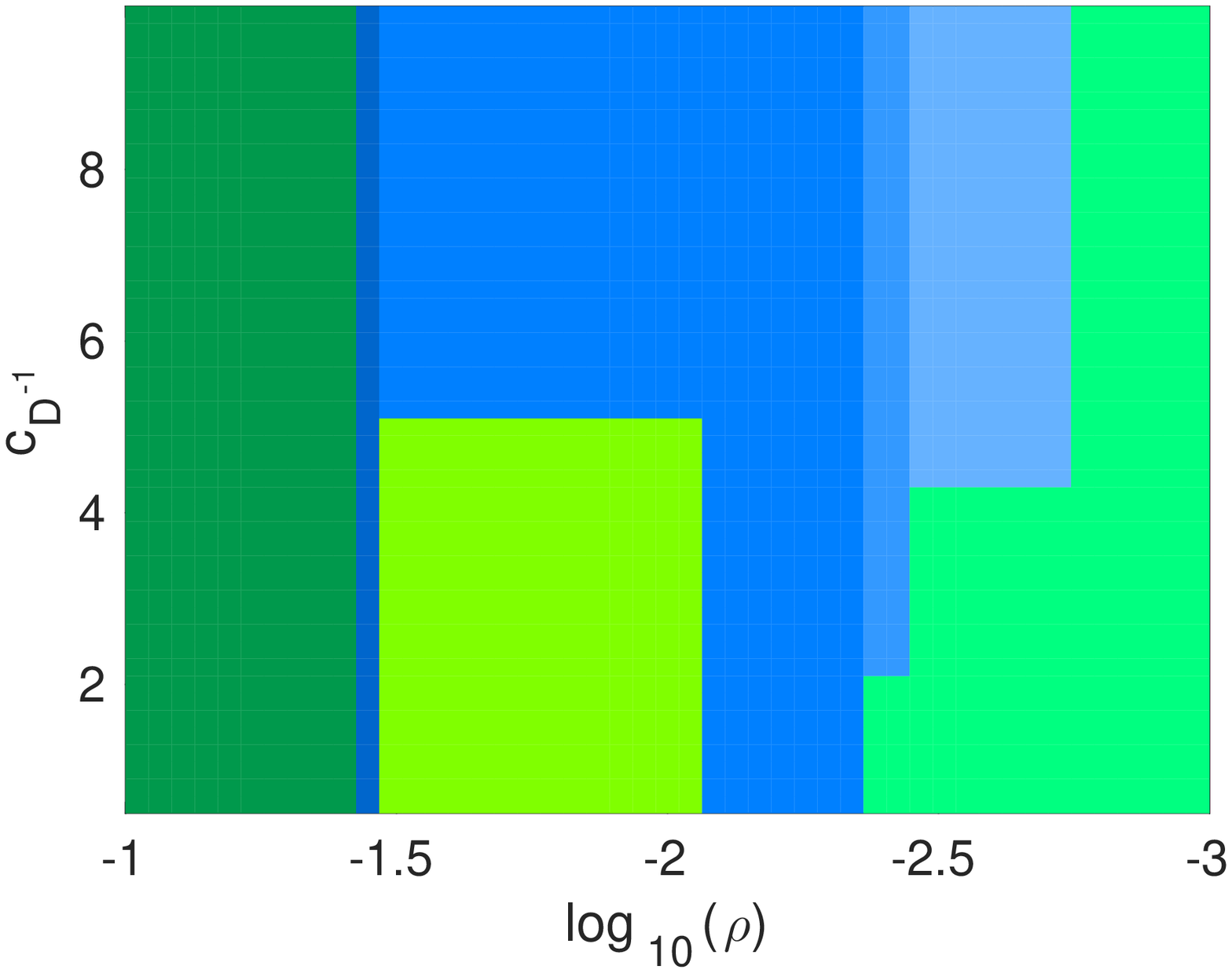} 
        \hspace*{2mm}
        \includegraphics[height=3.4cm,width=3.4cm,clip=true]{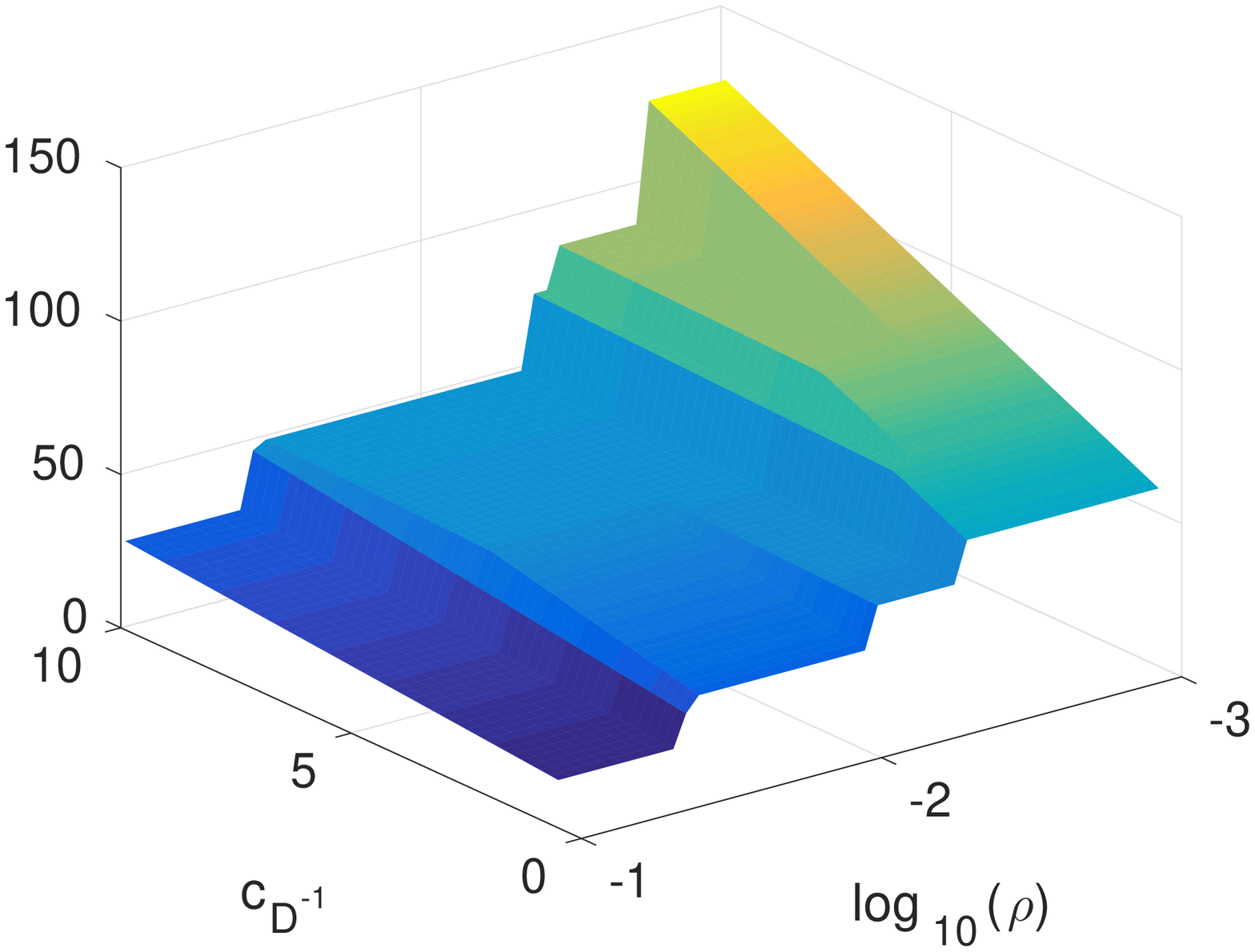}
      }}
    \fbox{\parbox[c]{7.2cm}{
        \centerline{$p=50$}
        \vspace*{2mm}
        \includegraphics[height=3.4cm,width=3.4cm,clip=true]{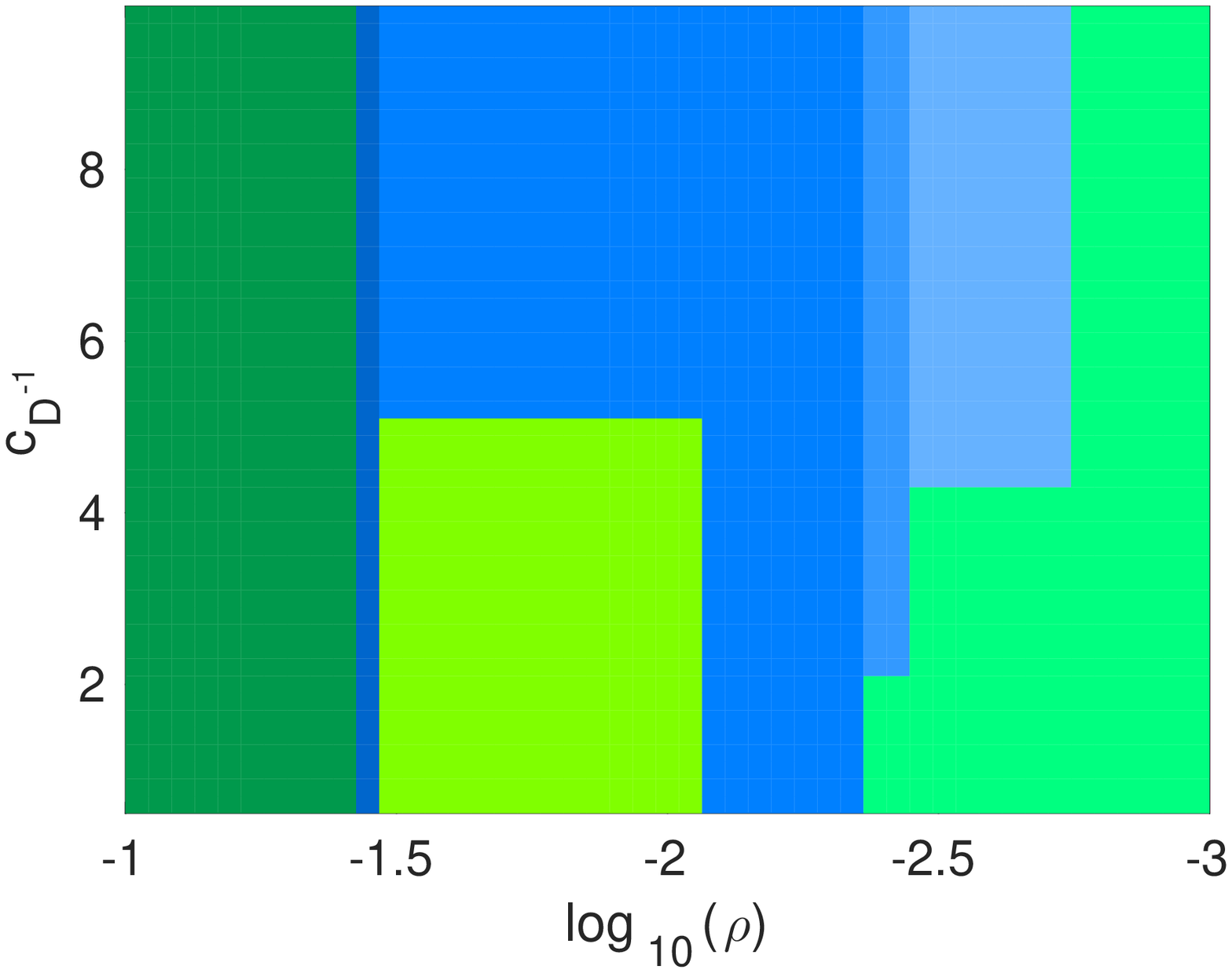} 
        \hspace*{2mm}
        \includegraphics[height=3.4cm,width=3.4cm,clip=true]{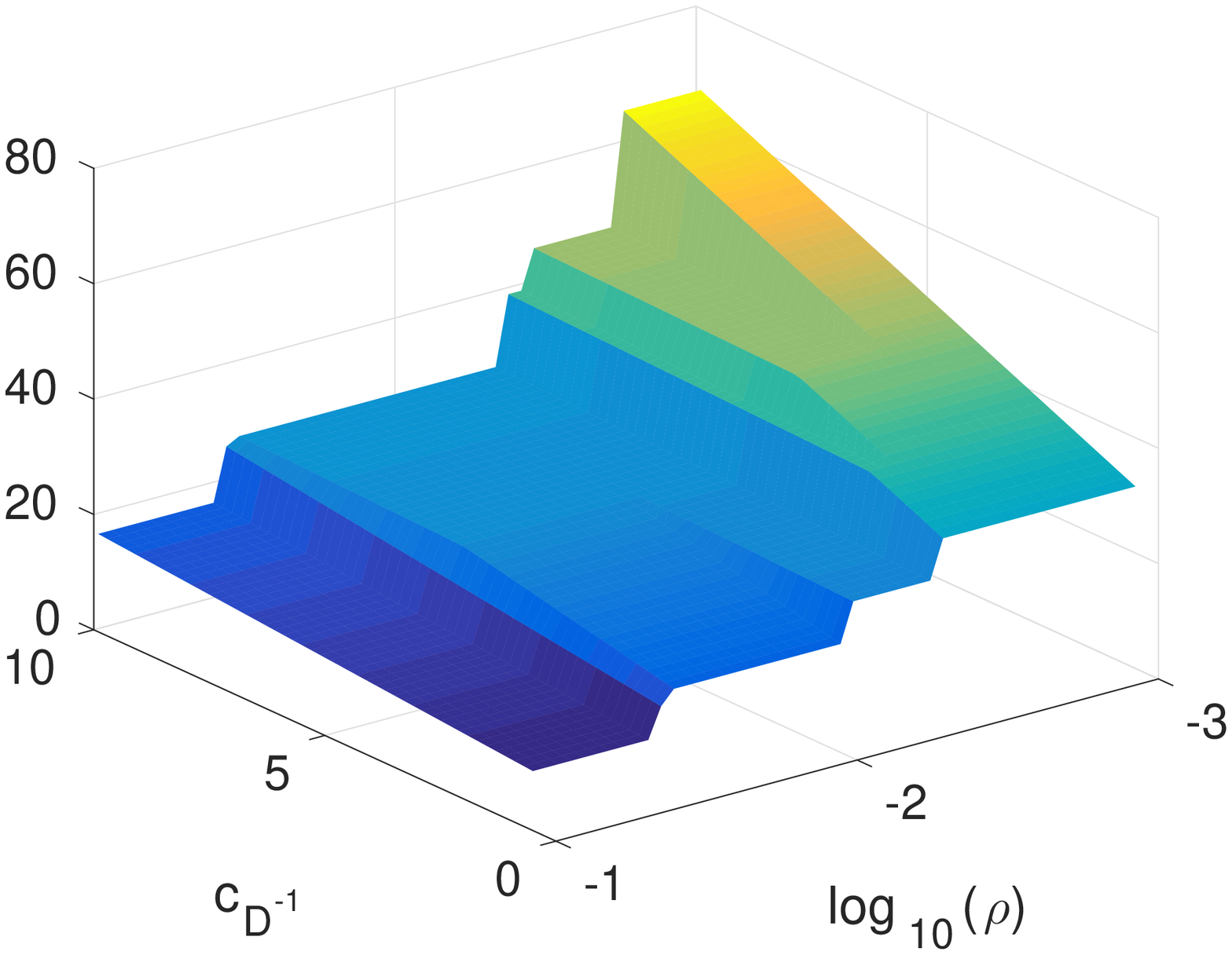}
      }}}
  \vspace*{1mm}
  \begin{center}
   Map colors: 
  \begin{tabular}{llllllll}
  \includegraphics[height=3mm,width=3mm]{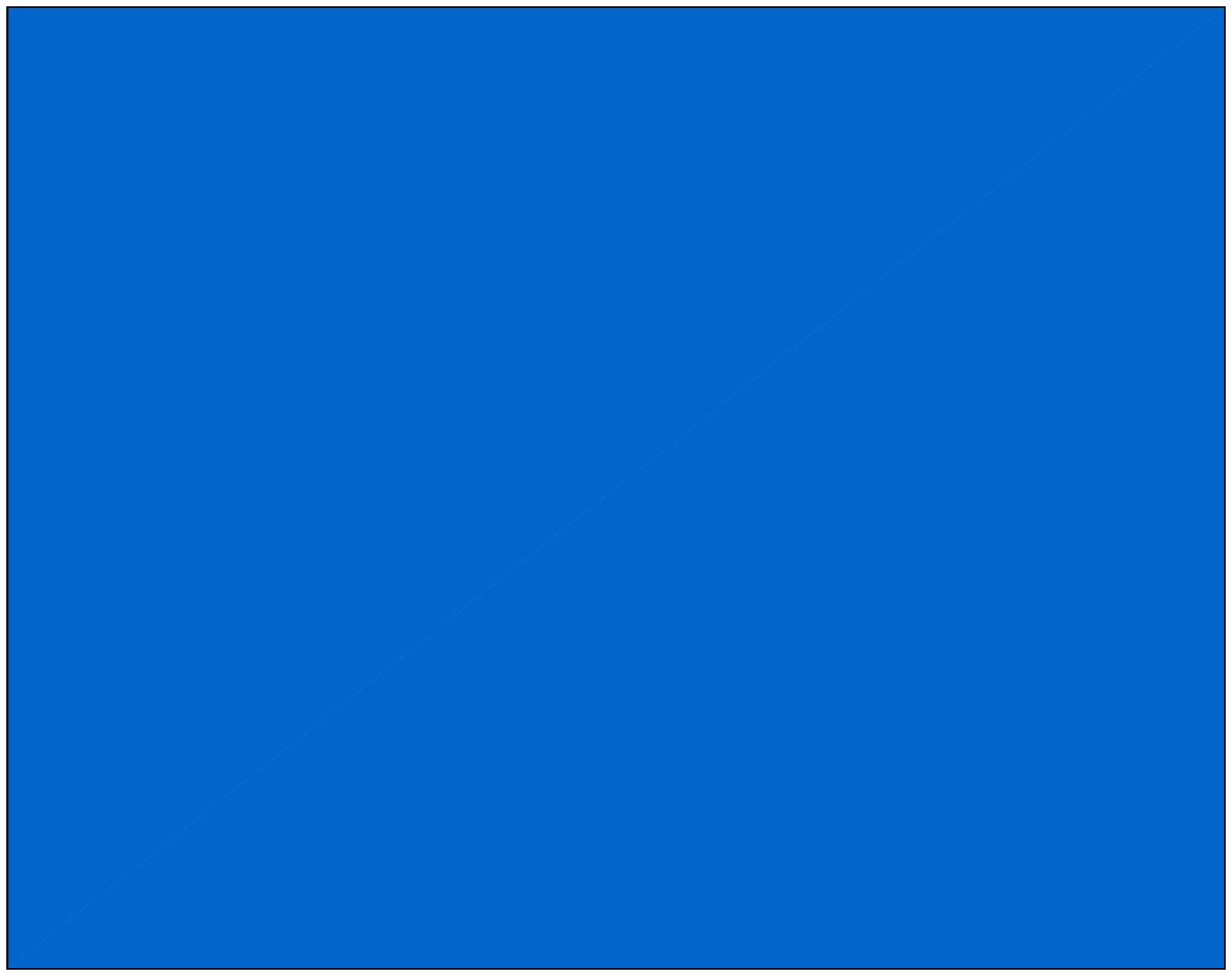}  \hspace*{1mm} & SAQ1-M-0  &
  \includegraphics[height=3mm,width=3mm]{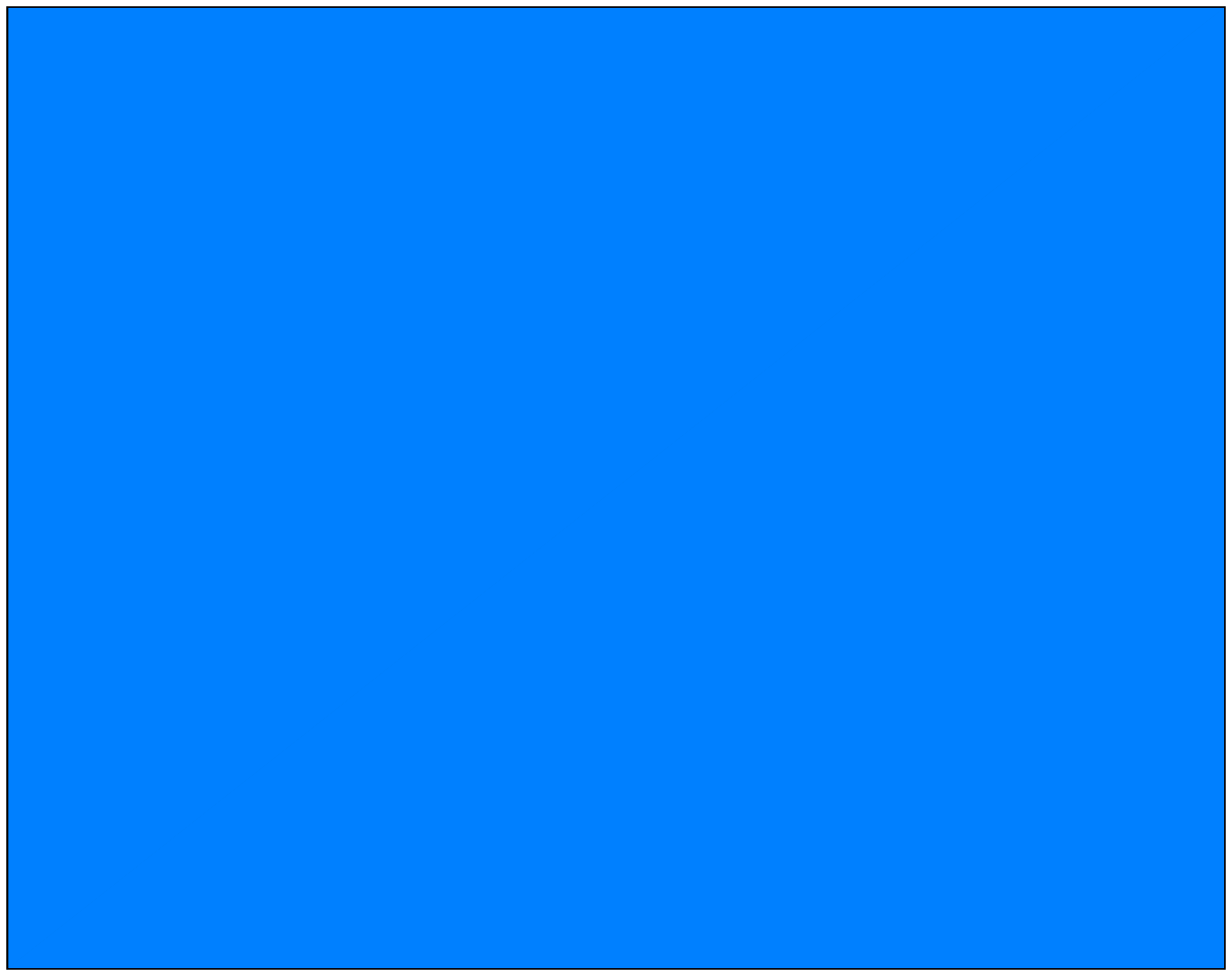} \hspace*{1mm} & SAQ1-M-0  &
  \includegraphics[height=3mm,width=3mm]{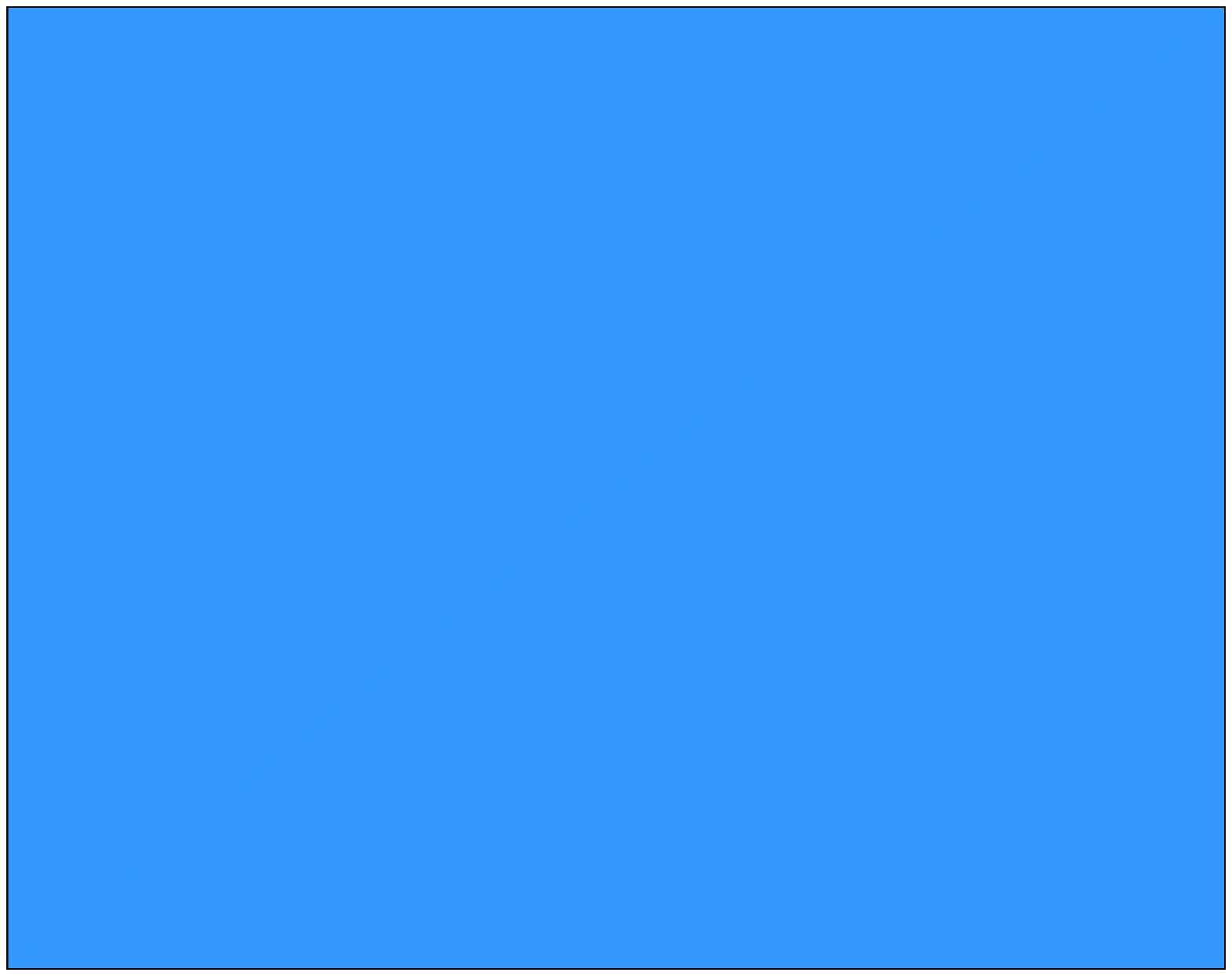} \hspace*{1mm} & SAQ25-M-0 &
  \includegraphics[height=3mm,width=3mm]{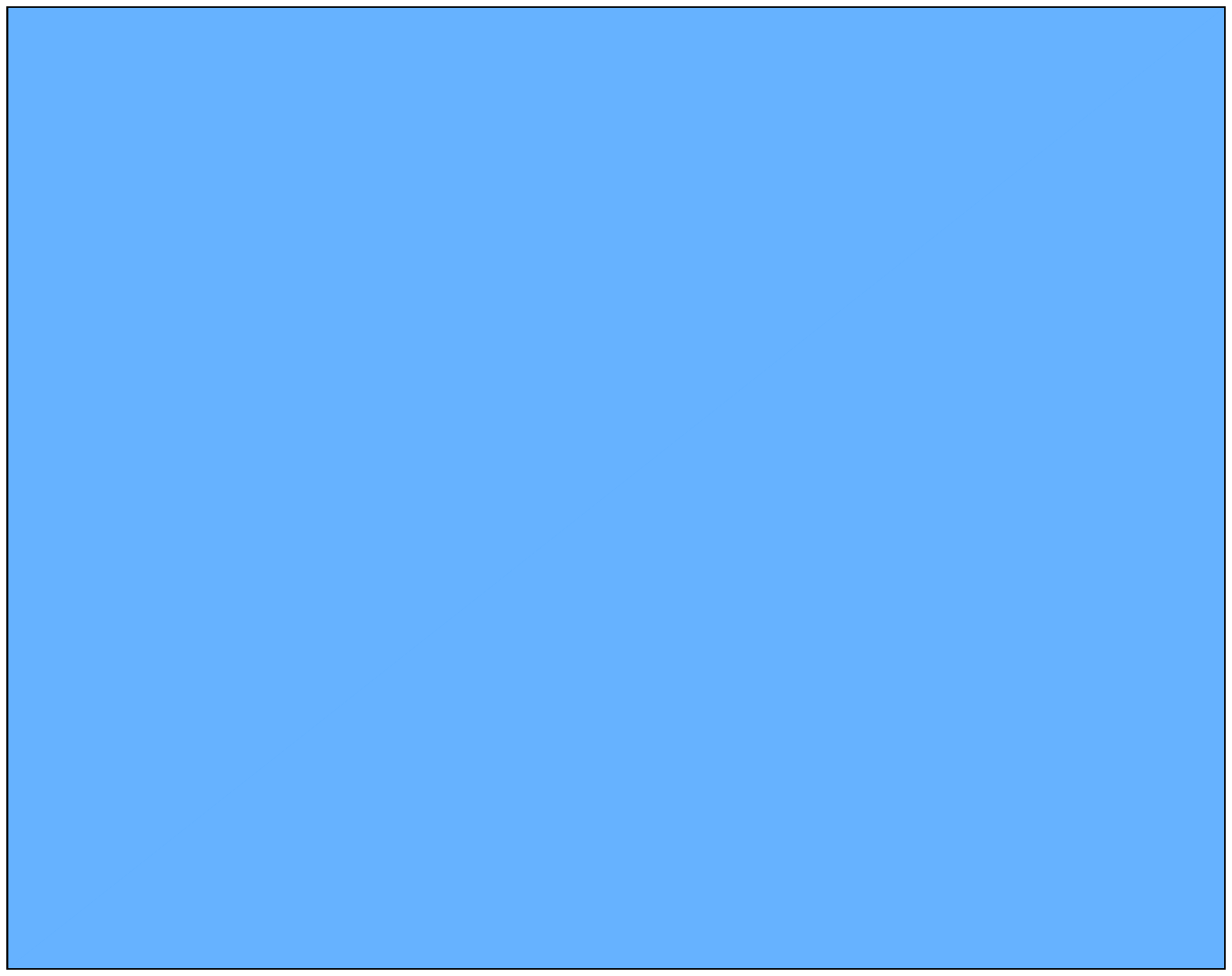} \hspace*{1mm} & SAQ40-M-0 \\
  \includegraphics[height=3mm,width=3mm]{sa1q50_Mi_color.eps}   \hspace*{1mm} & SAQ50-M-I &
  \includegraphics[height=3mm,width=3mm]{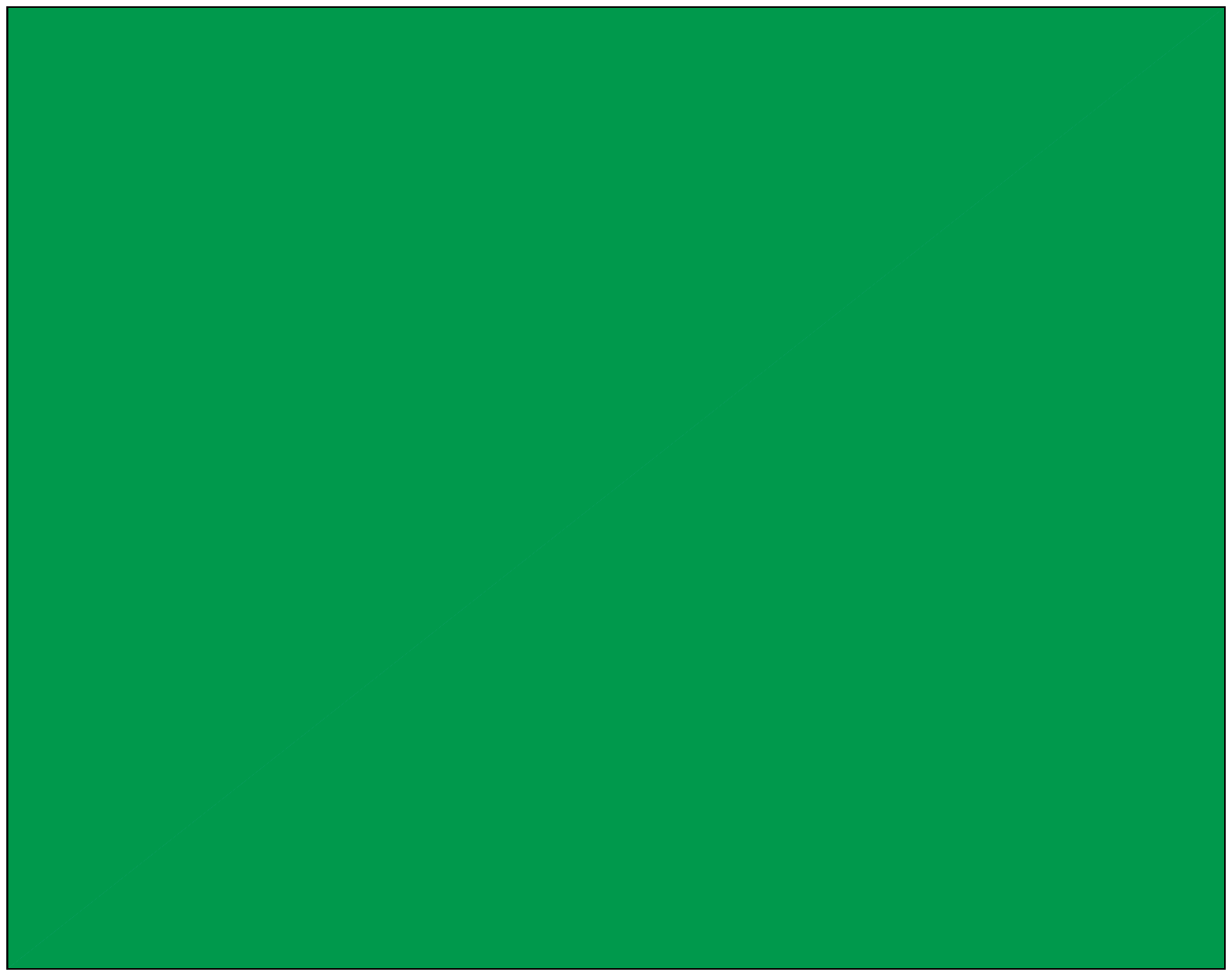}     \hspace*{1mm} & STQ1-S-I  &
  \includegraphics[height=3mm,width=3mm]{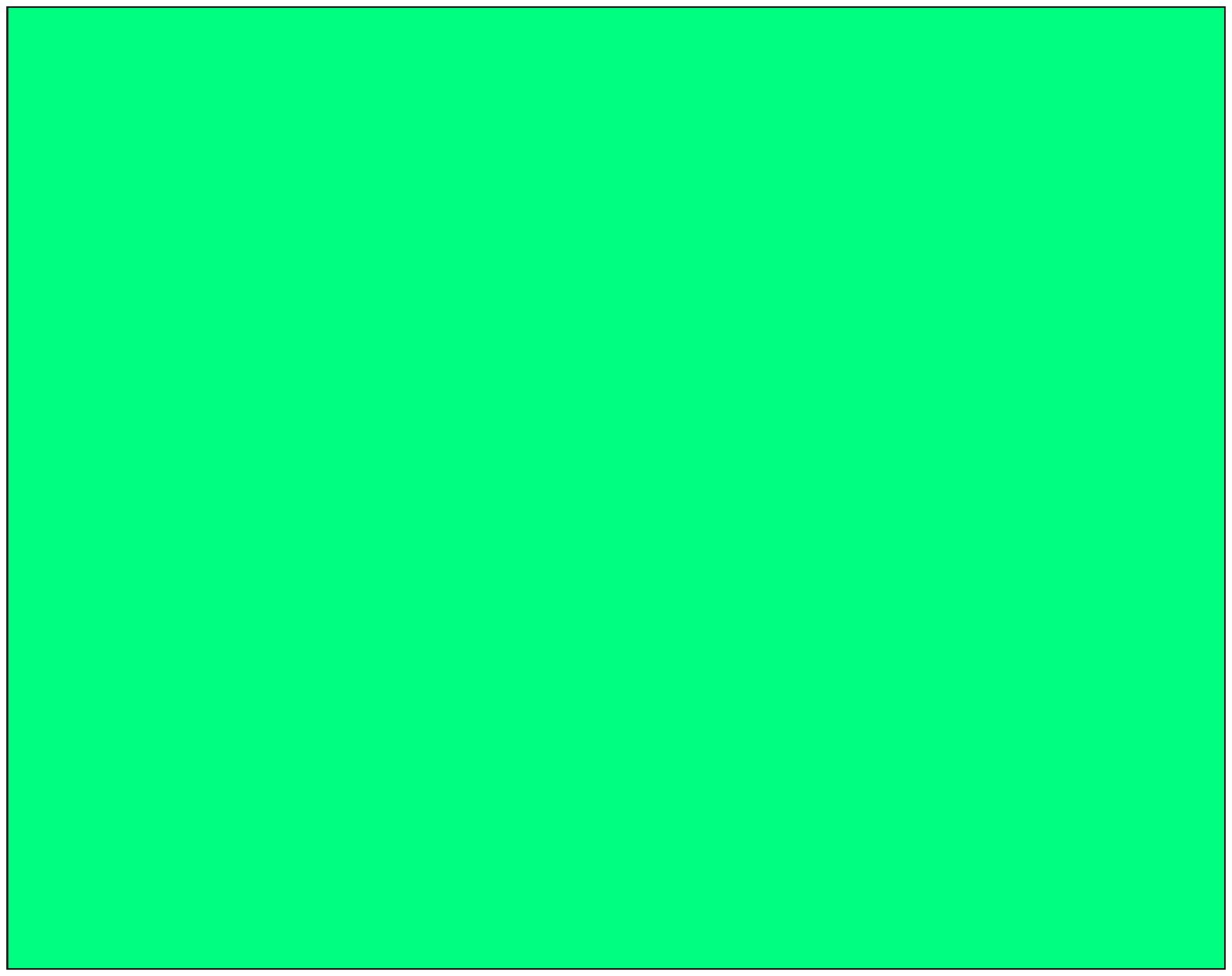}    \hspace*{1mm} & STQ15-S-I &
  \includegraphics[height=3mm,width=3mm]{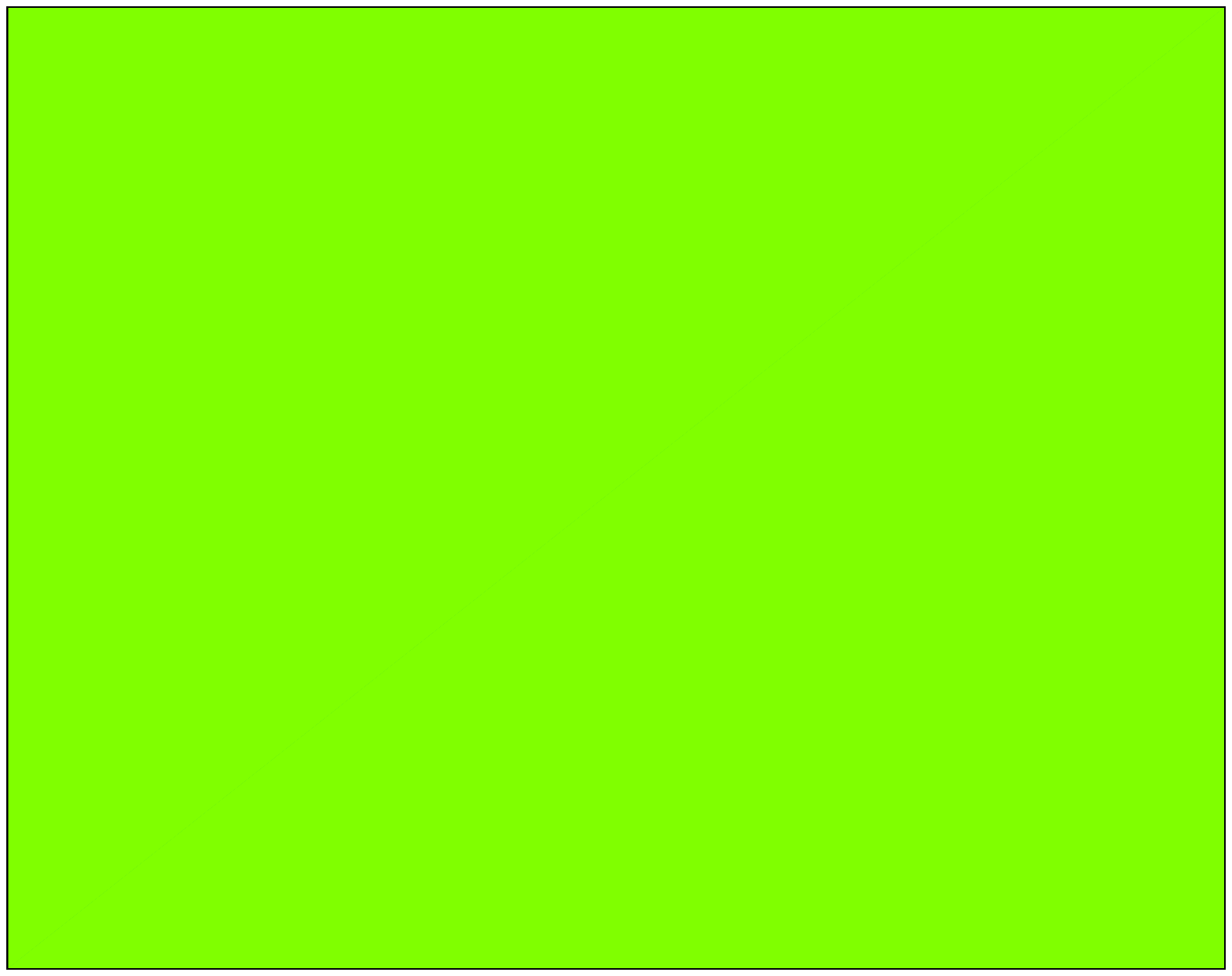}    \hspace*{1mm} & STQ15-S-0 \\
  \includegraphics[height=3mm,width=3mm]{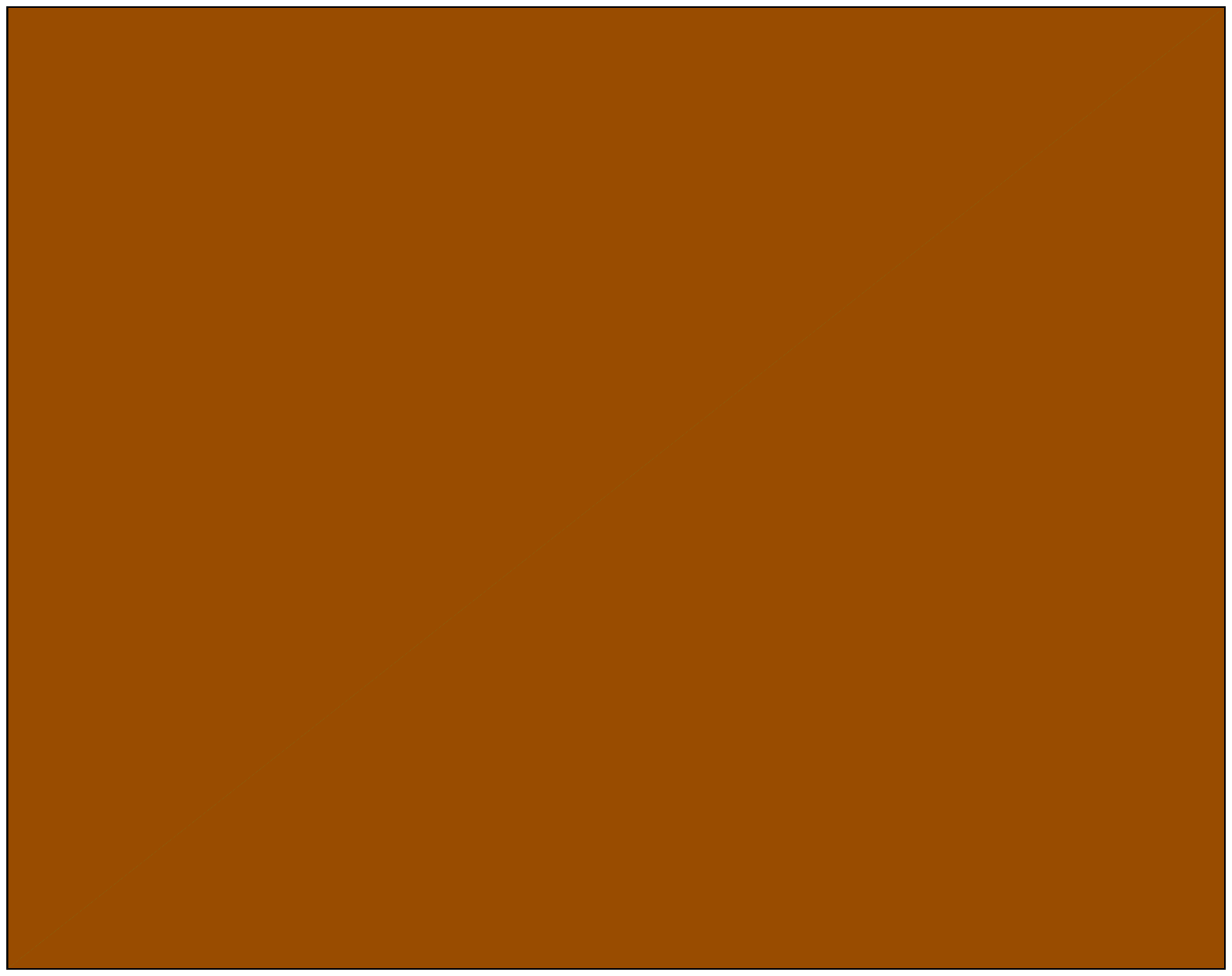}   \hspace*{1mm} & FOQ1-D  &
  \includegraphics[height=3mm,width=3mm]{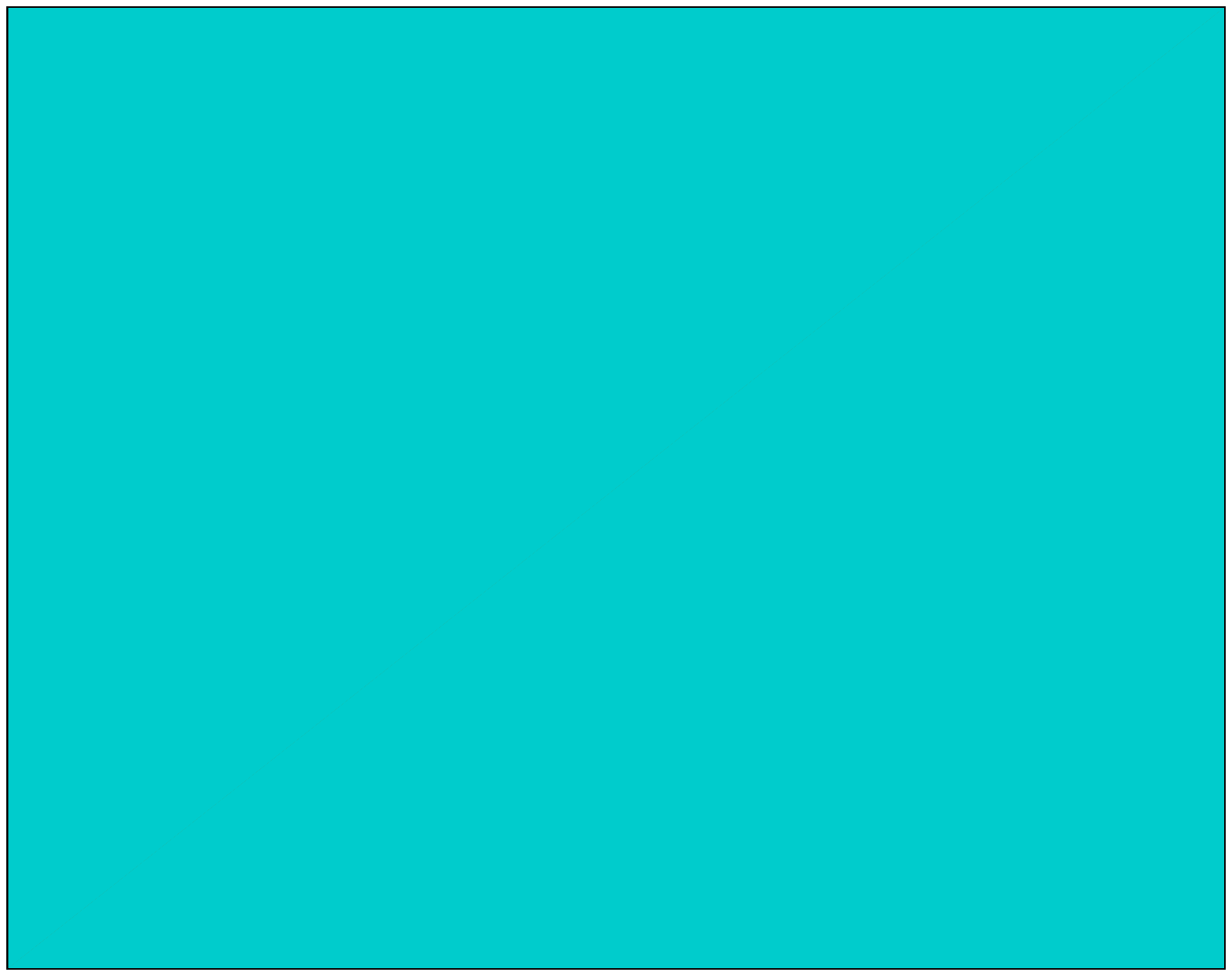}  \hspace*{1mm} & FOQ15-D &
  \includegraphics[height=3mm,width=3mm]{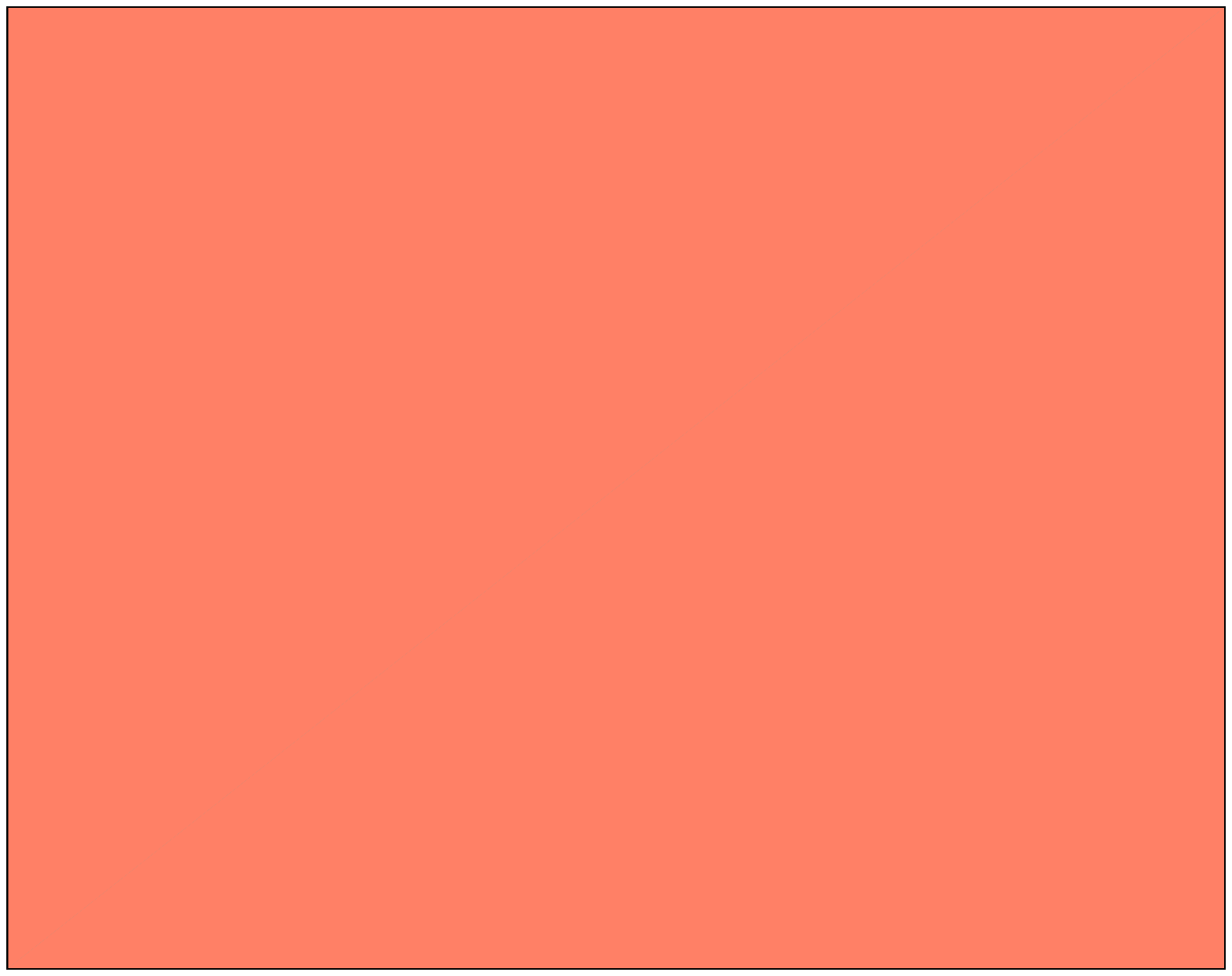}  \hspace*{1mm} & FOQ25-D & \\
  \end{tabular}
  \end{center}
  \vspace*{-4mm}
\caption{\label{fig-map-loc-oops}Best algorithmic variants as a function of
  $c_{D^{-1}}$  the number of computing processes $p$ and the reliability
  factor $\rho$ (QG example, fully MPI model)}
\end{figure}
        
\begin{figure}[htbp]
    \centerline{\fbox{\parbox[c]{7.2cm}{
        \centerline{$p=1$}
        \vspace*{2mm}
        \includegraphics[height=3.4cm,width=3.4cm,clip=true]{map-loc-1-rho-oops48.eps}
        \hspace*{2mm}
        \includegraphics[height=3.4cm,width=3.4cm,clip=true]{minc-loc-1-rho-oops48.eps}
      }}
    \fbox{\parbox[c]{7.2cm}{
        \centerline{$p=15$}
        \vspace*{2mm}
        \includegraphics[height=3.4cm,width=3.4cm,clip=true]{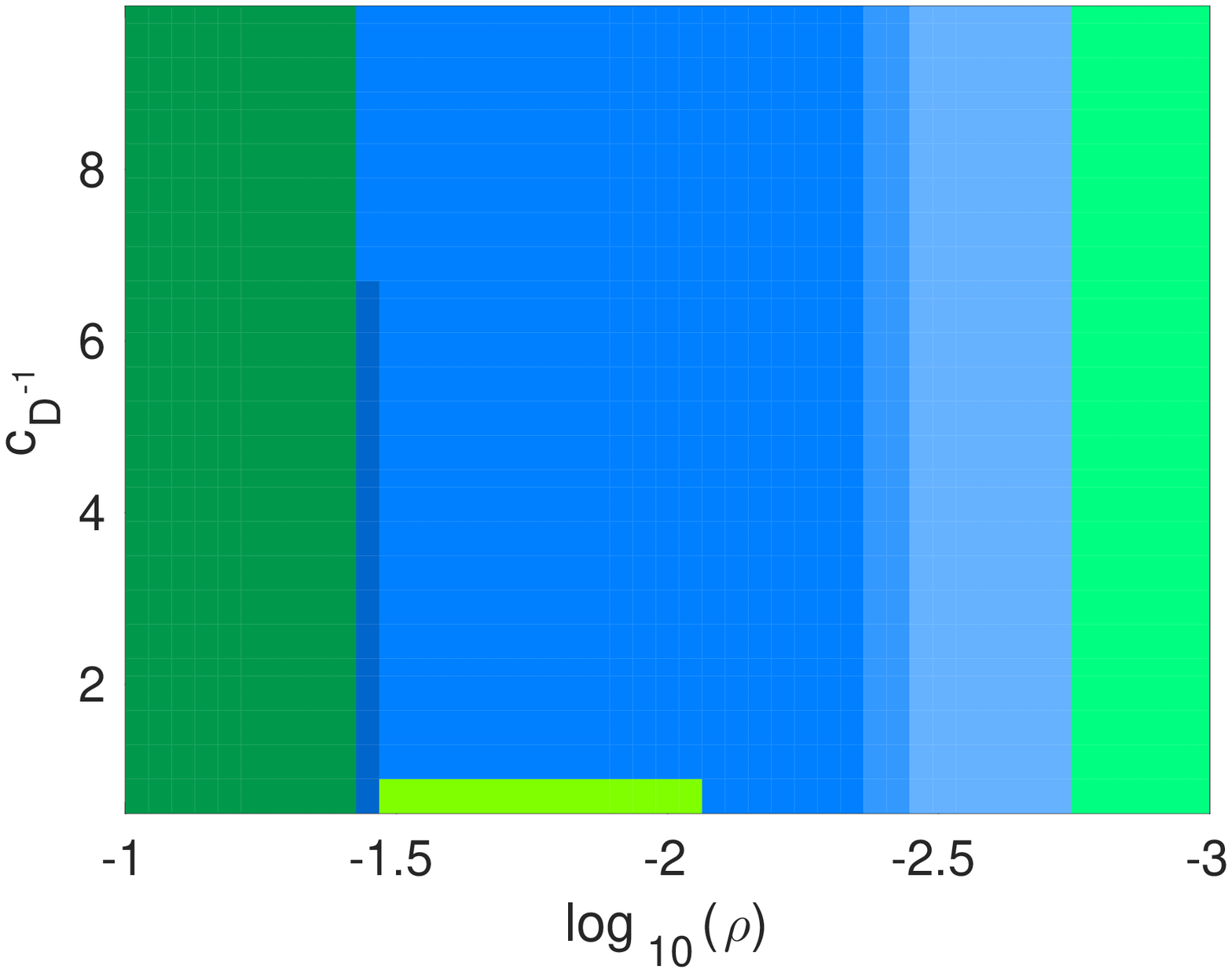} 
        \hspace*{2mm}
        \includegraphics[height=3.4cm,width=3.4cm,clip=true]{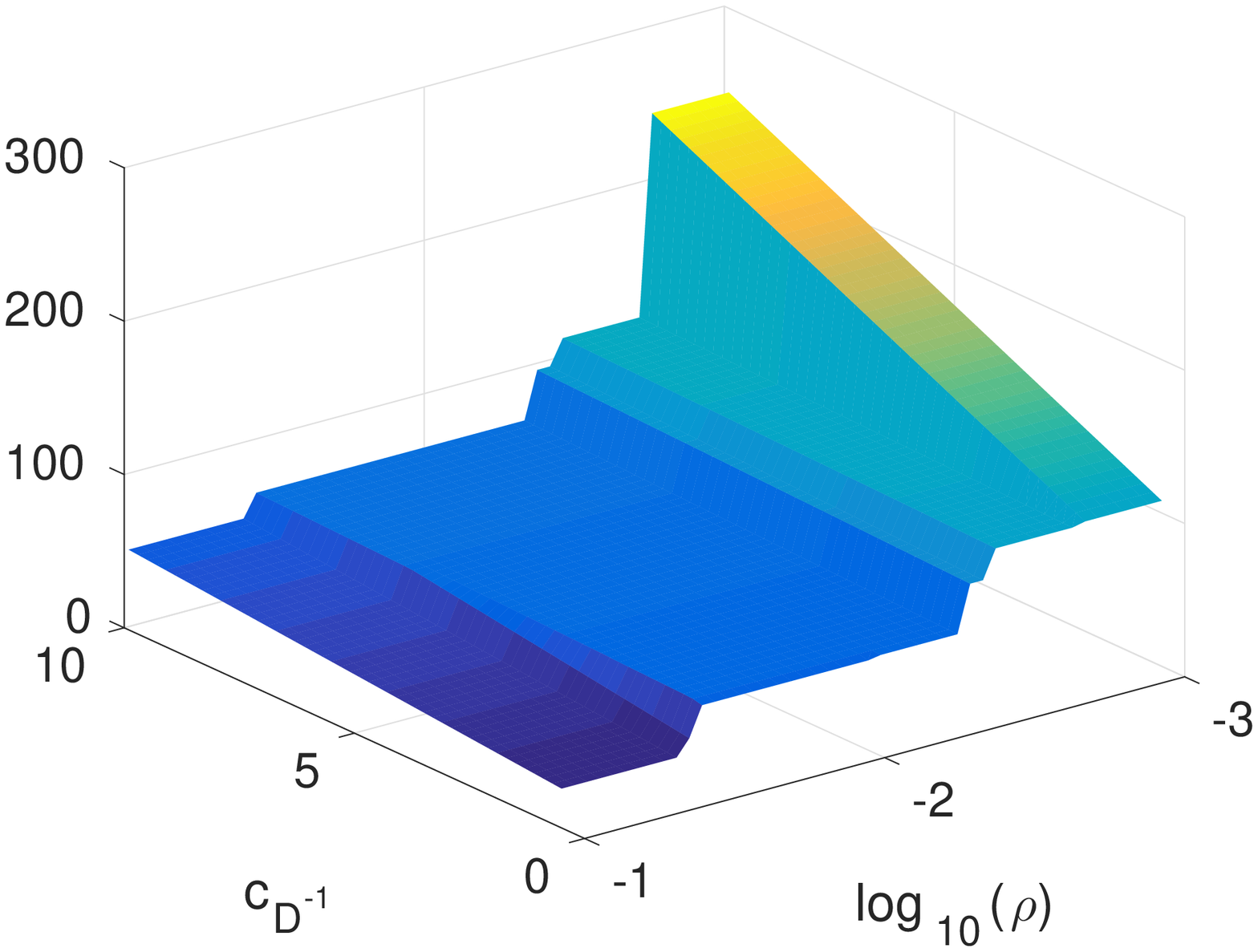}
      }}}
    \vspace*{1mm}
    \centerline{\fbox{\parbox[c]{7.2cm}{
        \centerline{$p=25$}
        \vspace*{2mm}
        \includegraphics[height=3.4cm,width=3.4cm,clip=true]{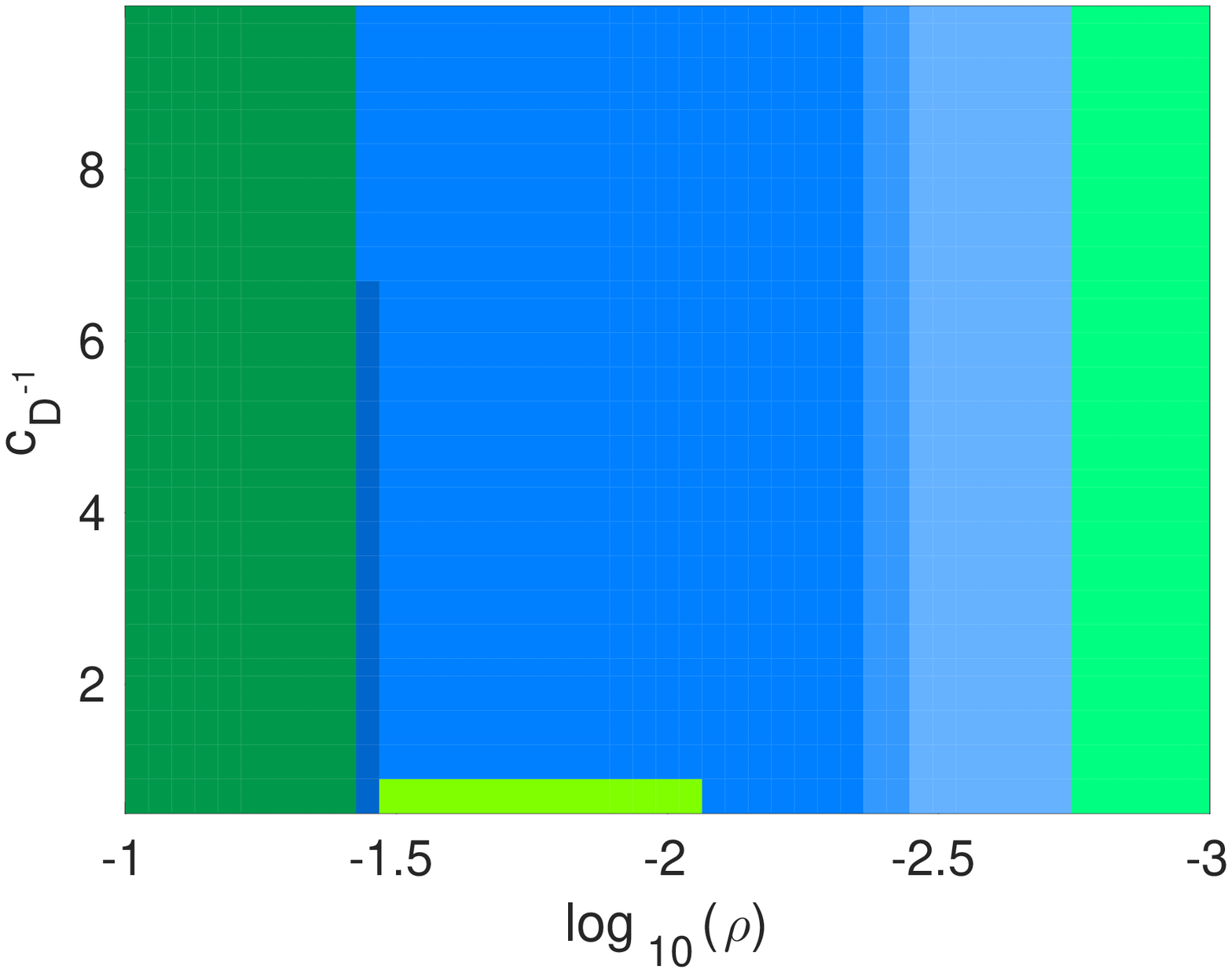} 
        \hspace*{2mm}
        \includegraphics[height=3.4cm,width=3.4cm,clip=true]{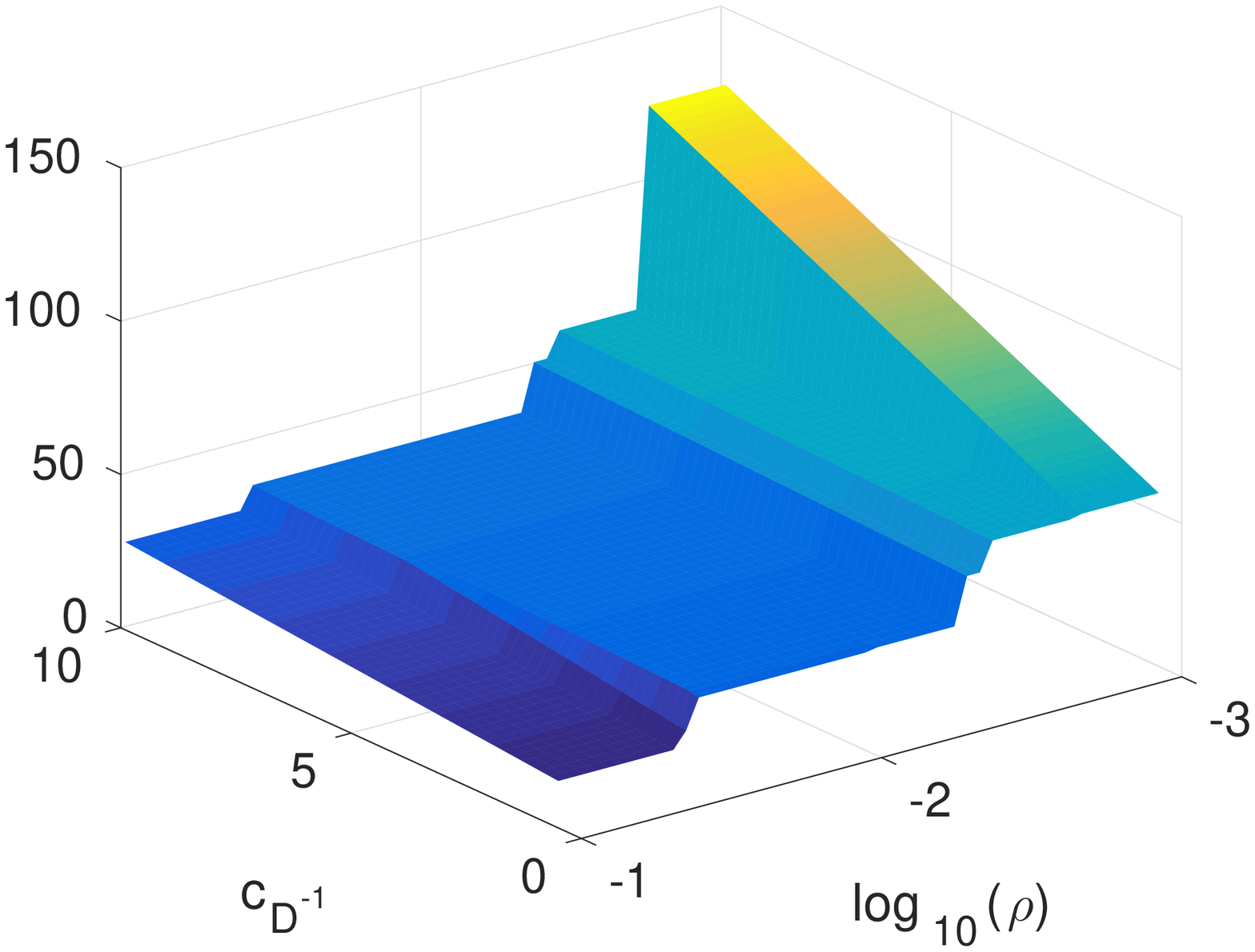}
      }}
    \fbox{\parbox[c]{7.2cm}{
        \centerline{$p=50$}
        \vspace*{2mm}
        \includegraphics[height=3.4cm,width=3.4cm,clip=true]{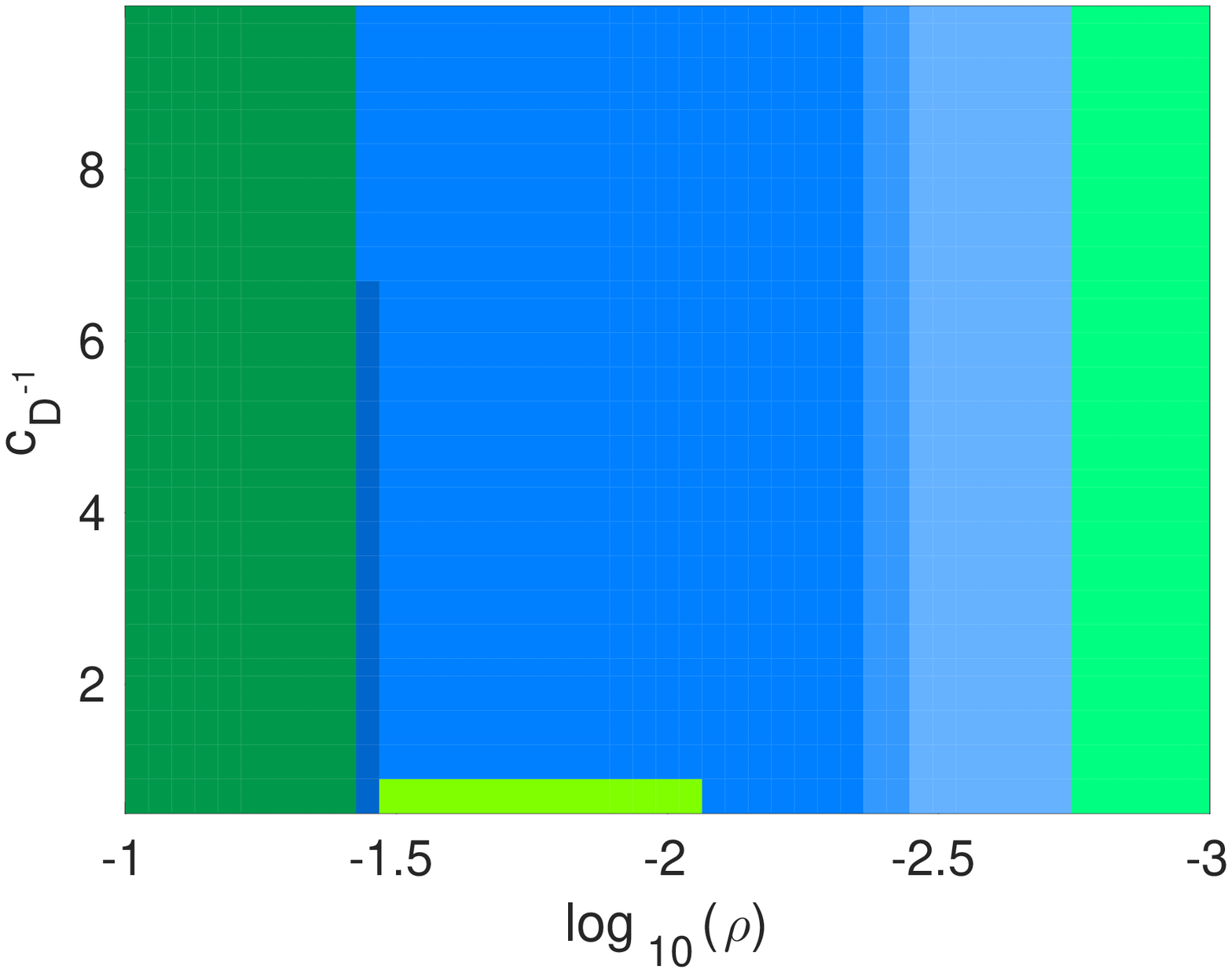} 
        \hspace*{2mm}
        \includegraphics[height=3.4cm,width=3.4cm,clip=true]{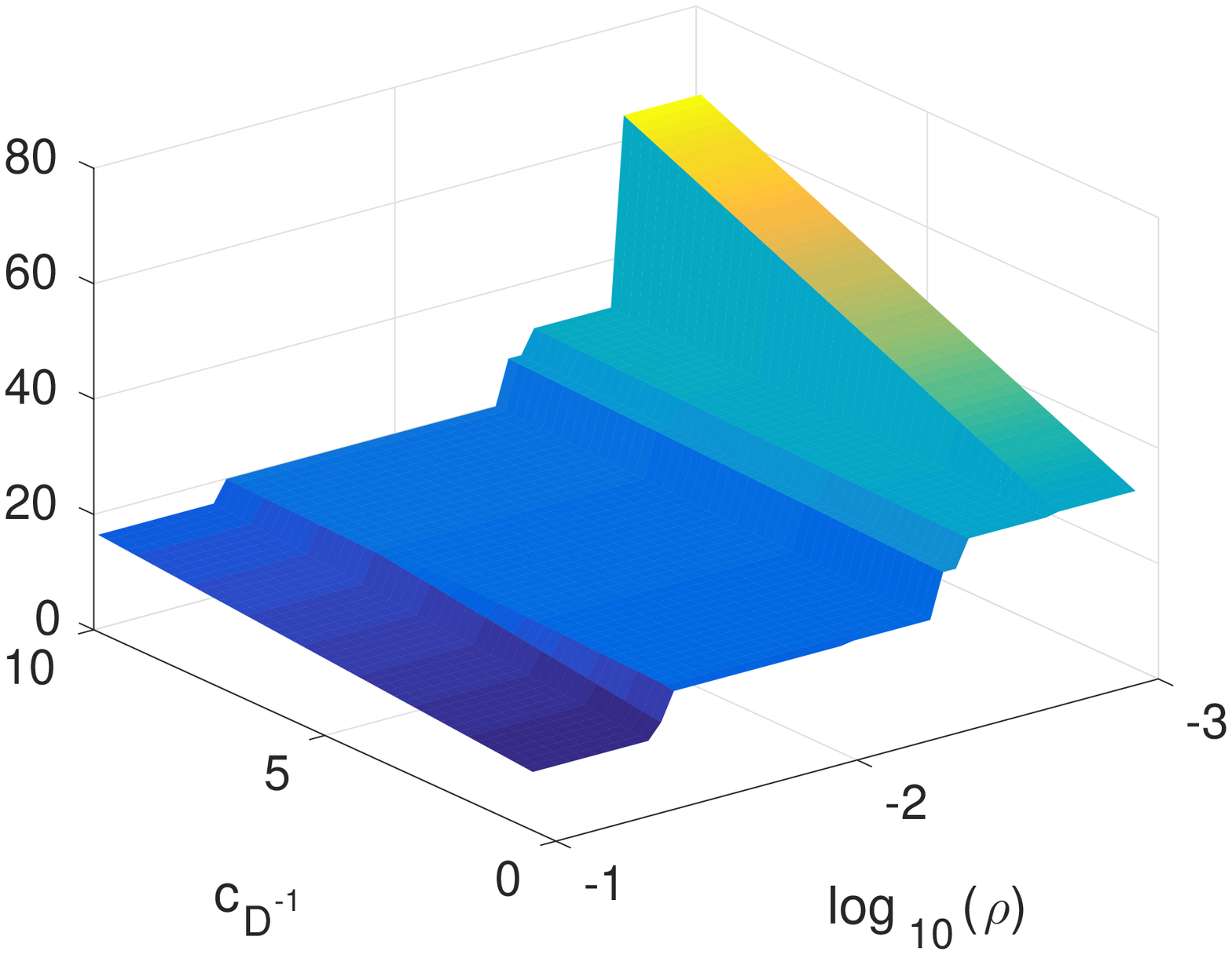}
      }}}
    \vspace*{2mm}
  \vspace*{1mm}
  \begin{center}
   Map colors: 
  \begin{tabular}{llllllll}
  \includegraphics[height=3mm,width=3mm]{sa1q1_Li_color.eps}  \hspace*{1mm} & SAQ1-M-0  &
  \includegraphics[height=3mm,width=3mm]{sa1q15_Li_color.eps} \hspace*{1mm} & SAQ1-M-0  &
  \includegraphics[height=3mm,width=3mm]{sa1q25_Li_color.eps} \hspace*{1mm} & SAQ25-M-0 &
  \includegraphics[height=3mm,width=3mm]{sa1q50_Li_color.eps} \hspace*{1mm} & SAQ40-M-0 \\
  \includegraphics[height=3mm,width=3mm]{sa1q50_Mi_color.eps} \hspace*{1mm} & SAQ50-M-I &
  \includegraphics[height=3mm,width=3mm]{stq1_Mi_color.eps}   \hspace*{1mm} & STQ1-S-I &
  \includegraphics[height=3mm,width=3mm]{stq15_Mi_color.eps}  \hspace*{1mm} & STQ15-S-I &
  \includegraphics[height=3mm,width=3mm]{stq15_Li_color.eps}  \hspace*{1mm} & STQ15-S-0 \\
  \includegraphics[height=3mm,width=3mm]{sfo1q1_ex_color.eps}   \hspace*{1mm} & FOQ1-D  &
  \includegraphics[height=3mm,width=3mm]{sfo1q15_ex_color.eps}  \hspace*{1mm} & FOQ15-D &
  \includegraphics[height=3mm,width=3mm]{sfo1q25_ex_color.eps}  \hspace*{1mm} & FOQ25-D & \\
  \end{tabular}
  \end{center}
  \vspace*{-4mm}
\caption{\label{fig-map-nloc-oops}Best algorithmic variants as a function of
  $c_{D^{-1}}$  the number of computing processes $p$ and the reliability
  factor $\rho$ (QG example, hybrid MPI/OpenMP model)}
\vspace*{-8mm}
\end{figure}
        
\noindent
The situation is much more intricate here than for the Burgers example.
In particular, the competition between state and forcing formulations is very
tight in the sequential case. As above, the latter stops being advantageous
when the number of computing processes increases. The maximal computational
costs now range from approximately 1977 ($p=1$) to 65 ($p=50$), giving
an excellent speed-up of 30.  Looking at the ``best-method'' maps in more
detail, we see that state-based algorithms seems to obtain a better decrease
quickly, then being outperformed by saddle-based methods, for finally
nevertheless taking over (using $\tilde{M}_i = I$) for the more stringent accuracy
requirements. We also note by looking at the minimum-cost surfaces that such
requirements come at a significant computational cost (irrespective of $p$),
in contrast with what was observed for the Burgers example. The original
saddle algorithm, which we kept in the comparison here, is never the best method.

\numsection{Avoiding the use of $D^{-1}$}\label{approxDinv-s}

As already noted, the cost of the $D^{-1}$ operator may vary considerably from
application to application.  In some oceanographic models, $D$ is computed
using a diffusion operator which is integrated using an implicit scheme
\cite{WeavCour01,MiroWeav10,GratToinTshi12,WeavTshiPiac16}. This makes the
cost of $D^{-1}$ very comparable to that of $D$.  In some other applications,
such as atmospheric modelling and weather forecasting, the $D$ operator may
involve more complicated elements, such as localization schemes
\cite{WangHamiWhitBish07,ClaytLoreBark13,Lore17}, which makes applying
$D^{-1}$ potentially more costly. The authors are aware of the strong
reluctance of practitioners in these areas to even provide the $D^{-1}$ operator
at all.  If that is the case, and if one nevertheless desires to enjoy the
security of a solid global convergence theory at an acceptable cost, it is
possible to use an approximate $D^{-1}$ operator, for instance by computing
$y=D^{-1}x$ by approximately solving the linear system $Dy = x$.  Various
iterative methods can be considered for this task, including
conjugate-gradients or FOM, which should already be available in the data
assimilation system. In addition, as in the Burgers and QG examples,
parallelism can often be exploited to make this computation efficient. We have
however met two difficulties when experimenting with the idea.  The first is
that a very inaccurate solution of the linear system may result in an
unsymmetric $D^{-1}$ operator, which then empties the very state formulation
of its meaning (in addition to causing numerical havoc).  The second is that
too inexact solutions may also slow down the convergence of STQ$\ell$
significantly, although the saddle-based algorithms seem more robust. We have
however found that a relatively modest number of conjugate gradient iterations
is very often sufficient to reach this accuracy level. In preliminary tests on
the Burgers example with 25 unpreconditoned CG iterations, the SAQ50-M-I
algorithm turned out to be the best choice for both parallel computing models,
and the parallel computing cost ($p=50$) increased by less than a factor three
compared to using the exact $D^{-1}$ (with $c_{D^{-1}}=c_D \in
[\half,10]$). Very similar conclusions can be reached when applying the same
strategy (with 20 unpreconditioned CG iterations) to the QG example. As number of
outer and inner iterations differ only very marginally from that observed when
using the exact $D^{-1}$ and since the cost of 20 products with $D$ brings
$c_{D^{-1}}=10$ to the top of the range considered for exact $D^{-1}$, the
results are essentially identical to those obtained for $c_{D^{-1}}=10$ in
Section~\ref{numerics-qg-s}.

There is little doubt that preconditioning and problem specific tuning would
reduce the cost of the approximate $D^{-1}$ even further.  Our experiments
therefore show that, even if the $D^{-1}$ operator is unavailable, considering
the SAQ$\ell$ or STQ$\ell$ algorithms may be the best option, irrespective of
the number of computing processes and of the parallel computing model.

We finally note that using the forcing formulation makes it possible to avoid
using $B^{-1}$ (not $D^{-1}$) altogether by starting the Gauss-Newton
algorithm with $x_0= x_b$ and recurring $B^{-1}(x_k-x_b)$ over successive
major iterations from by-products of the FOM or CG algorithms. This technique
however suffers from the same parallelization problems as the standard FOM and
does not avoid using the (possibly approximate) $Q_i^{-1}$ operators (see
\req{H-D-def}). This is why we haven't considered this variant in detail in
our parallel computing assessment.

\numsection{Conclusion and perspectives}\label{concl-s}

In this paper, we have exposed the problematic behaviour of the original
saddle formulation as a general method for solving the weakly constrained 4D-Var
problem.  Its undesirable features are caused by the very poor correlation,
for approximate solutions, between quadratic model decrease (the objective)
and reduction of the residual of the associated optimality conditions (the
mean).  This mismatch in turn causes the values of the cost function(s) to behave
chaotically and makes terminating the inner iteration too much dependent on
chance, potentially resulting in divergence of the whole process. We have
nevertheless proposed a strategy (and a corresponding class of algorithms)
which cures the problem and for which strong global convergence results can be
proved.

We have then experimented with this new class of saddle-based algorithms and
compared their performance with that of methods associated with alternative
variational formulations of the problem. This comparison was conducted on two
different and complementary examples of data assimilation, taking into account
not only performance in terms of number of iterations, but also considering
two more elaborate approximations of computational costs in a two different
parallel computing models. A parametric study of the sequential and parallel
computing cost as a function of the costs of applying the $D^{-1}$ operator
and the accuracy obtained has been conducted, showing the relative merits of
the new saddle algorithms and the more classical CG/FOM solvers for the state
formulation. Both appear to have their place in the data assimilation
toolbox. We have also provided a preliminary discussion of the application of
both classes of algorithms in the case where the $D^{-1}$ operator is
unavailable, indicating that similar conclusions hold if it is approximated.

Several issues remain to be explored further, one of which is the use of
approximate operators: we only briefly touched the question in our discussion
of the use of the approximate $D^{-1}$, and further elaboration
including preconditioning and the possible use of inexact products
\cite{GratGuroToinTshiWeav12} might be of interest. The second and most
important one is the translation of our conclusions, drawn in a relatively
controlled context, to the more complex environments of truly parallel operational
systems.


\appendix

\appnumsection{A1. The Burgers assimilation problem}\label{burger-a}

We consider the one-dimensional Burgers equation on the spatio-temporal domain
$\Omega=[0,T]\times[0,1]$ whose governing equation is
\begin{equation}
\displaystyle{ \frac{\partial u}{\partial t}+ u \frac{\partial u}{\partial
    x}-\nu \frac{\partial^2 u}{\partial x^2}=g}
\mbox{ on } \Omega\\
\end{equation}
with Dirichlet boundary conditions 
\begin{equation}
u(0,t)=u(1,t)=0 \quad \forall t \geq 0. \\
\end{equation}
The forcing term is then given by
\begin{equation}
\begin{array}{rcl}
g(x,t)&=&\pi k [x+ k(t+1)\sin \pi (1-x)(t+1)] \cos \pi x(t+1) \sin \pi (1-x)(t+1)\\
& &+\pi k [1-x - k(t+1)\sin \pi x(t+1)] \sin \pi x(t+1) \cos \pi(1-x)(t+1)\\
& &+2\nu k^2 \pi^2(t+1)^2 [\sin \pi x(t+1) \sin \pi(1-x)(t+1)\\
& &\hspace*{5.2cm} +\cos \pi x(t+1) \cos \pi(1-x)(t+1)]
\end{array}
\end{equation}
The discretization uses a first-order upwind scheme in time and a
second-order centered scheme in space. With a space step $\Delta x$ and a
time step $\Delta t$, this gives
\begin{equation}
\frac{1}{\Delta t}(u_i^{n+1}-u_i^n)+\frac{u_i^n}{2 \Delta
  x}(u_{i+1}^n-u_{i-1}^n)-\frac{\nu}
     {(\Delta x)^2}(u_{i+1}^n-2u_i^n+u_{i-1}^n)=g(i\Delta x,n \Delta t)
\end{equation}
with $u_i^n=u(i\Delta x, n \Delta t)$. We choose $\Delta x=0.01$ leading to a
state vector of dimension $n=100$, $\Delta t= 10^{-5}$  and a diffusion
coefficient $\displaystyle{\nu=0.25}$. The length of the assimilation window
$T$ is equal to $0.03$. It is divided into $N_{sw}=50$ subwindows of equal length. 

The reference solution is built by running the model with the initial condition
\begin{equation}
u_{\rm true}(x,0)=k \sin(2\pi x), \quad \forall x \in [0,1]
\label{eq:true_init}
\end{equation}
with $k=0.1$ as for the forcing  function $g$, and by adding a Gaussian random
variable at the end of each subwindow:
\[
x^t_j=\mathcal{M}_j(x^t_{j-1})+\epsilon^m_j, \quad
\epsilon^m_j\sim\mathcal{N}(0,\sigma_m^2I_n)
\ms \ms (j=1,\ldots,N_{sw})
\]
with $x^t_j=[u_1^j, \cdots, u_n^j]^T$ the state vector at time $t_j$,
and $\mathcal{M}_j$ the integration of the numerical model from time $t_{j-1}$
to $t_j$. We choose $\sigma_m^2=1.10^{-4}T/N_{sw}$. 

At the end of each subwindow (time $t_j$), $m_j=20$ observations are built
by randomly selecting  components of the reference solution and then
adding a Gaussian random variable:
\[
\forall j=1:N_{sw} \quad y_j=\mathcal{H}_j(x^t_{j})+\epsilon^o_j, \quad
\epsilon^o_j\sim\mathcal{N}(0,\sigma_o^2I_{m_j})
\]
with $\mathcal{H}_j$ the observation operator at time $t_j$ (basically the
random selection of $m_j$ components of the reference solution). We choose
$\sigma_o^2=10^{-3}$. This strategy results in a total of $1000$ assimilated
observations for the whole assimilation window. The reference solution and
the observations  at the end of the first and last subwindows are shown in Figure \ref{fig:obs}. 

\begin{figure}[t]
\begin{center}
\hspace{-2cm}\begin{minipage}[c]{7cm}
\begin{center}
\includegraphics[width=8cm,keepaspectratio=true]{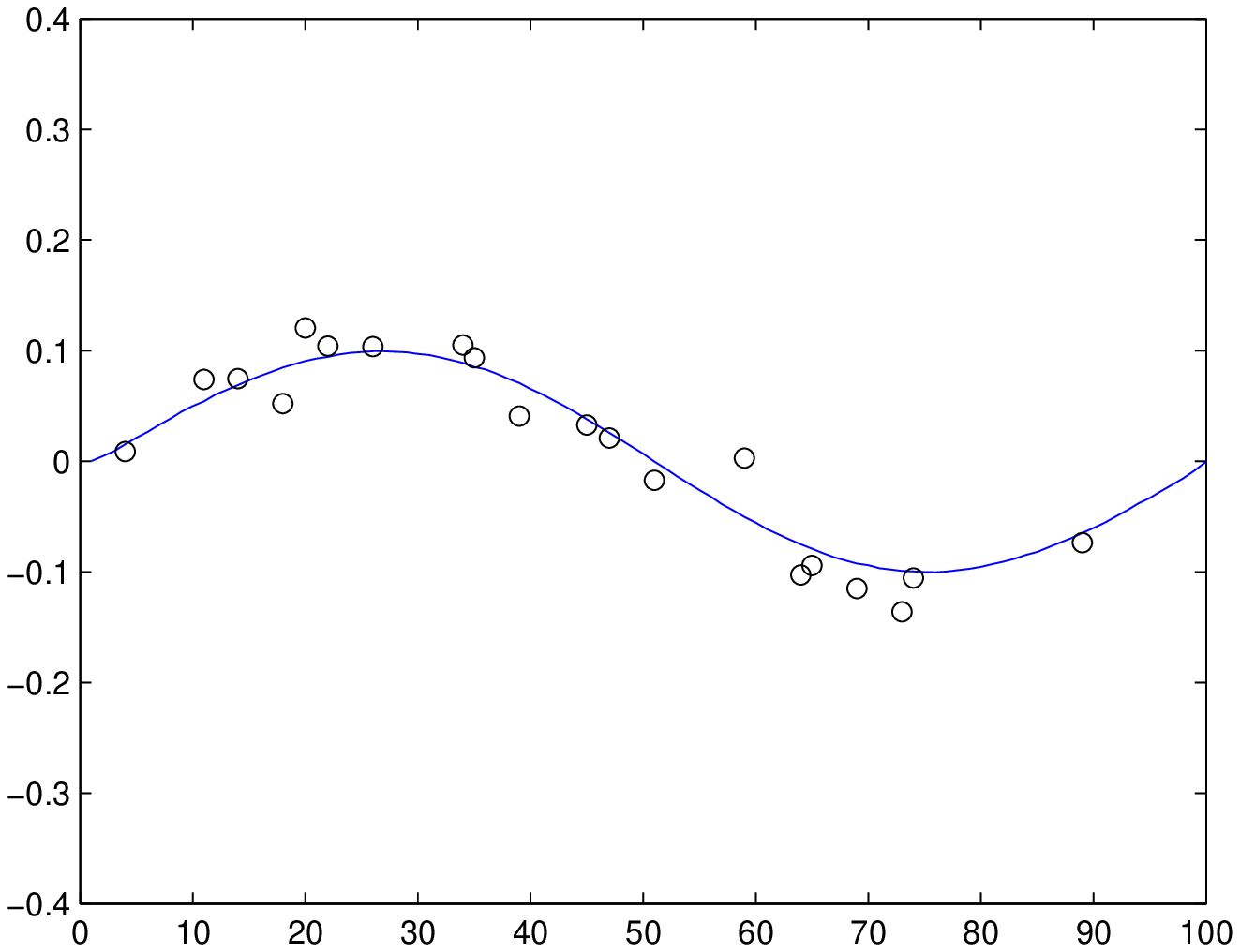}
\end{center}
\end{minipage}
\hspace*{2mm}
\begin{minipage}[c]{7.0cm}
\begin{center}
\includegraphics[width=8cm,keepaspectratio=true]{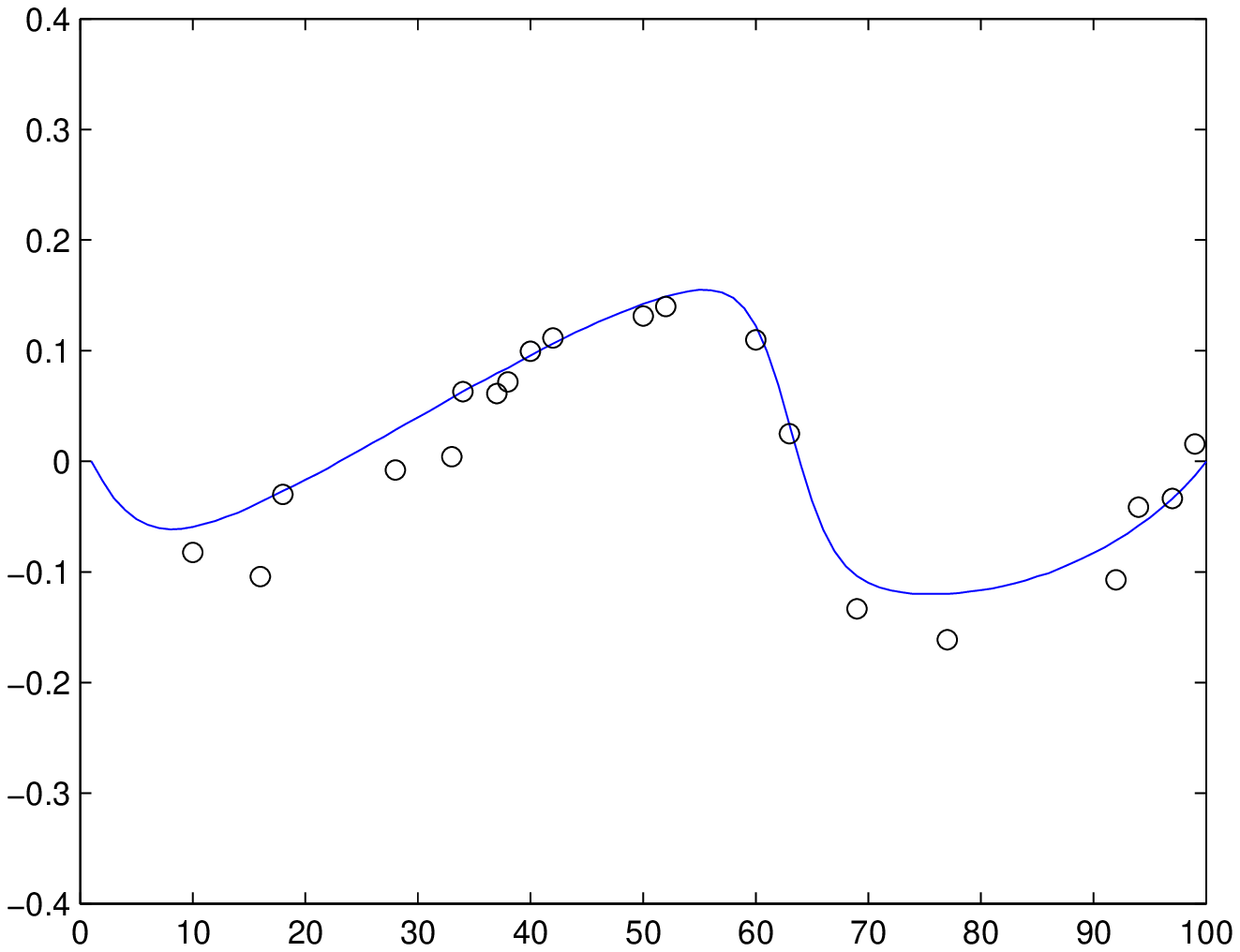}
\end{center}
\end{minipage}
\end{center}
\label{fig:obs}
\caption{Reference trajectory and observations at the end of the first
  (left) and last (right) subwindow}
\end{figure}

The observation error covariance matrices $R_j$, with $j=1:N_{sw}$,
introduced in the definition of the 4D-Var cost function are diagonal and
their diagonal entries are chosen such that they are positive, the largest one
is equal to one and the condition number of $R_j$ is equal to $10^3$.

The background solution $x_b$ corresponds to the sum of the reference solution
at initial time \req{eq:true_init} and a random variable $\epsilon^b\sim
\mathcal{N}(0,\sigma_b^2 I_n)$, with $\sigma_b^2=10^{-2}$. A model error is
introduced during the numerical integration of the model by adding a Gaussian
random variable at the end of each subwindow as for the reference
solution. The background error covariance matrix corresponds to the weighted
sum of the squared exponential covariance and the identity matrix given by
\begin{equation}
\label{eq:cov_mat}
\displaystyle{B=\sigma_b^2(\alpha I_n+(1-\alpha) \tilde{B}),
\quad \mbox{ with }
\tilde{B}_{i,j}= e^{-\frac{d(i,j)^2}{L^2}}}
\end{equation}
with $d(i,j)$ the distance between the spatial grid points $i$ and $j$,
$L=0.25$ a specified length scale, and $\alpha\in[0,1]$ the weight associated
with both matrices. This last parameter is a simple way to allow variation in
the condition number of the matrix $B$. We choose $\alpha=0.001$ which
results in a condition number of $1.10^5$. The model error covariance matrices
$Q_j$, with $ j=1:N_{sw}$, are built using the same strategy except that
$L=0.05$ and $\alpha=0.01$. The condition number of these matrices is close to
$1.65$ $10^3$.

\appnumsection{A2. The ECMWF QG problem}\label{burger-b}

In the quasi Geostrophic (QG) problem, a reference stream function is
generated from a~model with layer depths of $H_1=6000\,m$ and $H_2=4000\,m$,
and the time step is set to $600\,s$, whereas the assimilating model has layer
depths of $H_1=5500\,m$ and $H_2=4500\,m$, and the time step is set to
$3600\,s$. These differences in the layer depths and the time steps provide
a~source of model error.

Observations of the non-dimensional stream function, vector wind and wind
speed were taken from the reference of the model at $100$ points randomly
distributed over both levels for each hour.  Observation errors were assumed
to be independent from each others and uncorrelated in time, the standard
deviations were chosen to be $0.4$ for the stream function observation error,
$0.6$ for the vector wind and $1.2$ for the wind speed. The observation
operator is the bi-linear interpolation of the model fields to horizontal
observation locations.

The background error covariance matrix ($B$ matrix) and the model error
covariances (matrices $Q_{i}$) correspond to vertical and horizontal
correlations. The vertical and horizontal structures are assumed to be
separable.  In the horizontal plane, covariance matrices correspond to
isotropic, homogeneous correlations of stream function with Gaussian spatial
structure.  For the background error covariance matrix $B$, the standard
deviation and the horizontal correlation length scale are set to $0.8$ and
$10^6\,m$ respectively. For the model error covariance matrices $Q_i$, the
standard deviation and the horizontal correlation length scale are set to
$0.6$ and $2\times10^5\,m$ respectively.  The vertical correlation is assumed
to be constant over the horizontal grid and the correlation coefficient value
between the two layers was taken as $0.2$ for $B$ and $0.5$ for $Q_i$.

The length of the assimilation window is set to $48$ hours, divided into $48$
equal sub-windows of $1$ hour each.

\appnumsection{A3. The left-preconditioned FOM algorithm and its application to the
  forcing formulation}\label{burger-c}

We first state, as Algorithm~\ref{foms} \vpageref{foms} the
left-preconditioned FOM algorithm for general symmetric positive definite
system $Ax=b$ with symmetric positive definite preconditioning matrix
$M$. This algorithm uses the inner product induced by $M^{-1}$. In the
description, we use the notation $1\!:\!k$ as a short-hand for $\ii{k}$.

\algo{foms}{Left-preconditioned FOM algorithm for solving $MAx=Mb$}{
  \begin{description}
  \item[1. Initialization.] Symmetric positive definite matrices $A,M \in
    \Re^{n\times n}$ are given, as well as a right-hand side $b\in \Re^n$.
    \begin{enumerate}
    \item[1.1 ] $w \leftarrow Mb$
    \item[1.2 ] $\beta \leftarrow \sqrt{ b^Tw }$
    \item[1.3 ] $U_{1:n,1} \leftarrow w/\beta$
    \item[1.4 ] $Q_{1:n,1} \leftarrow b/\beta$
    \item[1.5 ] $z_1 \leftarrow U_{1:n,1}^Tb$
    \end{enumerate}
  \item[2. Main loop.] For $i = 1, \ldots, {\rm maxit}$,
    \begin{enumerate}
    \item[2.1 ] $w \leftarrow AU_{1:n,k}$
    \item[2.2 ] $v \leftarrow Mw$
    \item[2.3 ] for $j = 1,\ldots, k$
    \begin{enumerate}
    \item[2.3.1 ] $H_{j,k} \leftarrow Q_{1:n,j}^Tv$
    \item[2.3.2 ] $v \leftarrow v - U_{1:n,j}H_{j,k}$
    \item[2.3.3 ] $w \leftarrow w - Q_{1:n,j}H_{j,k}$
    \end{enumerate}
    \item[2.4 ] $H_{k+1,k} \leftarrow \sqrt{w^Tv}$
    \item[2.5 ] $y \leftarrow \beta H_{1:k,1:k}^{-1} e_1$
    \item[2.6 ] $\gamma_k \leftarrow | H_{k+1,k}y_k|$
    \item[2.7 ] $q_k  \leftarrow - \half z^Ty$
    \item[2.8 ] In view of $q_k$ and $\gamma$, terminate with $x = Uy$?
    \item[2.9 ] $U_{1:n,k+1}  \leftarrow v / H_{k+1,k}$
    \item[2.10 ] $Q_{1:n,k+1} \leftarrow w / H_{k+1,k}$
    \item[2.11 ] $z_{k+1} \leftarrow U_{1:n,k+1}^Tb$
    \end{enumerate}     
  \end{description}
}

Note that Steps 1.5, 2.7 and 2.11 are only necessary if the value of the model
quadratic
$
q(x) = \half x^TAx - b^Tx
$
must be tracked in the course of the inner iterations.  If this is the case
$q_k$ is the value of $q(x)= q(Uy)$ at iteration $k$, while $\gamma_k =
\|\nabla_x q(x)\|_M$ (the preconditioned norm of the system's residual).  These values may then
be used to decide on termination in Step~2.8 (for instance according to
\req{qm-termination}). Also note that the operator $M^{-1}$ is only used
implictly and never appears in the algorithm.

We next consider applying this algorithm (with $M=D$) to find $\delta x = L^{-1} \delta p$,
where $\delta p$ (approximately) solves \req{forcing-formulation}.  This gives
Algorithm~\ref{foms2} \vpageref{foms2}. Storing the matrix $P$ in Step~2.1 of
this algorithm allows avoiding the backsolve \req{fo-backsolve} by returning
$\delta x = Py$ instead of $\delta p = Uy$.

\algo{foms2}{Specialized FOM algorithm for \req{forcing-formulation}-\req{fo-backsolve}}{
  \begin{description}
  \item[1. Initialization.] The matrices $L, D, H, R \in
    \Re^{n\times n}$ are given, as well as a right-hand side $r = D^{-1}b + L^{-T}H^{T}R^{-1}d \in \Re^n$.
    \begin{enumerate}
    \item[1.1 ] $w \leftarrow Dr$
    \item[1.2 ] $\beta \leftarrow \sqrt{ w^Tr }$
    \item[1.3 ] $U_{1:n,1} \leftarrow w/\beta$
    \item[1.4 ] $Q_{1:n,1} \leftarrow r/\beta$
    \item[1.5 ] $z_1 \leftarrow U_{1:n,1}^Tr$
    \end{enumerate}
  \item[2. Main loop.] For $i = 1, \ldots, {\rm maxit}$,
    \begin{enumerate}
    \item[2.1 ] $P_{1:n,k} \leftarrow L^{-1}U_{1:n,k}$
    \item[2.2 ] $v \,\leftarrow L^{-T}H^{T}R^{-1}HP_{1:n,k}$
    \item[2.3 ] $w \leftarrow Q_{1:n,k}+ v$
    \item[2.4 ]  $v \,\leftarrow U_{1:n,k}+ Dv$
    \item[2.5 ]  for $j = 1,\ldots, k$
    \begin{enumerate}
    \item[2.6.1 ] $T_{j,k} \leftarrow Q_{1:n,j}^Tv$
    \item[2.6.2 ] $v \,\leftarrow \,v - U_{1:n,j}T_{j,k}$
    \item[2.6.3 ] $w \leftarrow w - Q_{1:n,j}T_{j,k}$
    \end{enumerate}
    \item[2.7 ] $T_{k+1,k} \leftarrow \sqrt{w^Tv}$
    \item[2.8 ] $y \leftarrow \beta T_{1:k,1:k}^{-1} e_1$
    \item[2.9 ] $\gamma_k \leftarrow | T_{k+1,k}y_k|$
    \item[2.10 ] $q_k  \leftarrow - \half z^Ty$
    \item[2.11 ] In view of $q_k$ and $\gamma$, terminate with $\delta x = Py$?
    \item[2.12 ] $U_{1:n,k+1}  \leftarrow v / T_{k+1,k}$
    \item[2.13 ] $Q_{1:n,k+1} \leftarrow w / T_{k+1,k}$
    \item[2.14 ] $z_{k+1} \leftarrow  U_{1:n,k+1}^Tr$
    \end{enumerate}     
  \end{description}
}


{\footnotesize
\begin{thebibliography}{10}

\bibitem{BenzWath08}
M.~Benzi and A.~Wathen.
\newblock Some preconditioning techniques for saddle point problems.
\newblock In {\em Model Order Reduction}, number~13 in Mathematics in Industry,
  pages 195--211, Heidelberg, Berlin, New York, 2008. Springer Verlag.

\bibitem{BergGondVentZill07}
L.~Bergamaschi, J.~Gondzio, M.~Venturin, and G.~Zilli.
\newblock Inexact constraint preconditioners for linear systems arising in
  interior point methods.
\newblock {\em Computational Optimization and Applications}, 36(2-3):136--147,
  2007.

\bibitem{BergGondVentZill11}
L.~Bergamaschi, J.~Gondzio, M.~Venturin, and G.~Zilli.
\newblock Erratum to: {I}nexact constraint preconditioners for linear systems
  arising in interior point methods.
\newblock {\em Computational Optimization and Applications}, 49(2):401--406,
  2011.

\bibitem{Bout99}
F.~Bouttier and P.~Courtier.
\newblock Data assimilation concepts and methods.
\newblock Technical report, ECMWF, Reading, England, 1999.
\newblock ECMWF Meteorological Training Course Lecture Series.

\bibitem{ClaytLoreBark13}
A.~M. Clayton, A.~C. Lorenc, and D.~M. Barker.
\newblock Operational implementation of a hybrid ensemble/{4D-Var} global data
  assimilation system at the {M}et {O}ffice.
\newblock {\em Quarterly Journal of the Royal Meteorological Society},
  139:1445--1461, 2013.

\bibitem{ConnGoulToin00}
A.~R. Conn, N.~I.~M. Gould, and Ph.~L. Toint.
\newblock {\em Trust-Region Methods}.
\newblock MPS-SIAM Series on Optimization. SIAM, Philadelphia, USA, 2000.

\bibitem{Cour97}
Ph. Courtier.
\newblock Dual formulation of four-dimensional variational assimilation.
\newblock {\em Quarterly Journal of the Royal Meteorological Society},
  123:2449--2461, 1997.

\bibitem{CourThepHoll94}
Ph. Courtier, J.-N. Th\'{e}paut, and A.~Hollingsworth.
\newblock A strategy for operational implementation of {4D-Var} using an
  incremental approach.
\newblock {\em Quarterly Journal of the Royal Meteorological Society},
  120:1367--1388, 1994.

\bibitem{ElSa15}
A.~El{-}Said.
\newblock {\em Variational Data Assimilation Problem for Numerical Weather
  Prediction}.
\newblock PhD thesis, University of Reading, Reading, UK, 2015.

\bibitem{ElSaNichLawl17}
A.~El{-}Said, N.~K. Nichols, and A.~S. Lawless.
\newblock Conditioning of the weak-constraint {4DVAR} problem.
\newblock Technical report, University of Reading, Reading, UK, 2017.

\bibitem{FishGratGuroVassTrem17}
M.~Fisher, S.~Gratton, S.~G\"{u}rol, Y.~Tr\'{e}molet, and X.~Vasseur.
\newblock Low rank updates in preconditioning the saddle point systems arising
  from data assimilation problems.
\newblock {\em Optimization Methods and Software}, \textmd{(to appear)}, 2017.

\bibitem{FishGuro17}
M.~Fisher and S.~G\"{u}rol.
\newblock Parallelisation in the time dimension of four-dimensional variational
  data assimilation.
\newblock {\em Quarterly Journal of the Royal Meteorological Society},
  143(703):1136--1147, 2017.

\bibitem{FishTremAuviTanPoli11}
M.~Fisher, Y.~Tr{\'e}molet, H.~Auvinen, D.~Tan, and P.~Poli.
\newblock Weak-constrained and long window {4D-V}ar.
\newblock Technical Report 655, ECMWF, 2011.

\bibitem{FreiGree17}
M.~A. Freitag and D.~L.~H. Green.
\newblock A low-rank approach to the solution of weak constraint variational
  data assimilation problems.
\newblock arXiv:1702.07278v1, 2017.

\bibitem{GratGuroSimoToin17}
S.~Gratton, S.~G\"{u}rol, E.~Simon, and Ph.~L. Toint.
\newblock Issues in making the weakly-constrained {4DVar} formulation
  computationally efficient.
\newblock Oberwolfach Reports 47, 2017.

\bibitem{GratGuroSimoToin17c}
S.~Gratton, S.~G\"{u}rol, E.~Simon, and Ph.~L. Toint.
\newblock Preconditioning weighted linear least-squares with an application to
  weakly constrained variational data assimilation.
\newblock (in preparation), 2017.

\bibitem{GratGuroToinTshiWeav12}
S.~Gratton, S.~G\"{u}rol, Ph.~L. Toint, J.~Tshimanga, and A.~Weaver.
\newblock Krylov methods in the observation space for data assimilation.
\newblock Oberwolfach Reports, 2012.

\bibitem{GratLawlNich07}
S.~Gratton, A.~Lawless, and N.~K. Nichols.
\newblock Approximate {G}auss-{N}ewton methods for nonlinear least-squares
  problems.
\newblock {\em SIAM Journal on Optimization}, 18:106--132, 2007.

\bibitem{GratToinTshi11}
S.~Gratton, Ph.~L. Toint, and J.~Tshimanga.
\newblock Range-space variants and inexact matrix-vector products in {K}rylov
  solvers for linear systems arising from inverse problems.
\newblock {\em SIAM Journal on Matrix Analysis}, 32(3):969--986, 2011.

\bibitem{GratToinTshi12}
S.~Gratton, Ph.~L. Toint, and J.~Tshimanga.
\newblock Conjugate-gradients versus multigrid solvers for diffusion-based
  correlation models in data assimilation.
\newblock {\em Quarterly Journal of the Royal Meteorological Society},
  139:1481--1487, 2013.

\bibitem{GratTshi09}
S.~Gratton and J.~Tshimanga.
\newblock An observation-space formulation of variational assimilation using a
  restricted preconditioned conjugate-gradient algorithm.
\newblock {\em Quarterly Journal of the Royal Meteorological Society},
  135:1573--1585, 2009.

\bibitem{HestStie52}
M.~R. Hestenes and E.~Stiefel.
\newblock Methods of conjugate gradients for solving linear systems.
\newblock {\em Journal of the National Bureau of Standards}, 49:409--436, 1952.

\bibitem{LeDiTala86}
F.-X. {Le Dimet} and O.~Talagrand.
\newblock Variational assimilation of meteorological observations: Theoretical
  aspect.
\newblock {\em Tellus}, 38A:97--110, 1986.

\bibitem{Lore17}
A.~C. Lorenc.
\newblock Improving ensemble covariances in hybrid variational data
  assimilation without increasing ensemble size.
\newblock {\em Quarterly Journal of the Royal Meteorological Society},
  143:1062--1072, 2017.

\bibitem{MiroWeav10}
I.~Mirouze and A.~T. Weaver.
\newblock Representation of correlation functions in variational assimilation
  using an implicit diffusion operator.
\newblock {\em Quarterly Journal of the Royal Meteorological Society},
  136:1421--1443, 2010.

\bibitem{NoceWrig99}
J.~Nocedal and S.~J. Wright.
\newblock {\em Numerical Optimization}.
\newblock Series in Operations Research. Springer Verlag, Heidelberg, Berlin,
  New York, 1999.

\bibitem{Saad96}
Y.~Saad.
\newblock {\em Iterative Methods for Sparse Linear Systems}.
\newblock PWS Publishing Company, Boston, USA, 1996.

\bibitem{SaadSchu86}
Y.~Saad and M.~Schultz.
\newblock {GMRES}: A generalized minimal residual algorithm for solving
  nonsymmetric linear systems.
\newblock {\em SIAM J. Sci. Stat. Comput.}, 7:856--869, 1986.

\bibitem{Trem06}
Y.~Tr\'{e}molet.
\newblock Accounting for an imperfect model in {4D-Var}.
\newblock {\em Quarterly Journal of the Royal Meteorological Society},
  132(621):2483--2504, 2006.

\bibitem{Trem07}
Y.~Tr\'{e}molet.
\newblock Model error estimation in {4D-Var}.
\newblock {\em Quarterly Journal of the Royal Meteorological Society},
  133(626):1267--1280, 2006.

\bibitem{VidaPiacleDi04}
P.~A. Vidard, A.~Piacentini, and F.-X. {Le Dimet}.
\newblock Variational data analysis with control of the forecast bias.
\newblock {\em Tellus}, 56A:177--188, 2004.

\bibitem{WangHamiWhitBish07}
X.~Wang, T.~M. Hamill, J.~S. Whitaker, and C.~H. Bishop.
\newblock A comparison of hybrid ensemble transform kalman filter - optimum
  interpolation and ensemble square root filter analysis scheme.
\newblock {\em Monthly Weather Review}, 135:1055--1076, 2007.

\bibitem{Wath15}
A.~Wathen.
\newblock Preconditioning.
\newblock {\em Acta Numerica}, 24:329--376, 2015.

\bibitem{WeavCour01}
A.~T. Weaver and {Ph}. Courtier.
\newblock Correlation modelling on the sphere using a generalized diffusion
  equation.
\newblock {\em Quarterly Journal of the Royal Meteorological Society},
  127:1815--1846, 2001.

\bibitem{WeavTshiPiac16}
A.~T. Weaver, J.~Tshimanga, and A.~Piacentini.
\newblock Correlation operators based on implicitly formulated diffusion
  equation solved with the {C}hebyshev iteration.
\newblock {\em Quarterly Journal of the Royal Meteorological Society},
  142:455--471, 2016.

\bibitem{Zupa97}
D.~Zupanski.
\newblock A general weak constraint applicable to operational {4DVAR} data
  assimilation systems.
\newblock {\em Monthly Weather Review}, 125, 1997.

\end{thebibliography}
}
\end{document}